\documentclass{amsart}
\usepackage{amssymb,amscd,verbatim,amsthm,amsmath,amsfonts,amsgen,xspace}
\usepackage{graphicx, color}

\usepackage{epic}
\usepackage[vcentermath]{youngtab}
\usepackage[english]{babel}
\usepackage{amsmath}
\usepackage{amsthm}
\usepackage{amssymb}
\usepackage[all]{xy}
\usepackage{xypic}

\usepackage{latexsym}
\usepackage{exscale}
\usepackage{comment}

\usepackage{amssymb}
\usepackage{amscd}
 
 \usepackage{amssymb,amscd,verbatim, amsthm,amsmath,amsgen,xspace}
\usepackage{latexsym}
\usepackage{amsfonts}
\usepackage{color}
\usepackage{ulem}
\usepackage{float}
\usepackage[all]{xy}

\usepackage{graphicx}

\usepackage[all]{xy}

\usepackage{amsbsy}

\let\ssize\scriptstyle
\newif\ifFIRST\newdimen\MAXright\MAXright0pt
\def\sdynkin{\bgroup\eightpoint\dynkin}
\def\endsdynkin{\enddynkin\egroup}
\def\dynkin{\bgroup\FIRSTtrue\hskip.5em\setbox1\hbox{$\diagup$}%
	\setbox2\hbox{$\diagdown$}%
	\setbox0\hbox to2\wd1{\hrulefill}%
	\setbox3\hbox{$\bullet$}%
	\setbox4\hbox{$\times$}%
	\setbox7\hbox{$\circ$}%       (L.K.)
	\def\whiteroot##1{\ifFIRST\setbox5\hbox{$##1$}\ifdim\wd5>1.3em%       (L.K.)
		\hskip-.5em\hskip.5\wd5\fi\fi\FIRSTfalse%                             (L.K.)
		\hskip-.25em\raise1.5\wd3\hbox to0pt{\hss\hskip.45em$%                (L.K.)
			\ssize##1$\hss}\copy7\hskip-.25em\setbox6\hbox{$##1$}%                (L.K.)
		\MAXright\wd6}%                                                       (L.K.)
	\def\root##1{\ifFIRST\setbox5\hbox{$##1$}\ifdim\wd5>1.3em%
		\hskip-.5em\hskip.5\wd5\fi\fi\FIRSTfalse%
		\hskip-.25em\raise1.5\wd3\hbox to0pt{\hss\hskip.45em$%
			\ssize##1$\hss}\copy3\hskip-.25em\setbox6\hbox{$##1$}%
		\MAXright\wd6}%
	\def\whitedroot##1{\ifFIRST\setbox5\hbox{$##1$}\ifdim\wd5>1.3em% (L.K.)
		\hskip-.5em\hskip.5\wd5\fi\fi\FIRSTfalse% (L.K.)
		\hskip-.25em\lower1.8\wd3\hbox to0pt{\hss\hskip.45em$%  (L.K.)
			\ssize##1$\hss}\copy7\hskip-.25em\setbox6\hbox{$##1$}% (L.K.)
		\MAXright\wd6}%
	\def\whiterroot##1{\hskip-.25em\copy7\hbox to0pt{\hskip.3em$\ssize##1$\hss}%
		\hskip-.25em\setbox6\hbox{\hskip.6em$##1##1$}%
		\MAXright\wd6}%
	\def\droot##1{\ifFIRST\setbox5\hbox{$##1$}\ifdim\wd5>1.3em%
		\hskip-.5em\hskip.5\wd5\fi\fi\FIRSTfalse%
		\hskip-.25em\lower1.8\wd3\hbox to0pt{\hss\hskip.45em$%
			\ssize##1$\hss}\copy3\hskip-.25em\setbox6\hbox{$##1$}%
		\MAXright\wd6}%
	\def\rroot##1{\hskip-.25em\copy3\hbox to0pt{\hskip.3em$\ssize##1$\hss}%
		\hskip-.25em\setbox6\hbox{\hskip.6em$##1##1$}%
		\MAXright\wd6}%
	\def\norroot##1{\hskip-.36em\copy4\hbox to0pt{\hskip.3em$\ssize##1$\hss}%
		\hskip-.48em\setbox6\hbox{\hskip.6em$##1##1$}%
		\MAXright\wd6}%
	\def\noroot##1{\ifFIRST\setbox5\hbox{$##1$}\ifdim\wd5>1.3em%
		\hskip-.5em\hskip.5\wd5\fi\fi\FIRSTfalse%
		\hskip-.36em\raise1.5\wd3\hbox to0pt{\hss\hskip.6em$%
			\ssize##1$\hss}\copy4\hskip-.38em\setbox6\hbox{$##1$}%
		\MAXright\wd6}%
	\def\nodroot##1{\ifFIRST\setbox5\hbox{$##1$}\ifdim\wd5>1.3em%
		\hskip-.5em\hskip.5\wd5\fi\fi\FIRSTfalse%
		\hskip-.36em\lower1.8\wd3\hbox to0pt{\hss\hskip.6em$%
			\ssize##1$\hss}\copy4\hskip-.38em\setbox6\hbox{$##1$}%
		\MAXright\wd6}%
	\def\nolink{\hskip\wd0}%      (L.K.)
	\def\link{\raise.22em\copy0}%
	\def\llink##1{\raise.32em\copy0\hskip-\wd0%
		\raise.12em\copy0\hskip-.5\wd0\hbox to0pt{\hss$##1$\hss}\hskip.5\wd0}%
	\def\lllink##1{\raise.22em\copy0\hskip-\wd0\raise.32em\copy0\hskip0.5\wd0%
		\raise.12em\copy0\hskip-.5\wd0\hbox to0pt{\hss$##1$\hss}\hskip1.5\wd0}%
	\def\rootupright##1{\hbox to0pt{\raise.45em\copy1\hskip-.55em\raise6.3\ht1%
			\hbox{\copy3\hskip.3em$\ssize##1$}\hss}%
		\setbox6\hbox{\hskip.6em\copy1\copy1$##1##1$}%
		\ifdim\MAXright<\wd6\MAXright\wd6\fi}%
	\def\whiterootupright##1{\hbox to0pt{\raise.45em\copy1\hskip-.25em\raise1.3\ht1% (L.K.)
			\hbox{\copy7\hskip.3em$\ssize##1$}\hss}% (L.K.)
		\setbox6\hbox{\hskip.6em\copy1\copy1$##1##1$}% (L.K.)
		\ifdim\MAXright<\wd6\MAXright\wd6\fi}% (L.K.)
	\def\norootupright##1{\hbox to0pt{\raise.45em\copy1\hskip-.36em\raise1.3\ht1%
			\hbox{\copy4\hskip.3em$\ssize##1$}\hss}%
		\setbox6\hbox{\hskip.6em\copy1\copy1$##1##1$}%
		\ifdim\MAXright<\wd6\MAXright\wd6\fi}%
	\def\rootdownright##1{\hbox to0pt{\raise-.5em\copy2\hskip-.25em\raise-1.35\ht1%
			\hbox{\copy3\hskip.3em$\ssize##1$}\hss}\setbox6%
		\hbox{\hskip.6em\copy2\copy2$##1##1$}%
		\ifdim\MAXright<\wd6\MAXright\wd6\fi}%
	\def\whiterootdownright##1{\hbox to0pt{\raise-.5em\copy2\hskip-.25em\raise-1.35\ht1% (L.K.)
			\hbox{\copy7\hskip.3em$\ssize##1$}\hss}\setbox6% (L.K.)
		\hbox{\hskip.6em\copy2\copy2$##1##1$}% (L.K.)
		\ifdim\MAXright<\wd6\MAXright\wd6\fi}% (L.K.)
	\def\rootdown##1{\hbox to0pt{\hskip-.05em\vrule height.25em depth.65em%
			\hskip-.25em\raise-.95em\hbox{\copy3\hskip.3em$\ssize##1$}\hss}%
		\setbox6\hbox{$##1$}%
		\ifdim\MAXright<\wd6\MAXright\wd6\fi}%
		\def\rootup##1{\hbox to0pt{\hskip-.05em\vrule height2.6em depth-.5em%
			\hskip-.25em\raise2.5em\hbox{\copy7\hskip.3em$\ssize##1$}\hss}%
		\setbox6\hbox{$##1$}%
		\ifdim\MAXright<\wd6\MAXright\wd6\fi}%
	\def\whiterootdown##1{\hbox to0pt{\hskip-.05em\vrule height.25em depth.65em% (L.K.)
			\hskip-.25em\raise-.95em\hbox{\copy7\hskip.3em$\ssize##1$}\hss}% (L.K.)
		\setbox6\hbox{$##1$}% (L.K.)
		\ifdim\MAXright<\wd6\MAXright\wd6\fi}% (L.K.)
	\def\dots{\hskip.5em\cdots\hskip.5em}}%
\def\enddynkin{\ifdim\MAXright>1em\hskip.5\MAXright\else\hskip.5em\fi\egroup}%

\newtheorem{theorem}[equation]{Theorem}
\newtheorem{definition}[equation]{Definition}

\newtheorem{lemma}[equation]{Lemma}

\newtheorem{proposition}[equation]{Proposition}
\newtheorem{remark}[equation]{Remark}

\numberwithin{equation}{section}

\newcommand{\be}{\begin{eqnarray}}
\newcommand{\ee}{\end{eqnarray}}
\newcommand{\bd}{\begin{definition}}
\newcommand{\ed}{\end{definition}}

\newcommand{\br}{\begin{remark}}
\newcommand{\er}{\end{remark}}

\newcommand{\bt}{\begin{tabular}}
\newcommand{\et}{\end{tabular}}

\newcommand \unb[2]{\underset{#1}{{\underbrace{#2}}}}

\def\ad{\operatorname{ad}}

\def\ker{\operatorname{Ker}}
\def\im{\operatorname{Im}}

\def\top{\operatorname{top}}

\def\Sp{\operatorname{Sp}}

\def\rank{\operatorname{rank}}

\newcommand{\Kt}{\tilde{K}}
\newcommand{\wti}{\widetilde}

\newcommand{\pf}{\begin{proof}}
\newcommand{\epf}{\end{proof}}
\newcommand{\eq}{\begin{equation}}
\newcommand{\eeq}{\end{equation}}
\newcommand{\eqn}{\begin{equation*}}
\newcommand{\eeqn}{\end{equation*}}

\newcommand{\eps}{\varepsilon}
\newcommand{\half}{\frac{1}{2}}

\newcommand{\frc}{\mathfrak{c}}
\newcommand{\frd}{\mathfrak{d}}

\newcommand{\frg}{\mathfrak{g}}
\newcommand{\frh}{\mathfrak{h}}

\newcommand{\frk}{\mathfrak{k}}
\newcommand{\frl}{\mathfrak{l}}

\newcommand{\fro}{\mathfrak{o}}
\newcommand{\frp}{\mathfrak{p}}
\newcommand{\frqqq}{\mathfrak{q}}

\newcommand{\frs}{\mathfrak{s}}

\newcommand{\fru}{\mathfrak{u}}

\newcommand{\bbC}{\mathbb{C}}

\newcommand{\bbR}{\mathbb{R}}

\newcommand{\caO}{\mathcal{O}}

\newcommand{\Lamdba}{\Lambda}

\input xypic
\input epsf
\input xy
\xyoption{all}

\begin{document}

\title[Computing the associatied cycles]{Computing the associatied cycles\\ of certain Harish-Chandra modules}
\author{Salah Mehdi}
\author{Pavle Pand\v zi\'c}
\author{David Vogan}
\author{Roger Zierau}
\address[Mehdi]{Institut Elie Cartan de Lorraine, CNRS - UMR 7502, Universit\'e de Lorraine, Metz, F-57045, France.}
\email{salah.mehdi@univ-lorraine.fr}
\address[Pand\v zi\'c]{Department of Mathematics, Faculty of Science, University of Zagreb, Zagreb, Croatia.}
\email{pandzic@math.hr}
\address[Vogan]{Department of Mathematics, Massachusetts Institute of Technology, Cambridge, MA 02139, USA}
\email{dav@math.mit.edu}
\address[Zierau]{Department of Mathematics, Oklahoma State University, Stillwater, Oklahoma 74078, USA}
\email{roger.zierau@okstate.edu}

\thanks{The second named author was supported by grant no. 
4176 of the Croatian Science Foundation and 
by the QuantiXLie Centre of Excellence, a project
cofinanced by the Croatian Government and European Union through the European Regional Development Fund - the Competitiveness and Cohesion Operational Programme (KK.01.1.1.01.0004). The third named author was supported 
in part by NSF grant DMS 0967272.}

\begin{abstract}
Let $G_{\mathbb{R}}$ be a simple real linear Lie group with 
maximal compact subgroup $K_{\mathbb{R}}$ and assume that $\rank(G_\bbR)=\rank(K_\bbR)$. In \cite{MPVZ} we proved that for any representation $X$ of Gelfand-Kirillov dimension $\frac{1}{2}\dim(G_{\mathbb{R}}/K_{\mathbb{R}})$, the polynomial on the dual of a compact Cartan subalgebra  given by the dimension of the Dirac index of members of the coherent family containing $X$ is a linear combination, with integer coefficients,  of the multiplicities of the irreducible components occurring in the associated cycle.
In this paper we compute these coefficients explicitly. 
\end{abstract}

\maketitle

%%%%%%%%%%%%%%%%%%%%%%%%%%%%%%%%%%%%%%%%%%%%%%%%%%%%%%%%%%%%%%%%%%%%%%%%%%%%%%
%%%%%%%%%%%%%%%%%%%%%%%%%%%%%%%%%%%%%%%%%%%%%%%%%%%%%%%%%%%%%%%%%%%%%%%%%%%%%%

\section{Introduction}

Let $G_{\mathbb{R}}$ be a simple real linear Lie group with a Cartan involution $\theta$  and maximal compact subgroup $K_{\mathbb{R}}=G_\bbR^\theta$. Let $\frg=\frk\oplus\frp$ be the Cartan decomposition of the  complexified Lie algebra $\frg$ of $G_\bbR$; this decomposition is orthogonal with respect to the Killing form $B$. Let $K$ be the complexification of $K_\bbR$ and $G$ a complex Lie group (with Lie algebra $\frg$) containing $K$ as the set of fixed points of the complex extension of $\theta$. We assume throughout the paper that $\frg$ and $\frk$ have equal rank, i.e., there is a Cartan subalgebra $\frh$ of $\frg$ contained in $\frk$. We fix such $\frh$ and write $W$ for the Weyl group of $(\frg,\frh)$. 

In this paper we are concerned with comparing two important invariants of $(\frg,K_\bbR)$-modules. One is the Dirac index studied in \cite{MPV}. It is defined using the Dirac operator $D\in U(\frg)\otimes C(\frp)$, where $U(\frg)$ is the universal enveloping algebra of $\frg$ and $C(\frp)$ is the Clifford algebra of $\frp$ with respect to $B$. If $M$ is a $(\frg,K_\bbR)$-module, then $D$ acts on $M\otimes S$ where $S$ is a spin module for $C(\frp)$. The Dirac cohomology of $M$ is defined as
\[
H_D(M)=\ker(D)/(\im(D)\cap\ker(D);
\]
it is a module for the spin double cover $\Kt$ of $K_\bbR$ (finite-dimensional if $M$ is admissible). This invariant was introduced in \cite{V2}, it turned out to be very interesting and also quite computable; see for example \cite{HP1}, \cite{HP2}, \cite{HKP}, \cite{HPR}, \cite{HPP}, \cite{HPZ}, \cite{BP1}, \cite{BP2}, \cite{BPT}, \cite{MP}, \cite{MZ}.

Decomposing the $\Kt$-module $S$ as $S=S^+\oplus S^-$ induces a decomposition of Dirac cohomology
\[
H_D(M)=H_D(M)^+\oplus H_D(M)^-.
\]
The Dirac index of $M$ is then defined as the virtual $\Kt$-module
\[
DI_v(M)=H_D(M)^+ - H_D(M)^-.
\]
It is proved in \cite{MPV} that Dirac index varies nicely over coherent families of $(\frg,K_\bbR)$-modules. In particular, if $\{M_\lambda\}$ is such a coherent family, attached to a module $M$, then the function
\[
\lambda\mapsto \dim DI_v(M_\lambda)
\]
extends to a polynomial on $\frh^*$, which we denote by $DI_p(M)$.

Another very useful invariant of a Harish-Chandra module $M$ is its associated cycle $AC(M)$, defined in \cite{V1}. See \cite{MPVZ} for a short review of the definition. 

In concrete terms, for irreducible $M$, $AC(M)$ can be written as the formal sum
\[
AC(M)=\sum_i m_i(M)\overline{\caO_i},
\]
where $\caO_i\subset\frp$ are the real forms of a complex nilpotent $G$-orbit $\caO^\bbC\subset\frg$, and the multiplicities $m_i(M)$ are nonnegative integers. The orbit $\caO^\bbC$  is specified by the requirement that $\overline{\caO^\bbC}$ is the associated variety of the annihilator of $M$. 

If $M$ is put into a coherent family $\{M_\lambda\}$, then the corresponding multiplicities extend to polynomials $m_i(M)$ on $\frh^*$. It was conjectured in \cite{MPV}, and proved in \cite{MPVZ}, that in certain special circumstances these multiplicity polynomials are related to the Dirac index polynomial by 
\[
DI_p(M)=\sum_i c_i m_i(M)
\]
for some integers $c_i$. 
Such a relationship is true when the associated variety of the annihilator of $M$ is contained in $\overline{\caO^\bbC}$, with 
$\caO^\bbC$ corresponding via Springer correspondence to the $W$-representation generated by the Weyl dimension polynomial $P_K$ for $K$ ($P_K$ is defined by \eqref{def pk}).

The purpose of this paper is to complement \cite{MPVZ} by explicitly computing the constants $c_i$ in the classical cases other than $SU(p,q)$. The case $G_\bbR=SU(p,q)$ as well as the case of exceptional groups are done in \cite{MPVZ}

We start by reviewing some facts about real forms of nilpotent orbits in Section 2, and assembling a few useful general facts about the computations in Section 3. Then we do the case-by-case computations in Sections 4 -- 8.

\section{Nilpotent orbits and their real forms}

We recall that the list of the classical real groups for which the conjecture from \cite{MPV} applies is given in \cite{MPVZ}, Section 6, Table 1, along with the relevant explanations. The groups on the list are the connected classical equal rank groups such that the $W$-representation $\sigma_K$ generated by the Weyl dimension polynomial $P_K$ for $K$ is Springer. 
The list consists of 
\begin{itemize}
\item $SU(p,q),\quad q\geq p\geq 1$;
\item $SO_e(2p,2q+1),\quad q\geq p-1\geq 0$;
\item $\Sp(2n,\bbR),\quad n\geq 1$;
\item $SO^*(2n),\quad n\geq 1$;
\item $SO_e(2p,2q),\quad q\geq p\geq 1$.
\end{itemize}
The table in \cite{MPVZ} also includes the nilpotent orbits $\caO^\bbC$ corresponding to $\sigma_K$ in each of the cases, as well as the number of real forms of these orbits. Here we explain how to get these real forms, and in particular how to write down the semisimple elements $h$ of the corresponding $\frs\frl_2$-triples, which we need to begin our computations.

We start by recalling that complex nilpotent orbits in classical Lie algebras are in
one-to-one correspondence with the set of partitions  
$\lbrack d_1,\cdots,d_k\rbrack$ with $d_1\geq d_2\geq\cdots\geq
d_k\geq 1$ (if $d_j$ occurs $m$ times, we will simply write $d_j^m$) 
such that (see \cite{CM}, Chapter 5): 
\begin{itemize}
%  \item[]
\item[$\bullet$] $d_{1}+d_{2}+\cdots+d_{k}=n$, when
  ${\mathfrak g}\simeq{\mathfrak s}{\mathfrak l}(n,{\mathbb C})$; 
\item[$\bullet$] $d_{1}+d_{2}+\cdots+d_{k}=2n+1$ and the even $d_j$
  occur with even multiplicity, when ${\mathfrak g}\simeq{\mathfrak s}{\mathfrak o} (2n+1,{\mathbb C})$; 
\item[$\bullet$] $d_{1}+d_{2}+\cdots+d_{k}=2n$ and the odd $d_j$ occur
  with even multiplicity, when ${\mathfrak g}\simeq{\mathfrak s}{\mathfrak p}(2n,{\mathbb C})$; 
\item[$\bullet$] $d_{1}+d_{2}+\cdots+d_{k}=2n$ and the even $d_j$
  occur with even multiplicity, when ${\mathfrak g}\simeq{\mathfrak s}{\mathfrak o} (2n,{\mathbb C})$;
  except that the partitions having all the
  $d_j$ even and occurring with even multiplicity are each 
  associated to {\em two} orbits. 
\end{itemize}

We now recall the procedure which attaches ${\mathfrak s}{\mathfrak l}_2$-triples to complex nilpotent orbits (see Chapter 3 in \cite{CM}). For a positive integer $i$, define the Jordan block $J_i$ to be the 
$i\times i$ matrix 
\vspace*{0.2cm}
\begin{equation*}
J_i=\begin{pmatrix} 0 & 1 & 0 & 0 &\cdots & 0\\ 0 & 0 & 1 & 0 & \cdots &0\\
\cdot & \cdot & \cdot & \cdot & \cdots & \cdot\\
\cdot & \cdot & \cdot & \cdot & \cdots & \cdot\\
0 & 0 & 0 & \cdots & 0 & 1\\
0 & 0 & 0 & 0 & \cdots & 0\end{pmatrix}
\end{equation*}
For a postive integer $n$, write $\lbrack d_1,d_2,\cdots,d_k\rbrack$ for a partition of $n$. Define the $n\times n$ matrix 
\vspace*{0.2cm}
\begin{equation*}
X_{\lbrack d_1,d_2,\cdots,d_k\rbrack}=\begin{pmatrix} J_{d_1} & 0 & 0 & 0 &\cdots & 0\\ 
0 & J_{d_2} & 0 & 0 &\cdots &0\\
\cdot & \cdot & \cdot & \cdot & \cdots & \cdot\\
\cdot & \cdot & \cdot & \cdot & \cdots & \cdot\\
0 & 0 & 0 & 0 & \cdots & J_{d_k}\end{pmatrix}
\end{equation*}
Then $X_{\lbrack d_1,d_2,\cdots,d_k\rbrack}$ is a nilpotent element in the complex Lie algebra ${\mathfrak s}{\mathfrak l}_n$. Write ${\mathcal O}_{\lbrack d_1,d_2,\cdots,d_k\rbrack}$ for the complex nilpotent orbit under the adjoint group $PSL_n$ of ${\mathfrak s}{\mathfrak l}_n$. It is convenient to attach to ${\mathcal O}_{\lbrack d_1,d_2,\cdots,d_k\rbrack}$ a Young tableau, i.e a left-justified arrangement of empty boxes of rows with size in the non-increasing order $d_1$, $d_2$, ..., $d_k$. 

Let $H=\begin{pmatrix} 1&0\\ 0&-1\end{pmatrix}$, $X=J_2=\begin{pmatrix} 0&1\\ 0&0\end{pmatrix}$ and $Y=J_2^t=\begin{pmatrix} 0&0\\ 1&0\end{pmatrix}$. Then $\lbrack H,X\rbrack=2X$, $\lbrack H,Y\rbrack=-2Y$ and $\lbrack X,Y\rbrack=H$, so that $H,X,Y$ span, over ${\mathbb C}$, the simple Lie algebra ${\mathfrak s}{\mathfrak l}_2$ of $2\times 2$ complex matrices with zero trace. For a non-negative integer $r$, define the linear map $\rho_r:{\mathfrak s}{\mathfrak l}_2\rightarrow{\mathfrak s}{\mathfrak l}_{r+1}$ by 
\begin{eqnarray*}
&&\\
\rho_r(H)&=&\begin{pmatrix} r & 0 & 0 & 0 &\cdots & 0\\ 
0 & r-2 & 0 & 0 & \cdots &0\\
\cdot & \cdot & \cdot & \cdot & \cdots & \cdot\\
\cdot & \cdot & \cdot & \cdot & \cdots & \cdot\\
0 & \cdot & \cdot & \cdot &-r+2 & 0\\
0 & \cdot & \cdot & \cdot & 0 & -r\end{pmatrix}\\
&&\\
\rho_r(X)&=&J_{r+1}\\
&&\\
\rho_r(Y)&=&\begin{pmatrix} 0 & 0 & 0 & 0 &\cdots & 0\\ 
\mu_1 & 0 & 0 & 0 & \cdots &0\\
0 & \mu_2 & 0 & 0 & \cdots &0\\
\cdot & \cdot & \cdot & \cdot & \cdots & \cdot\\
\cdot & \cdot & \cdot & \cdot & \cdots & \cdot\\
0 & \cdot & \cdot & \cdot &\mu_{r}& 0\end{pmatrix}\\
&&\\
&&\text{with }\mu_i=i(r+1-i)\text{ for }1\leq i\leq r
\end{eqnarray*}
$\rho_r$ defines an irreducible representation of ${\mathfrak s}{\mathfrak l}_2$ of dimension $r+1$, and any finite dimensional irreducible 
representation of ${\mathfrak s}{\mathfrak l}_2$ is equivalent to $\rho_r$ for some $r$. The map $\rho_r$  induces the homomorphism 
$\Phi_{\mathcal O}:{\mathfrak s}{\mathfrak l}_2\rightarrow{\mathfrak s}{\mathfrak l}_{n}$ defined by:
\begin{equation*}
\Phi_{\mathcal O}=\bigoplus_{1\leq j\leq k}\rho_{d_j-1}
\end{equation*}
so that $\Phi_{\mathcal O}(X)=X_{\lbrack d_1,d_2,\cdots,d_k\rbrack}$. The standard ${\mathfrak s}{\mathfrak l}_2$-triple associated with the complex nilpotent orbit ${\mathcal O}_{\lbrack d_1,d_2,\cdots,d_k\rbrack}$ is $\{H_{\lbrack d_1,d_2,\cdots,d_k\rbrack};X_{\lbrack d_1,d_2,\cdots,d_k\rbrack};Y_{\lbrack d_1,d_2,\cdots,d_k\rbrack}\}$ where $H_{\lbrack d_1,d_2,\cdots,d_k\rbrack}:=\Phi_{\mathcal O}(H)$, 
$X_{\lbrack d_1,d_2,\cdots,d_k\rbrack}:=\Phi_{\mathcal O}(X)$ and 
$Y_{\lbrack d_1,d_2,\cdots,d_k\rbrack}:=\Phi_{\mathcal O}(Y)$. Choose the Cartan subalgebra consisting of $n\times n$ diagonal matrices $\text{diag}(a_1,a_2,\cdots,a_n)$ with zero trace, and fix the positive system of roots $\{\epsilon_i-\epsilon_j\mid 1\leq i<j\leq n\}$ whose corresponding Borel subalgebra consists of the upper tringular matrices of zero trace. Here $\epsilon_i$ is the complex linear form defined on the Cartan subalgebra such that $\epsilon_i(\text{diag}(a_1,a_2,\cdots,a_n))=a_i$. Then up to a Weyl group element of ${\mathfrak s}{\mathfrak l}_{n}$, the element $H_{\lbrack d_1,d_2,\cdots,d_k\rbrack}$  is conjugate to a dominant element 
\begin{equation*}
(h_1,h_2,\cdots,h_n):=\text{diag}(h_1,h_2,\cdots,h_n)
\end{equation*}
with $h_1\geq h_2\geq \cdots\geq h_n$ and $h_1+h_2+\cdots + h_n=0$. Associated with the orbit ${\mathcal O}_{\lbrack d_1,d_2,\cdots,d_k\rbrack}$ is the weighted Dynkin diagram 
\vskip 2mm
\[
\dynkin \whiteroot{h_1-h_2}\hspace*{0.5mm}\link\link\hspace*{0.5mm}\whiteroot{h_2-h_3}\hspace*{0.5mm}\link\link\cdots\cdots
\link\link\hspace*{0.5mm}\whiteroot{h_{n-1}-h_{n}}
\enddynkin
\]
\vskip 3mm

Suppose now that $G_{\mathbb R}=SU(p,q)$, $K_{\mathbb R}=S(U(p)\times U(q))$ and 
${\mathfrak g}={\mathfrak s}{\mathfrak l}_{p+q}$, with $q\geq p\geq 1$ and $p+q=n$. The dominant $h$ associated with the complex nilpotent orbit ${\mathcal O}^{\mathbb C}={\mathcal O}_{\lbrack 2^p;1^{q-p}\rbrack}$ is given by 
\begin{equation*}
h=(\underbrace{1,1,\cdots,1}_{p},\underbrace{0,0,\cdots,0}_{q-p},\underbrace{-1,-1,\cdots,-1}_{p})
\end{equation*} 
along with the weighted Dynkin diagram
\vskip 2mm
\[
\dynkin \underset{\epsilon_1-\epsilon_2}{\whiteroot{0}}\hspace*{-3mm}\link\cdots\cdots\link\hspace*{-4.5mm}\underset{\epsilon_p-\epsilon_{p+1}}{\whiteroot{1}}\hspace*{-4.5mm}\link\link\hspace*{-6mm}\underset{\epsilon_{p+1}-\epsilon_{p+2}}{\whiteroot{0}}\hspace*{-6mm}\link\cdots\cdots
\link\hspace*{-4.5mm}\underset{\epsilon_q-\epsilon_{q+1}}{\whiteroot{1}}\hspace*{-4.5mm}\link\link\hspace*{-6mm}\underset{\epsilon_{q+1}-\epsilon_{q+2}}{\whiteroot{0}}\hspace*{-6mm}\link\cdots\cdots
\link\hspace*{-7.5mm}\underset{\epsilon_{p+q}-\epsilon_{p+q-1}}{\whiteroot{0}}
\enddynkin
\]
\vskip 3mm

Moreover, by Kostant-Sekiguchi, it is known that the nilpotent orbit 
${\mathcal O}^{\mathbb C}$ has $p+1$ real forms which are in one to one correspondence with nilpotent $K$-orbits in ${\mathfrak p}$. More precisely, 
for $k=0,1,2,\cdots,p$, let $I_0=0$ and
\vspace*{0.1cm}\\
\[
\sbox0{$\begin{matrix}I_k&0\\0&0\end{matrix}$}
\sbox1{$\begin{matrix}0&0\\0&I_{p-k}\end{matrix}$}
e_k=\left[
\begin{array}{c|c}
\makebox[\wd0]{\large $0_p$}&\usebox{0}\\
\hline
  \vphantom{\usebox{0}}\makebox[\wd0]{\usebox{1}}&\makebox[\wd0]{\large $0_q$}
\end{array}
\right]
\]
\vspace*{0.1cm}\\
For each $k$, the element $e_k$ belongs to ${\mathcal O}^{\mathbb C}$. On the other hand, the $K$-orbit of $e_k$ consists of matrices of the form 
\begin{equation*}
\begin{pmatrix}
0&A\\ B&0
\end{pmatrix}\;\;\text{ with }\text{rank}(A)=k\text{ and }\text{rank}(B)=p-k
\end{equation*}
In particular, if $k\neq k^{\prime}$ then the $K$-orbits of $e_k$ and 
$e_{k^{\prime}}$ are disjoint. Choosing the positive system $\{\epsilon_i-\epsilon_j\mid 1\leq i<j\leq p\text{ or }p+1\leq i<j\leq p+q\}$, the (dominant) neutral element of the ${\mathfrak s}{\mathfrak l}_2$-triple corresponding to the real form $K\cdot e_k$ is 
\vskip 1mm
\begin{equation*}
h_k=(\underbrace{1,1,\cdots,1}_{k},\underbrace{-1,-1,\cdots,-1}_{p-k},\underbrace{1,1,\cdots,1}_{p-k},\underbrace{0,0,\cdots,0}_{q-p},\underbrace{-1,-1,\cdots,-1}_{k})
\end{equation*} 
\vskip 1mm
\noindent
The description in terms of Young tableaux of the complex orbit ${\mathcal O}^{\mathbb C}$ and of its real forms is as follows:
\vspace{0.2cm}\\

%\[
$\begin{array}{rl}
{\mathcal O}^{\mathbb C}\hspace*{0.2cm}
\left.\young(\empty\;\empty,::,\empty\;\empty)\right\} & \hspace{-0.3cm}p
\end{array}$
\\[-1.17mm]
\hspace*{1.049cm}
$\begin{array}{rl}
\left.\young(\empty,:,\empty)\right\} & \hspace{-0.3cm}q-p\\[-0.8mm]
\end{array}$
%\]

\vspace*{-6.8em}
%\[
$\begin{array}{rl}
\hspace*{5cm}\text{ real forms for }{\mathcal O}^{\mathbb C}\hspace*{0.2cm}\left.\young(+-,::,+-)\right\} & \hspace{-0.3cm}k\in\{0,\ldots,p\} \\[-0.75mm]
\left.\young(-+,::,-+)\right\} & \hspace{-0.3cm}p-k \\[-0.8mm]
\end{array}$
\\[-0.439mm]
\hspace*{8.355cm}
$\begin{array}{rl}
\left.\young(-,:,-)\right\} & \hspace{-0.3cm}q-p
\end{array}$
%\]
\vspace*{0.2cm}\\

Consider $G_{\mathbb R}=\Sp(2n,{\mathbb R})$ and $K_{\mathbb R}=U(n)$. The complexification ${\mathfrak g}={\mathfrak s}{\mathfrak p}_{2n}$ of $G_{\mathbb R}$ is realized as the following set of matrices 
\begin{equation*}
\Big\{\begin{pmatrix} Z_1 & Z_2\\ Z_3&-Z_1^t\end{pmatrix}\mid Z_1\text{ $n\times n$ complex matrices, }Z_2,Z_3\text{ symmetric complex matrices}\Big\}
\end{equation*}
A Cartan subalgebra in ${\mathfrak g}$ consists of diagonal complex matrices of the form\\$\text{diag}(a_1,a_2,\cdots,a_n,-a_1,-a_2,\cdots,-a_n)$. Fix the standard system of positive roots $\{\epsilon_i\pm\epsilon_j\mid 1\leq i<j\leq n\}\cup\{2\epsilon_k\mid 1\leq k\leq n\}$. As in type A, there is an explicit recipe which attaches an ${\mathfrak s}{\mathfrak l}_2$-triple to a complex nilpotent orbit (see 5.2.2 in \cite{CM}). We apply this recipe to the nilpotent orbit ${\mathcal O}^{\mathbb C}={\mathcal O}_{\lbrack 2^n\rbrack}$, using $n$ chunks coninciding with $\{2\}$. We obtain (viewing ${\mathfrak s}{\mathfrak p}_{2n}$ as a subalgebra of 
${\mathfrak s}{\mathfrak l}_{2n}$) 
\vskip 1mm
\begin{equation*}
h:=\text{diag}(\underbrace{1,1,\cdots,1}_{n},\underbrace{-1,-1,\cdots,-1}_{n})
\end{equation*} 
\vskip 1mm
\noindent
which we will simply write 
\vskip 1mm
\begin{equation*}
h=(\underbrace{1,1,\cdots,1}_{n})
\end{equation*} 
along with the weighted Dynkin diagram
\vskip 2mm
\[
\dynkin \underset{\epsilon_1-\epsilon_2}{\whiteroot{0}}\hspace*{-3mm}\link\link\hspace*{-2.7mm}\underset{\epsilon_2-\epsilon_3}{\whiteroot{0}}\hspace*{-3mm}\link\hspace*{0.5mm}\cdots\cdots
\link\hspace*{-4.5mm}\underset{\epsilon_{n-1}-\epsilon_{n}}{\whiteroot{0}}\hspace*{-5.1mm}\link\link\hspace*{-1.5mm}\underset{2\epsilon_{n}}{\whiteroot{2}}
\enddynkin
\]
\vskip 3mm
The same argument as in type $A$ shows that ${\mathcal O}^{\mathbb C}$ possesses $n+1$ real forms with 
\vskip 1mm
\begin{equation*}
h_k=(\underbrace{1,1,\cdots,1}_{k},\underbrace{-1,-1,\cdots,-1}_{n-k},\underbrace{-1,-1,\cdots,-1}_{k},\underbrace{1,1,\cdots,1}_{n-k})
\end{equation*} 
\vskip 1mm
\noindent
The description in terms of Young tableaux of the complex orbit ${\mathcal O}^{\mathbb C}$ and of its real forms is as follows:
\vspace{0.1cm}\\

%\[
$\begin{array}{rl}
{\mathcal O}^{\mathbb C}\hspace*{0.2cm}
\left.\young(\empty\;\empty,\empty\;\empty,::,\empty\;\empty,\empty\;\empty)\right\} & \hspace{-0.3cm}n
\end{array}$
%\]

\vspace*{-2cm}
%\[
$\begin{array}{rl}
\vspace*{-0.7cm}\hspace*{4cm}\text{ real forms }
\end{array}$
\vspace*{-0.2cm}
\hspace*{0.1cm}
$\begin{array}{rl}
\hspace*{-0.1cm}\left.\young(+-,::,+-)\right\} & \hspace{-0.3cm}k\in\{0,\ldots,n\} \\[-0.8mm]
\end{array}$
\hspace*{-3.94cm}
$\begin{array}{rl}
\vspace*{-2.45cm}\left.\young(-+,::,-+)\right\} & \hspace{-2mm}n-k
\end{array}$
%\]
\vspace*{2cm}\\

Consider $G_{\mathbb R}=SO_e(2p,2q+1)$ and $K_{\mathbb R}=S(O(2p)\times O(2q+1))$. The complexification ${\mathfrak g}={\mathfrak s}{\mathfrak o}_{2n+1}$ of $G_{\mathbb R}$, with $n=p+q$ and $q\geq p\geq 1$, is realized as the following set of matrices 
\vskip1mm
\begin{equation*}
\Big\{\begin{pmatrix} 0&u&v\\ -v^t&Z_1 & Z_2\\ -u^t & Z_3 & -Z_1^t\end{pmatrix}\mid u,v\in{\mathbb C}^n, Z_1\text{ $n\times n$ complex matrices, }Z_2,Z_3\text{ skew-symmetric}\Big\}
\end{equation*}
\vskip1mm
\noindent
A Cartan subalgebra in ${\mathfrak g}$ consists of diagonal complex matrices of the form\\$\text{diag}(0,a_1,a_2,\cdots,a_n,-a_1,-a_2,\cdots,-a_n)$ (first row and column of zeros). Fix the standard system of positive roots $\{\epsilon_i\pm\epsilon_j\mid 2\leq i<j\leq n+1\}\cup\{\epsilon_k\mid 2\leq k\leq n+1\}$. There is an explicit recipe which attaches an ${\mathfrak s}{\mathfrak l}_2$-triple to a complex nilpotent orbit (see 5.2.4 in \cite{CM}). We apply this recipe to the nilpotent orbit ${\mathcal O}^{\mathbb C}={\mathcal O}_{\lbrack 3,2^{2p-2},1^{2(q-p+1)}\rbrack}$, using the following chunks : $\{3\}$, $p-1$ $\{2;2\}$'s and $q-p+1$ $\{1;1\}$'s. We obtain (viewing ${\mathfrak s}{\mathfrak o}_{2n+1}$ as a subalgebra of 
${\mathfrak s}{\mathfrak l}_{2n+1}$) 
\vskip 1mm
\begin{equation*}
h:=\text{diag}(0,2,\underbrace{1,1,\cdots,1}_{2p-2},\underbrace{0,0,\cdots,0}_{2(q-p+1)},\underbrace{-1,-1,\cdots,-1}_{2p-2},-2)
\end{equation*} 
\vskip 1mm
\noindent
which we will simply write (dropping the first zero coordinate and shifting indices of $\epsilon_i$'s) 
\vskip 1mm
\begin{equation*}
h=(2,\underbrace{1,1,\cdots,1}_{2p-2},\underbrace{0,0,\cdots,0}_{q-p+1})
\end{equation*} 
along with the weighted Dynkin diagram
\vskip 2mm
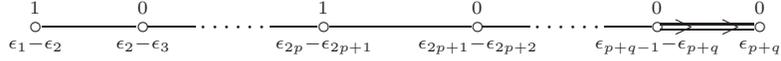
\[
\dynkin \underset{\epsilon_1-\epsilon_2}{\whiteroot{1}}\hspace*{-2.8mm}\link\link\hspace*{-2.7mm}\underset{\epsilon_2-\epsilon_{3}}{\whiteroot{0}}\hspace*{-3mm}\link\cdots\cdots\link\hspace*{-5.5mm}\underset{\epsilon_{2p}-\epsilon_{2p+1}}{\whiteroot{1}}\hspace*{-5.7mm}\link\link\link\hspace*{-7mm}\underset{\epsilon_{2p+1}-\epsilon_{2p+2}}{\whiteroot{0}}\hspace*{-7.3mm}\link\cdots\cdots\link\hspace*{-7.5mm}\underset{\epsilon_{p+q-1}-\epsilon_{p+q}}{\whiteroot{0}}\hspace*{-8mm}\llink>\llink>\hspace*{-2mm}\underset{\epsilon_{p+q}}{\whiteroot{0}}
\enddynkin
\]
\vskip 3mm
The nilpotent orbit ${\mathcal O}^{\mathbb C}$ possesses $2$ or $3$ real forms depending wether $q>p-1$ or not. The description in terms of Young tableaux of the complex orbit ${\mathcal O}^{\mathbb C}$ and of its real forms is given below. 
The recipe to produce the real $h$'s from the signed tableau can be stated as follows: the first row of lenght $3$ gives a $2$ in the first $p$ coordinates if the row starts with a "+" and a $2$ in the $p+1$ coordinate if it starts with a "-". For the rows of length two, a "+" (resp. "-") sign in the left-most box provides $+1$  in the first $p$ coordinates (resp. in the second group of coordinates $p+1,\cdots$); a "+" (resp. "-") sign in the right-most box provides $-1$ in the first $p$ coordinates (resp. in the second group of coordinates $p+1,\cdots$). In particular, we get
\vskip 1mm
\begin{eqnarray*}
h_1^{I}&=&(2,\underbrace{1,1,\cdots,1}_{p-1},\underbrace{1,1,\cdots,1}_{p-1},\underbrace{0,0,\cdots,0}_{q-p+1})\\
&&\\
h_1^{II}&=&(2,\underbrace{1,1,\cdots,-1}_{p-1},\underbrace{1,1,\cdots,1}_{p-1},\underbrace{0,0,\cdots,0}_{q-p+1})\\
&&\\
h_2&=&(\underbrace{1,1,\cdots,1}_{p-1},0,2,\underbrace{1,1,\cdots,1}_{p-1},\underbrace{0,0,\cdots,0}_{q-p})\text{ only if }q>p-1
\end{eqnarray*} 
\vskip 1mm
\noindent 
$h_1^{II}$ is obtained from $h_1^{I}$ by the outer automorphism 
$\epsilon_{p-1}+\epsilon_p\longleftrightarrow \epsilon_{p-1}-\epsilon_{p}$:
\vskip 2mm
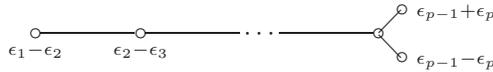
\[
\dynkin \underset{\epsilon_1-\epsilon_2}{\whiteroot{}}\hspace*{-3mm}\link\link\hspace*{-2.8mm}\underset{\epsilon_2-\epsilon_3}{\whiteroot{}}\hspace*{-3mm}\link\link\cdots\link\link\hspace*{0.7mm}\whiteroot{}\whiterootupright{\epsilon_{p-1}+\epsilon_{p}}
\whiterootdownright{\epsilon_{p-1}-\epsilon_{p}}
\enddynkin
\]
\vskip 3mm
The description of the orbit ${\mathcal O}^{\mathbb C}$ and its real forms in terms of Young tabelaux is as follows:
\newpage
\vspace*{0.5cm}
%\[
$\begin{array}{rl}
\hspace*{0.36cm}\young(\;\empty\;\empty)
\end{array}$
\\[-0.096cm]
$\begin{array}{rl}
{\mathcal O}^{\mathbb C}\hspace*{0.2cm}
\left.\young(\empty\;\empty,::,\empty\;\empty)\right\} & \hspace{-0.3cm}2p-2
\end{array}$
\\[-1.17mm]
\hspace*{0.63cm}
$\begin{array}{rl}
\left.\young(\empty,:,\empty)\right\} & \hspace{-0.3cm}2(q-p+1)\\[-0.8mm]
\end{array}$
%\]

\vspace*{-2.78cm}
%\vspace*{-1cm}
%\[
$\begin{array}{rl}
\hspace*{7.67cm}\young(+-+)
\end{array}$
\\[-0.096cm]
$\begin{array}{rl}
\hspace*{5cm}\text{ real forms for }{\mathcal O}^{\mathbb C}\hspace*{0.2cm}\left.\young(+-,::,+-)\right\} & \hspace{-0.3cm}2p-2
\hspace{0.5cm}\text{I, II}\\[-0.8mm]
\end{array}$
\\[-0.420mm]
\hspace*{7.935cm}
$\begin{array}{rl}
\left.\young(-,:,-)\right\} & \hspace{-0.3cm}2(q-p+1)
\end{array}$
%\]

\vspace*{0.5cm}
%\vspace*{-1cm}
%\[
$\begin{array}{rl}
\hspace*{7.67cm}\young(-+-)
\end{array}$
\\[-0.096cm]
$\begin{array}{rl}
\hspace*{7.85cm}\hspace*{0.2cm}\left.\young(+-,::,+-)\right\} & \hspace{-0.3cm}2p-2
\hspace{0.5cm}\text{only if }q>p-1\\[-0.8mm]
\end{array}$
\vspace*{0.02cm}
\hspace*{7.975cm}
$\begin{array}{rl}
\young(+)
\end{array}$
\vspace*{-0.11cm}\\
\hspace*{7.94cm}
$\begin{array}{rl}
\left.\young(-,:,-)\right\} & \hspace{-0.3cm}2(q-p)+1
\end{array}$
\vspace*{1cm}\\

Consider $G_{\mathbb R}=SO_e(2p,2q)$ and $K_{\mathbb R}=S(O(2p)\times O(2q))$. The complexification ${\mathfrak g}={\mathfrak s}{\mathfrak o}_{2n}$ of $G_{\mathbb R}$, with $n=p+q$ and $q\geq p\geq 1$, is realized as the following set of matrices 
\vskip1mm
\begin{equation*}
\Big\{\begin{pmatrix} Z_1 & Z_2\\  Z_3 & -Z_1^t\end{pmatrix}\mid Z_i\text{ $n\times n$ complex matrices, }Z_2,Z_3\text{ skew-symmetric}\Big\}
\end{equation*}
\vskip1mm
\noindent
A Cartan subalgebra in ${\mathfrak g}$ consists of diagonal complex matrices of the form\\$\text{diag}(a_1,a_2,\cdots,a_n,-a_1,-a_2,\cdots,-a_n)$. Fix the standard system of positive roots $\{\epsilon_i\pm\epsilon_j\mid 1\leq i<j\leq n\}$. There is an explicit recipe which attaches an ${\mathfrak s}{\mathfrak l}_2$-triple to a complex nilpotent orbit (see 5.2.6 in \cite{CM}). We apply this recipe to the nilpotent orbit ${\mathcal O}^{\mathbb C}={\mathcal O}_{\lbrack 3,2^{2p-2},1^{2(q-p)+1}\rbrack}$, using the following chunks : $\{3;1\}$, $p-1$ $\{2;2\}$'s and $q-p$ $\{1;1\}$'s. We obtain (viewing ${\mathfrak s}{\mathfrak o}_{2n}$ as a subalgebra of 
${\mathfrak s}{\mathfrak l}_{2n}$) 
\begin{equation*}
h:=\text{diag}(2,\underbrace{1,1,\cdots,1}_{2p-2},\underbrace{0,0,\cdots,0}_{2(q-p+1)},\underbrace{-1,-1,\cdots,-1}_{2p-2},-2)
\end{equation*} 
which we will simply write 
\begin{equation*}
h=(2,\underbrace{1,1,\cdots,1}_{2p-2},\underbrace{0,0,\cdots,0}_{q-p+1})
\end{equation*} 
along with the weighted Dynkin diagram
\vskip 1mm
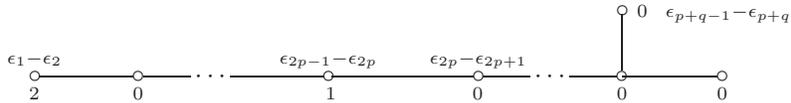
\[
\dynkin \underset{\;\;\;\; 2}{\whiteroot{\epsilon_1-\epsilon_2}}\hspace*{-0.1mm}\link\link\underset{0}{\whiteroot{}}\link\cdots\link\link\hspace*{-7mm}\underset{\;\;\;\;\;\;\;\;\;1}{\whiteroot{\epsilon_{2p-1}-\epsilon_{2p}}}\hspace*{-0.5mm}\link\link\link\hspace*{-6.5mm}\underset{\;\;\;\;\;\;\;\;0}{\whiteroot{\epsilon_{2p}-\epsilon_{2p+1}}}\link\cdots\link\underset{0}{\whiteroot{}}\hspace*{-0.6mm}\rootup{0\;\;\;\epsilon_{p+q-1}-\epsilon_{p+q}}\link\link\underset{0}{\whiteroot{}}
\enddynkin
\]
\vskip 4mm
The nilpotent orbit ${\mathcal O}^{\mathbb C}$ possesses $3$ or $4$ real forms depending wether $q>p$ or not. Using a recipe analogous to the one used for type $B$, we get
\vskip 1mm
\begin{eqnarray*}
h_1^{I}&=&(2,\underbrace{1,1,\cdots,1}_{p-1},\underbrace{1,1,\cdots,1}_{p-1},\underbrace{0,0,\cdots,0}_{q-p+1})\\
&&\\
h_1^{II}&=&(2,\underbrace{1,1,\cdots,-1}_{p-1},\underbrace{1,1,\cdots,1}_{p-1},\underbrace{0,0,\cdots,0}_{q-p+1})\\
&&\\
h_2^{I}&=&(\underbrace{1,1,\cdots,1}_{p-1},0,2,\underbrace{1,1,\cdots,1}_{p-1},\underbrace{0,0,\cdots,0}_{q-p})\\
&&\\
h_2^{II}&=&(\underbrace{1,1,\cdots,1}_{p-1},0,2,\underbrace{1,1,\cdots,-1}_{p-1},\underbrace{0,0,\cdots,0}_{q-p})\text{ only if }q=p
\end{eqnarray*} 
\vskip 1mm
\noindent
As before, $h_i^{II}$ is obtained from $h_i^{I}$ by the outer automorphism:
\vskip 2mm
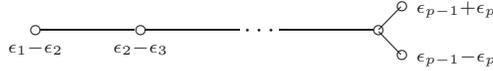
\[
\dynkin \underset{\epsilon_1-\epsilon_2}{\whiteroot{}}\hspace*{-3mm}\link\link\hspace*{-2.8mm}\underset{\epsilon_2-\epsilon_3}{\whiteroot{}}\hspace*{-3mm}\link\link\cdots\link\link\hspace*{0.7mm}\whiteroot{}\whiterootupright{\epsilon_{p-1}+\epsilon_{p}}
\whiterootdownright{\epsilon_{p-1}-\epsilon_{p}}
\enddynkin
\]
\vskip 3mm
\noindent
The description of the complex orbit ${\mathcal O}^{\mathbb C}$ and its real forms in terms of Young tableaux is as follows:
\vspace{0.5cm}\\
%\newpage
%\[
$\begin{array}{rl}
\hspace*{0.79cm}\young(\;\empty\;\empty)
\end{array}$
\\[-0.096cm]
$\begin{array}{rl}
{\mathcal O}^{\mathbb C}\hspace*{0.2cm}
\left.\young(\empty\;\empty,::,\empty\;\empty)\right\} & \hspace{-0.3cm}2p-2
\end{array}$
\\[-1.17mm]
\hspace*{0.63cm}
$\begin{array}{rl}
\left.\young(\empty,:,\empty)\right\} & \hspace{-0.3cm}2(q-p)+1\\[-0.8mm]
\end{array}$
%\]

\vspace*{-2.78cm}
%\vspace*{-1cm}
%\[
$\begin{array}{rl}
\hspace*{7.67cm}\young(+-+)
\end{array}$
\\[-0.096cm]
$\begin{array}{rl}
\hspace*{5cm}\text{ real forms for }{\mathcal O}^{\mathbb C}\hspace*{0.2cm}\left.\young(+-,::,+-)\right\} & \hspace{-0.3cm}2p-2
\hspace{0.5cm}\text{I, II}\\[-0.8mm]
\end{array}$
\\[-0.420mm]
\hspace*{7.935cm}
$\begin{array}{rl}
\left.\young(-,:,-)\right\} & \hspace{-0.3cm}2(q-p)+1
\end{array}$
%\]

\vspace*{1cm}
%\vspace*{-1cm}
%\[
$\begin{array}{rl}
\hspace*{7.67cm}\young(-+-)
\end{array}$
\\[-0.096cm]
$\begin{array}{rl}
\hspace*{7.85cm}\hspace*{0.2cm}\left.\young(+-,::,+-)\right\} & \hspace{-0.3cm}2p-2
\hspace{0.5cm}\text{(I, II if $q=p$)}\\[-0.8mm]
\end{array}$
\vspace*{0.02cm}
\hspace*{7.975cm}
$\begin{array}{rl}
\young(+)
\end{array}$
\vspace*{-0.11cm}\\
\hspace*{7.94cm}
$\begin{array}{rl}
\left.\young(-,:,-)\right\} & \hspace{-0.3cm}q-p
\end{array}$

\vspace*{1cm}
Consider $G_{\mathbb R}=SO^*(2n)=SO(2n,{\mathbb C})\cap{\mathfrak g}{\mathfrak l}(n,{\mathbb H})$, $K_{\mathbb R}=U(n)$ and 
${\mathfrak g}={\mathfrak s}{\mathfrak o}^*_{2n}$.
\vskip1mm
\noindent
For ${\mathfrak s}{\mathfrak o}_{2n}$, a Cartan subalgebra in ${\mathfrak g}$ consists of diagonal complex matrices of the form $\text{diag}(a_1,a_2,\cdots,a_n,-a_1,-a_2,\cdots,-a_n)$. Fix the standard system of positive roots $\{\epsilon_i\pm\epsilon_j\mid 1\leq i<j\leq n\}$. Using an recipe analogous to that of type $D$, one can attache an ${\mathfrak s}{\mathfrak l}_2$-triple to a complex nilpotent orbit. We apply this recipe to the nilpotent orbit ${\mathcal O}^{\mathbb C}={\mathcal O}_{\lbrack 2^{n}\rbrack}$ to obtain (viewing ${\mathfrak s}{\mathfrak o}^*_{2n}$ as a subalgebra of ${\mathfrak s}{\mathfrak l}_{2n}$) 
\begin{equation*}
h:=\text{diag}(\underbrace{1,1,\cdots,1}_{n},\underbrace{-1,-1,\cdots,-1}_{n})
\end{equation*} 
which we will simply write 
\begin{equation*}
h=(\underbrace{1,1,\cdots,1}_{n})
\end{equation*} 
along with the weighted Dynkin diagram
\vskip 1mm
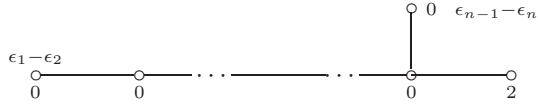
\[
\dynkin \underset{\;\;\;\; 0}{\whiteroot{\epsilon_1-\epsilon_2}}\hspace*{-0.1mm}\link\link\underset{0}{\whiteroot{}}\link\cdots\link\link\hspace*{-7mm}\link\cdots\link\underset{0}{\whiteroot{}}\hspace*{-0.6mm}\rootup{0\;\;\;\epsilon_{n-1}-\epsilon_{n}}\link\link\underset{2}{\whiteroot{}}
\enddynkin
\]
\vskip 4mm
The nilpotent orbit ${\mathcal O}^{\mathbb C}$ possesses $\frac{n}{2}+1$ real forms if $n$ is even, and $\frac{n+1}{2}$ real forms otherwise. Using a recipe analogous to the one used for type $D$, we get
\vskip 1mm
\begin{eqnarray*}
h_k&=&(\underbrace{1,1,\cdots,1}_{2k},\underbrace{-1,-1,\cdots,-1}_{n-2k})\text{ for $n$ even, and }k=0,\cdots,\frac{n}{2}\\
h_k&=&(\underbrace{1,1,\cdots,1}_{2k},0,\underbrace{-1,-1,\cdots,-1}_{n-2k-1})\text{ for $n$ odd, and }k=0,\cdots,\frac{n-1}{2}
\end{eqnarray*} 
\vskip 1mm
\noindent
Finally, the description of the complex orbit ${\mathcal O}^{\mathbb C}$ and its real forms in terms of Young tableaux is as follows:
\vspace{0.5cm}\\
\underline{$n$ even:}
\vspace{0.5cm}\\
%\[
$\begin{array}{rl}
{\mathcal O}^{\mathbb C}\hspace*{0.2cm}
\left.\young(\empty\;\empty,\empty\;\empty,::,\empty\;\empty,\empty\;\empty)\right\} & \hspace{-0.3cm}n
\end{array}$
%\]

\vspace*{-2cm}
%\[
$\begin{array}{rl}
\vspace*{-0.8cm}\hspace*{6cm}\text{ real forms }
\end{array}$
\vspace*{-0.4cm}
\hspace*{0.2cm}
$\begin{array}{rl}
\hspace*{-0.3cm}\left.\young(+-,::,+-)\right\} & \hspace{-0.34cm}k\in\{0,\ldots,\frac{n}{2}\} \\[-0.8mm]
\end{array}$
\hspace*{-3.94cm}
$\begin{array}{rl}
\vspace*{-2.45cm}\left.\young(-+,::,-+)\right\} & \hspace{-0.3cm}\frac{n}{2}-k
\end{array}$
%\]

\vspace*{2cm}
\underline{$n$ odd:}
\vspace*{0.5cm}\\
%\[
$\begin{array}{rl}
{\mathcal O}^{\mathbb C}\hspace*{0.2cm}
\left.\young(\empty\;\empty,::,\empty\;\empty)\right\} & \hspace{-0.3cm}n-1
\end{array}$
\vspace*{-0.515cm}\\
$\begin{array}{rl}
\hspace{1.455cm}\\\left.\young(\empty,\empty)\right\} & \hspace{-0.3cm}2\\[-0.8mm]
\end{array}$
%\]

\vspace*{-2cm}
%\[
$\begin{array}{rl}
\vspace*{0cm}\hspace*{4cm}\text{ real forms }
\end{array}$
\vspace*{-0.2cm}
\hspace*{0.1cm}
$\begin{array}{rl}
\hspace*{-0.1cm}\left.\young(+-,::,+-)\right\} & \hspace{-0.34cm}k\in\{0,\ldots,\frac{n-1}{2}\} \\[-0.8mm]
\end{array}$
\hspace*{-4.3cm}
$\begin{array}{rl}
\vspace*{-2.45cm}\left.\young(-+,::,-+)\right\} & \hspace{-0.3cm}\frac{n-1}{2}-k
\end{array}$
\vspace*{1.39cm}\\
$\begin{array}{rl}
\hspace*{6.815cm}\young(+)
\end{array}$
%\]

\section{Some general facts}

\bigskip

Let $\caO_i$ be a real form of the orbit $\caO^\bbC$, and denote the corresponding semisimple element of the (normal) $\frs\frl(2)$-triple by $h\in\frh$. As in \cite{MPVZ}, we attach to $h$ the $\theta$-stable parabolic subalgebra $\frqqq=\frl\oplus\fru$ such that $\frl$ is the centralizer of $h$ in $\frg$, and $\fru$ is the sum of negative eigenspaces for $\ad h$ on $\frg$. We fix a choice of $\Delta^+=\Delta^+(\frg,\frh)$; in examples, this will always be the standard positive root system. This defines a choice $\Delta^+_c=\Delta^+(\frk,\frh)=\Delta(\frk,\frh)\cap \Delta^+$ of positive compact roots. Let $\Delta^+_n:=\Delta^+\setminus\Delta^+_c$ be the set of positive noncompact roots. Denote by $\rho_c$ (resp. $\rho_n$) half the sum of positive compact (resp. noncompact) roots. The Weyl dimension polynomial
\eq
\label{def pk}
P_K(\lambda)=\prod_{\alpha\in\Delta_c^+}\frac{\langle\lambda,\alpha\rangle}{\langle\rho_c,\alpha\rangle}
\eeq
will always be defined with respect to this fixed positive root system $\Delta^+_c$. If $A$ is a subset of $\Delta(\frg,\frh)$, write $\rho(A)$ for half the sum of the roots in $A$.

We choose $\Delta^+(\frl)$ compatibly with $\Delta^+$, i.e., $\Delta^+(\frl)=\Delta(\frl)\cap\Delta^+$. This also gives a choice for positive roots of $\frl\cap\frk$, and fixes the Weyl dimension polynomial $P_{L\cap K}$. Denote by $\Delta^+_n(\frl)$ the set of noncompact roots in $\Delta^+(\frl)$, and by $\Delta(\frp_1)$ the set of noncompact roots that are $1$ on $h$.

The constant $c=c_i$ we are going to compute is 
attached to $\caO_i$ as in \cite{MPVZ}. 
It is defined by equation (6.4) of \cite{MPVZ}; this is 
up to sign the same equation as (5.9) of \cite{MPVZ}, but the sign is made precise using \cite{MPVZ}, Remark 3.8, equation (6.1) and the discussion around (6.1).
The equation is
\eq
\label{def const orig}
(-1)^N\sum_{\genfrac{}{}{0pt}{2}
{A\subseteq\Delta_n^+(\frl)}{C\subseteq\Delta(\frp_1)}} (-1)^{\#A+\#C}P_K(\lambda-\rho_n(\frl)+2\rho(A)-2\rho(C))=c P_{L\cap K}(\lambda),
\eeq
where
\eq
\label{def N}
N=\#\{\alpha\in\Delta^+\,\big|\,\alpha(h)>0\}.
\eeq

The computations we are going to make will be easier if  equation \eqref{def const orig} is turned into an analogue of equation (6.5) of \cite{MPVZ}:

\begin{proposition}
\label{eq const}
Assume that $\rho_n(\frl)$ is orthogonal to all roots of $\frl\cap\frk$. Then
\eq
\label{def const}
(-1)^{N+\#\Delta^+_n(\frl)}\sum_{\genfrac{}{}{0pt}{2}
{A\subseteq\Delta_n^+(\frl)}{C\subseteq\Delta(\frp_1)}} (-1)^{\#A+\#C}P_K(\lambda-2\rho(A)-2\rho(C))=c P_{L\cap K}(\lambda).
\eeq
\end{proposition}
\pf
This follows by passing from summation over $A$ to summation over the complement of $A$ in $\Delta_n^+(\frl)$. For any $A\subseteq\Delta_n^+(\frl)$, 
\[
-\rho_n(\frl)+2\rho(A)=\rho_n(\frl)-2\rho(\Delta_n^+(\frl)\setminus A)
\]
and
\[
(-1)^{\# A}=(-1)^{\#\Delta_n^+(\frl)}(-1)^{\#(\Delta_n^+(\frl)\setminus A)},
\]
so \eqref{def const orig} can be rewritten as
\[
(-1)^{N+\#\Delta^+_n(\frl)}\sum_{\genfrac{}{}{0pt}{2}
{A\subseteq\Delta_n^+(\frl)}{C\subseteq\Delta(\frp_1)}} (-1)^{\#A+\#C}P_K(\lambda+\rho_n(\frl)-2\rho(A)-2\rho(C))=c P_{L\cap K}(\lambda).
\]
We now replace $\lambda$ by $\lambda-\rho_n(\frl)$; since $\rho_n(\frl)$ is orthogonal to the roots of $\frl\cap\frk$,  
$P_{L\cap K}(\lambda-\rho_n(\frl))=P_{L\cap K}(\lambda)$, and the statement follows. 
\epf

In each of the examples we will consider, one can check directly that indeed $\rho_n(\frl)$ is orthogonal to all roots of $\frl\cap\frk$, and hence we can compute the constant $c$ using \eqref{def const}. A little more systematic way of checking this assumption, which will be easy to apply in all cases we consider, is given by the following lemma.

\begin{lemma}
\label{rhon}
Suppose that all simple factors of $\frl_0$ are either compact or noncompact Hermitian. Assume also that 
$\Delta^+(\frl)$ (induced by $\Delta^+$) is Borel - de Siebenthal for each noncompact factor of $\frl$. Then $\rho_n(\frl)$ is orthogonal to all roots of $\frl\cap\frk$. 
\end{lemma}
\pf
Since $\rho_n$ of any compact factor is 0, it is enough to prove the statement for each of the noncompact factors. Denote by $\frd$ one of these factors, and let 
$\frd=\frc\oplus\frs$ be its Cartan decomposition (so $\frc=\frd\cap\frk$ and $\frs=\frd\cap\frp$). Since $\frd$ is Hermitian,  
\[
\frs=\frs^+\oplus\frs^-
\]
as a $\frc$-module. If $\Delta^+(\frd)$ is a Borel-de Siebenthal positive root system with respect to a compact Cartan subalgebra of $\frd$, then $\Delta^+_n(\frd)$ must be equal to $\Delta(\frs^+)$ or $\Delta(\frs^-)$, and we can assume $\Delta^+_n(\frd)=\Delta(\frs^+)$. It follows that $\rho_n(\frd)$ is the weight of the one-dimensional $\frc$-module 
$\bigwedge^{\top}\frs^+$, and so it must be orthogonal to the roots of $\frc$.
\epf

\begin{remark}
\label{rmk nonherm}
{\rm
If $\frl$ has a simple noncompact factor that is not Hermitian, and if $\Delta^+(\frl)$ is any positive root system for $\frl$, then $\rho_n(\frl)$ is not orthogonal to all roots of $\frl\cap\frk$. Indeed, if $\frd=\frc\oplus\frs$ is the Cartan decomposition of one such factor, then $\frc$ is semisimple and hence has no nontrivial one-dimensional modules. So if $\rho_n(\frd)$ were orthogonal to all roots of $\frc$, it would have to be 0, but that is not possible since $\frd$ is noncompact.
}
\end{remark}

The following proposition will enable us to get our constants for some of the real forms of $\caO^\bbC$ without having to do computations. 

\begin{proposition}
\label{auto} 
Let $h_1$ and $h_2$ correspond to two real forms of $\caO^\bbC$. Assume that there is an automorphism $\sigma$ of $\frg$ such that
\begin{enumerate}
\item $\sigma$ preserves the compact Cartan subalgebra $\frh$ of $\frg$;
\item $\sigma$ commutes with the Cartan involution, so it preserves $\frk$ and $\frp$;
\item $\sigma(\Delta_c^+)=\Delta_c^+$;
\item $\sigma(h_1)=h_2$.
\end{enumerate}
Then the constants $c_1,c_2$ corresponding to $h_1,h_2$ are related by
\[
c_2=(-1)^{n+N_1+N_2}\,c_1, 
\]
where
\[
n=\#[\Delta_n^+(\frl_2)\cap(-\sigma(\Delta_n^+(\frl_1)))]=\#\{\alpha\in\Delta_n^+(\frl_1)\,\big|\,\sigma\alpha\in(-\Delta^+)\},
\]
and $N_1,N_2$ are defined as in \eqref{def N}, i.e., 
\[
N_i=\#\{\alpha\in\Delta^+\,\big|\,\alpha(h_i)>0\},\qquad i=1,2.
\]
\end{proposition}
\pf
If $\frl_i$ denotes the centralizer of $h_i$ in $\frg$, then it is clear that $\sigma(\frl_1)=\frl_2$. Moreover, the conditions we put on $\sigma$ ensure that $P_K\circ\sigma=P_K$, and also $P_{L_2\cap K}\circ\sigma=P_{L_1\cap K}$. 
 
We let $\sigma$ act on roots by
\[
\sigma(\alpha)=\alpha\circ\sigma^{-1}.
\]
Then it is clear that $\sigma$ takes $\Delta(\frp_1)_1$ (the set of noncompact roots that are 1 on $h_1$) to 
$\Delta(\frp_1)_2$ (the set of noncompact roots that are 1 on $h_2$). Furthermore, $\sigma$ takes 
$\Delta_n^+(\frl_1)$ to $\wti{\Delta_n^+}(\frl_2)$, where $\wti{\Delta_n^+}(\frl_2)$ is a positive root system for $\frl_2$, possibly different from $\Delta_n^+(\frl_2)$ which is defined using $\Delta^+$. It follows that $\sigma(\rho_n(\frl_1))=\tilde\rho_n(\frl_2)$, where $\tilde\rho_n(\frl_2)$ is the half sum of roots in $\wti{\Delta_n^+}(\frl_2)$. 

Also, for any $A\subseteq\Delta_n^+(\frl_1)$, $C\subseteq\Delta(\frp_1)_1$, we clearly have 
\[
2\rho(\sigma(A))=\sigma(2\rho(A)),\qquad 2\rho(\sigma(C))=\sigma(2\rho(C)).
\]
Writing the equation \eqref{def const orig} for $c_1$, we get
\[
(-1)^{N_1}\sum_{\genfrac{}{}{0pt}{2}
{A\subseteq\Delta_n^+(\frl_1)}{C\subseteq\Delta(\frp_1)_1}} (-1)^{\#A+\#C}P_K(\lambda-\rho_n(\frl_1)+2\rho(A)-2\rho(C))=c_1 P_{L_1\cap K}(\lambda).
\]
We now replace $\lambda$ by $\sigma^{-1}(\lambda)$ and use the equalities $P_{L_1\cap K}\circ\sigma^{-1}=P_{L_2\cap K}$, $P_K\circ\sigma^{-1}=P_K$. We also replace summing over $A$ and $C$ by summing over $\sigma(A)$ and $\sigma(C)$. We obtain
\eq
\label{auto calc}
(-1)^{N_1}\sum_{\genfrac{}{}{0pt}{2}
{\sigma(A)\subseteq\wti{\Delta_n^+}(\frl_2)}{\sigma(C)\subseteq\Delta(\frp_1)_2}} (-1)^{\#\sigma(A)+\#\sigma(C)}P_K(\lambda-\tilde\rho_n(\frl_2)+2\rho(\sigma(A))-2\rho(\sigma(C)))=c_1 P_{L_2\cap K}(\lambda).
\eeq
We now want to pass from summing over $\sigma(A)\subseteq\wti{\Delta_n^+}(\frl_2)$ to summing over $A'\subseteq\Delta_n^+(\frl_2)$. 
To do this, we define
\[
\Delta_1=\Delta^+_n(\frl_2)\cap {\wti{\Delta^+_n}(\frl_2)};\qquad \Delta_2=\Delta^+_n(\frl_2)\setminus\Delta_1=\Delta^+_n(\frl_2)\cap(-\wti{\Delta^+_n}(\frl_2));
\]
so 
\[
\Delta^+_n(\frl_2)=\Delta_1\cup\Delta_2;\qquad {\wti{\Delta^+_n}(\frl_2)}=\Delta_1\cup(-\Delta_2).
\]
It follows that for the half sums of roots $\rho_n(\frl_2),\tilde\rho_n(\frl_2)$ we have
\eq
\label{rho pm}
\tilde\rho_n(\frl_2)=\rho_n(\frl_2)-2\rho(\Delta_2).
\eeq
For any $A'\subseteq\Delta_n^+(\frl_2)$, let
\[
A'_1=A'\cap\Delta_1;\qquad A'_2=A'\cap \Delta_2;
\]
so $A'=A'_1\cup A'_2$. To each such $A'$ we attach
\[
\tilde A=A'_1 \cup (-(\Delta_2\setminus A'_2))\subseteq \wti{\Delta^+_n}(\frl_2).
\]
Then the correspondence $A'\leftrightarrow\tilde A$ defines a bijection between the subsets of $\Delta_n^+(\frl_2)$ and the subsets of $\wti{\Delta_n^+}(\frl_2)$. 
Using \eqref{rho pm}, we see that
\begin{align}
\label{2rhoA change}
2\rho(\tilde A)-\tilde\rho_n(\frl_2)=2\rho(A'_1)-(2\rho(\Delta_2)-2\rho(A'_2))-(\rho_n(\frl_2)-2\rho(\Delta_2))=\\ \nonumber 2\rho(A'_1)+2\rho(A'_2)-\rho_n(\frl_2)=2\rho(A')-\rho_n(\frl_2).
\end{align}
It follows that we can rewrite \eqref{auto calc} into a sum over $A'$ instead of a sum over $\tilde A=\sigma(A)$, taking into account that 
\[
(-1)^{\# \tilde A}=(-1)^{n}(-1)^{\# A'},
\]
with $n$ as in the statement of the proposition. If we also rename $A'$ by $A$ and $\sigma(C)$ by $C$, we get
\[
(-1)^{N_1+n}\sum_{\genfrac{}{}{0pt}{2}
{A\subseteq\Delta_n^+(\frl_2)}{C\subseteq\Delta(\frp_1)_2}} (-1)^{\# A+\# C}P_K(\lambda-\rho_n(\frl_2)+2\rho(A)-2\rho(C))=c_1 P_{L\cap K}(\lambda).
\]
If we compare this with the equation \eqref{def const orig} written for $c_2$, we immediately get the statement of the proposition.
\epf

\section{The case $G_\bbR=SU(p,q)$, $p\leq q$} 
\label{sec upq}

This case was treated in \cite{MPVZ} and we just record the results here.
The real forms of $\caO^\bbC$ correspond to
\eq
\label{h upq}
h_k=(\unb{k}{1,\dots,1},\,\unb{p-k}{-1,\dots,-1}\,|\,\unb{p-k}{1,\dots,1},\,\unb{q-p}{0,\dots,0},\,\unb{k}{-1,\dots,-1}),
\eeq
with $k=0,1,\dots,p$.
The corresponding constants $c=c_k^{p,q}$ can be computed from the formula \eqref{def const}.

The set $\Delta_n^+(\frl)$ is 
\begin{multline*}
\Delta^+_n(\frl)=\{\eps_i-\eps_j\,\big|\, 1\leq i\leq k,\, p+1\leq j\leq 2p-k\}\,\cup \\ 
\cup\,\{\eps_i-\eps_j\,\big|\, k+1\leq i\leq p,\, p+q-k+1\leq j\leq p+q\}.
\end{multline*}
The set $\Delta(\frp_1)$ is empty if $q=p$, and if $q>p$, then
\begin{multline*}
\Delta(\frp_1)=\{\eps_i-\eps_j\,\big|\, 1\leq i\leq k,\, 2p-k+1\leq j\leq p+q-k\}\,\cup \\ 
\cup\,\{\eps_i-\eps_j\,\big|\, 2p-k+1\leq i\leq p+q-k,\, k+1\leq j\leq p\}.
\end{multline*}
We evaluate \eqref{def const} at $\lambda=\lambda_0$, where
\[
\lambda_0 = (q,q-1,\dots,q-k+1,p,p-1,\dots,k+1 \,|\, p-k,\dots,1,q-k,\dots,p-k+1,k,\dots,1).
\]
For each choice of $A\subseteq\Delta_n^+(\frl)$ and $C\subseteq\Delta(\frp_1)$, we set 
\[
\Lambda=\lambda_0-2\rho(A)-2\rho(C).
\]
For $q>p$, we show that there is exactly one 
$C\subseteq\Delta(\frp_1)$ for which $P_K(\Lambda)$ can be nonzero: 
\[
C=\{\eps_i-\eps_j\,\big|\, 1\leq i\leq k; 2p-k+1\leq j\leq p+q-k\},
\]
with $\#C=k(q-p)$.

Then we show that all $\Lambda$ as above, with $P_K(\Lambda)\neq 0$, are of the form
\[
\Lambda=(i_1,\dots,i_k;j_1,\dots,j_{p-k}\,|\, j_1,\dots,j_{p-k}; q,\dots,p+1; i_1,\dots,i_k),
\]
with $i_1,\dots,i_k;j_1,\dots,j_{p-k}$ a shuffle of $p,\dots,1$, i.e., 
\[
i_1>\dots>i_k;\quad j_1>\dots>j_{p-k};\quad \{i_1,\dots,i_k;j_1,\dots,j_{p-k}\}=\{p,\dots,1\}.
\]
For each such $\Lambda$ there is a unique corresponding $A$, consisting of roots $\alpha_{a,b}$, $1\leq a\leq k$, $1\leq b\leq p-k$, where
\[
\alpha_{a,b}=\left\{\begin{matrix} \eps_a-\eps_{p+b}, &  \text{if } i_a<j_b;\cr 
                                                \eps_{k+a}-\eps_{p+q-k+b}, & \text{if } i_a>j_b.
                          \end{matrix}\right.
\]
In particular, for each $A$ involved, $\#A=k(p-k)$.

This leads to

\begin{theorem}
\label{thm A}
Let $G_\bbR=SU(p,q)$, and let $k\in\{0,1,\dots,p\}$ correspond to the real form of $\caO^\bbC$ given by (\ref{h upq}). Then
$c_k^{p,q}=(-1)^{k(p+q-k)}\binom{p}{k}$.
\end{theorem}

\bigskip

\section{The case $G_\bbR=SO_e(2p,2q+1)$, $q\geq p-1\geq 0$} 
\label{sec:B}

\bigskip

There are three real forms of $\caO^\bbC$ if $q\geq p\geq 1$, and two real forms if $q=p-1$. 

\subsection{The first real form}
\label{B1}  
 
This real form exists for all $q\geq p-1\geq 0$. 
The corresponding $h$ is
\[
h_1=(2,\unb{p-1}{1,\dots,1}\,|\,\unb{p-1}{1,\dots,1},\,\unb{q-p+1}{0,\dots,0}).
\]

Since $\frl=\frl_1$ is built from roots that vanish on $h_1$, we see that
\[
\Delta^+_n(\frl)=\{\eps_i-\eps_j\,\big|\, 2\leq i\leq p,\, p+1\leq j\leq 2p-1\}.
\]
It follows that for any $A\subseteq\Delta^+_n(\frl)$,
\eq
\label{2rhoA B1}
2\rho(A)=(0;a_1,\dots,a_{p-1}\,|\,-b_1,\dots,-b_{p-1};0,\dots,0),
\eeq
with 
\eq
\label{ineq ab B1}
0\leq a_i,b_j\leq p-1;\qquad \textstyle{\sum_i a_i=\sum_j b_j}.
\eeq
Furthermore, recall that $\Delta(\frp_1)$ consists of noncompact roots that are 1 on $h_1$. So 
\begin{multline*}
\Delta(\frp_1)=\{\eps_1-\eps_j\,\big|\, p+1\leq j\leq 2p-1\}\,\cup \\ 
\cup\,\{\eps_i\pm\eps_j\,\big|\, 2\leq i\leq p,\, 2p\leq j\leq p+q\}\cup \{\eps_i\,\big|\, 2\leq i\leq p\}.
\end{multline*}
It follows that for any $C\subseteq\Delta(\frp_1)$,
\eq
\label{2rhoC B1}
2\rho(C)=(c;d_1,\dots,d_{p-1}\,|\,-c_1,\dots,-c_{p-1};e_1,\dots,e_{q-p+1}),
\eeq
with 
\begin{align}
\label{ineq cde B1}
&0\leq c_j\leq 1;\qquad 0\leq c\leq p-1;\qquad \textstyle{c=\sum_j c_j;} \\ 
&0\leq d_i\leq 2(q-p+1)+1;\qquad -(p-1)\leq e_j\leq p-1. \nonumber
\end{align}
Note that for $q=p-1$, there are no coordinates after $2p-1$, so there are no zeros at the end of $2\rho(A)$, and there are no $e_j$. Otherwise, all of the above holds in this special case.

By \eqref{2rhoA B1},
\[
\rho_n(\frl)=(0,p-1,\dots,p-1\,|\,-p+1,\dots,-p+1,0,\dots,0).
\]
This is clearly orthogonal to all roots of $\frl\cap\frk$, which are equal to
\begin{align}
\label{roots lk B1}
\Delta(\frl\cap\frk)=\{\eps_i-\eps_j\,\big|\,2\leq i,j\leq p\}\cup\{\eps_i-\eps_j\,\big|\,p+1\leq i,j\leq 2p-1\}\cup\\
\nonumber \cup\{\eps_i\pm\eps_j\,\big|\,2p\leq i,j\leq p+q\}.
\end{align}
By Proposition \ref{eq const}, this means that we can determine the constant $c=c_1^{p,q}$ from the equation \eqref{def const}. To do this, we take $\lambda=\lambda_0$, where
\begin{align}
\label{la0 B1}
&\lambda_0=(\half;q+\half,q-\half\dots,q-p+\frac{5}{2}\,|\, -1,-2,\dots,-(p-1);q-p+1,q-p,\dots,1)\\ \nonumber
&\qquad\qquad\qquad\qquad\qquad\qquad \qquad\qquad\qquad \qquad\qquad\qquad\qquad\qquad\qquad   \text{if }\ q\geq p\geq 2; \\
&\lambda_0=(\half\,|\, q,q-1,\dots,1) \nonumber
\qquad \text{if }\ p=1, q\geq 1;\\
&\lambda_0=(\half;p-\half,p-\frac{3}{2}\dots,\frac{3}{2}\,|\, -1,-2,\dots,-(p-1))\nonumber
\qquad \text{if }\ p\geq 2,q=p-1;\\
&\lambda_0=(\half\,|\,) \qquad \text{if }\ p=1,q=0. \nonumber
\end{align}
\begin{proposition}
\label{prop B1} 
Let $\Lambda=\lambda_0-2\rho(A)-2\rho(C)$, with $\lambda_0$ given by \eqref{la0 B1}, and with $A\subseteq\Delta_n^+(\frl)$ and $C\subseteq\Delta(\frp_1)$. If $P_K(\Lambda)\neq 0$, then:
\begin{enumerate}
\item If $p=1$, then $A=C=\emptyset$ and $\Lambda=\lambda_0$.
\item If $p\geq 2$ and $q=p-1$, then $A=C=\emptyset$ and $\Lambda=\lambda_0$.
\item If $q\geq p\geq 2$, then
\begin{align*}
&A=\emptyset; \\
&C=\{\eps_i-\eps_j\,\big|\, 2\leq i\leq p,\, 2p\leq j\leq p+q\}; \\
&\Lambda=(\half;p-\half,p-\frac{3}{2},\dots,\frac{3}{2}\,|\, -1,-2,\dots,-(p-1);q,q-1,\dots,p).
\end{align*}
\end{enumerate}
\end{proposition}
\pf
We first note that if $p=1$, then 
\[
h_1=(2\,|\,0,\dots,0)\ \text{ or }\ h_1=(2),
\]
so $\Delta^+_n(\frl)=\Delta(\frp_1)=\emptyset$ and both $A$ and $C$ are automatically empty. It follows that the only possible $\Lambda$ is $\Lambda=\lambda_0$, and this proves the proposition for $p=1$. We continue by induction on $p$. Let us assume that $p\geq 2$, that $q\geq p-1$ is arbitrary, and that the statement of the proposition is true for $G_\bbR=SO_e(2p-2,2p-3)$, i.e., when $p,q$ are replaced by $p'=p-1$, and $q'=p-2$.

By \eqref{la0 B1}, \eqref{2rhoA B1} and \eqref{2rhoC B1}, we have
\begin{align}
\label{La B1}
\textstyle{\Lambda=(\half-c;q+\half-a_1-d_1,q-\half-a_2-d_2,\dots,q-p+\frac{5}{2}-a_{p-1}-d_{p-1}\,|}\\ 
|\,-1+b_1+c_1,-2+b_2+c_2,\dots,-(p-1)+b_{p-1}+c_{p-1}; \nonumber  \\ q-p+1-e_1,q-p-e_2,\dots,1-e_{q-p+1}). \nonumber 
\end{align}
(The third row of the above equation is not there if $q=p-1$.)

Using \eqref{ineq ab B1} and \eqref{ineq cde B1}, we see that the coordinates $\Lambda_{p+1},\dots,\Lambda_{p+q}$ are in the following intervals:
\begin{align*}
\Lambda=(\dots|\,\unb{[-1,p-1]}{-1+b_1+c_1},\unb{[-2,p-2]}{-2+b_2+c_2},\dots,\unb{[-(p-1),1]}{-(p-1)+b_{p-1}+c_{p-1}}; \\ \unb{[q-2p+2,q]}{q-p+1-e_1},\unb{[q-2p+1,q-1]}{q-p-e_2},\dots,\unb{[-(p-2),p]}{1-e_{q-p+1}}).  
\end{align*}
(The second row of the above equation is not there if $q=p-1$.)

So $\Lambda_{p+1},\dots,\Lambda_{p+q}$ are $q$ integers between $-(p-1)$ and $q$. Moreover, $P_K(\Lambda)\neq 0$ implies that these integers are nonzero, different from each other, and no two of them are opposite integers. If $q\geq p$, it follows that $q,q-1,\dots,p$ must each be equal to some $\Lambda_i$, and the only possibility for that is
\eqn
\Lambda_{2p}=q,\ \Lambda_{2p+1}=q-1,\ \dots,\ \Lambda_{p+q}=p.
\eeqn
So $e_1,\dots,e_{q-p+1}$ are all equal to $-(p-1)$, and hence
\eqn
\eps_i-\eps_j\in C,\ \eps_i+\eps_j\notin C,\qquad 2\leq i\leq p,\ 2p\leq j\leq p+q.
\eeqn
(If $q=p-1$, the above says nothing and should be skipped.)

This implies
\eq
\label{ineq d B1}
q-p+1\leq d_i\leq q-p+2,\qquad 1\leq i\leq p-1,
\eeq
with $d_i=q-p+1$ if $\eps_{i+1}\notin C$, and $d_i=q-p+2$ if $\eps_{i+1}\in C$. 
(If $q=p-1$, this gives no new information about the $d_i$.
The following arguments all work also in case $q=p-1$ if we delete the last group of coordinates, $q,q-1,\dots,p$.) 

Using \eqref{ineq d B1} together with the inequalities \eqref{ineq ab B1}, \eqref{ineq cde B1} for $a_i$ and $c$, we go back to \eqref{La B1} and conclude that 
$\Lambda_1,\dots,\Lambda_p$ are in the following intervals:
\begin{align*}
\Lambda=(\unb{[-p+\frac{3}{2},\half]}{\textstyle{\half-c}};\unb{[-\half,p-\half]}{\textstyle{q+\half-a_1-d_1}},\unb{[-\frac{3}{2},p-\frac{3}{2}]}{\textstyle{q-\half-a_2-d_2}},\dots,\unb{[-p+\frac{3}{2},\frac{3}{2}]}{\textstyle{q-p+\frac{5}{2}-a_{p-1}-d_{p-1}}}\,|\dots)
\end{align*}
So $\Lambda_1,\dots,\Lambda_p$ are $p$ half-integers between $-p+\frac{3}{2}$ and $p-\half$. Moreover, $P_K(\Lambda)\neq 0$ implies that these half-integers are different from each other, and no two of them are opposite. It follows that one of them must be equal to $p-\half$, and the only possibility is
\eqn
\Lambda_2=p-\half.
\eeqn
So $a_1=0$ and $d_1=q-p+1$. It follows that 
\begin{align*}
&\eps_2-\eps_j\notin A,\qquad p+1\leq j\leq 2p-1; \\ 
&\eps_2\notin C,
\end{align*}
and hence 
\eqn
0\leq b_j\leq p-2,\qquad 1\leq j\leq p-1.
\eeqn
These improved inequalities for the $b_j$ together with inequalities \eqref{ineq cde B1} for the $c_j$ imply 
\begin{align*}
\textstyle{\Lambda=(\dots|\,\unb{[-1,p-2]}{-1+b_1+c_1},\unb{[-2,p-3]}{-2+b_2+c_2},\dots,\unb{[-(p-1),0]}{-(p-1)+b_{p-1}+c_{p-1}};q,q-1,\dots,p).}  
\end{align*}
Since $\Lambda_{p+1},\dots,\Lambda_{2p-1}$ are $p-1$ nonzero integers between $-(p-1)$ and $p-2$, with no two of them equal or opposite to each other, we conclude that 
\eqn
\Lambda_{2p-1}=-(p-1).
\eeqn
This implies
\[
b_{p-1}=c_{p-1}=0,
\]
and hence
\begin{align*}
&\eps_i-\eps_{2p-1}\notin A,\quad 2\leq i\leq p;\\ 
&\eps_1-\eps_{2p-1}\notin C,
\end{align*}
and
\begin{align*}
&0\leq a_i\leq p-2,\quad 2\leq i\leq p-1;\\
&0\leq c\leq p-2.
\end{align*}
We see that
\begin{align*}
\textstyle{\Lambda=(\half-c;p-\half,q-\half-a_2-d_2,\dots,q-p+\frac{5}{2}-a_{p-1}-d_{p-1}\,|}\\ 
|\,-1+b_1+c_1,\dots,-(p-2)+b_{p-2}+c_{p-2},-(p-1);  \\ 
q,q-1,\dots,p). 
\end{align*}
(The third row is not there if $q=p-1$.)

We now consider the subalgebra $\frg'\cong\frs\fro(2p-2,2p-3)$ of $\frg$ built on coordinates 
\[
\eps_1,\eps_3,\dots,\eps_p;\eps_{p+1},\dots,\eps_{2p-2},
\]
so the coordinates 2 and $2p-1,2p,\dots,p+q$ are deleted. We also consider the real form of $\caO_{K'}$ given by 
\[
h'_1=(2,\unb{p-2}{1,\dots,1}\,|\,\unb{p-2}{1,\dots,1}),
\]
with centralizer $\frl'=\frl\cap\frg'$.
Then
\begin{align*}
&\Delta^+_n(\frl')=\Delta^+_n(\frl)\setminus\{\eps_i-\eps_j\,\big|\, i=2\text{ or } j=2p-1\};\\
&\Delta(\frp'_1)=\{\eps_1-\eps_{p+1},\dots,\eps_1-\eps_{2p-2};\,\eps_3,\dots,\eps_p\}.
\end{align*}
We set
\begin{align*}
& A'=A\cap\Delta^+_n(\frl')=A;\\
&C'=C\cap\Delta(\frp'_1)=C\setminus\{\eps_i-\eps_j\,\big|\,2\leq i\leq p,\ 2p\leq j\leq p+q\}.
\end{align*}
Then 
\begin{align*}
& 2\rho(A')=(0;a_2,\dots,a_{p-1}\,|\,-b_1,\dots,-b_{p-2})=    (0;a_1',\dots,a_{p-2}'\,|\,-b_1',\dots,-b_{p-2}') ;\\
&2\rho(C')=(c;d_2-(q-p+1),\dots,d_{p-1}-(q-p+1)\,|\,-c_1,\dots,-c_{p-2})=\\
&\qquad\qquad\qquad\qquad\qquad\qquad\qquad\qquad\qquad\qquad(c';d_1',\dots,d_{p-1}'\,|\,-c_1',\dots,-c_{p-2}'),
\end{align*}
where we define
\[
a_i'=a_{i+1};\quad b_i'=b_i;\quad c_i'=c_i;\quad c'=c;\quad d_i'=d_{i+1}-(q-p+1).
\]
The numbers $a_i',b_i',c_i',c',d_i'$ satisfy analogues of \eqref{ineq ab B1} and \eqref{ineq cde B1}.
We define $\lambda_0'$ by \eqref{la0 B1}, but for $G_\bbR=SO_e(2p-2,2p-3)$, i.e.,
\[
\lambda'_0=(\half;p-\frac{3}{2},\dots,\frac{3}{2}\,|\, -1,-2,\dots,-(p-2)).
\]
Then $A'$, $C'$ and 
\[
\Lambda'=\lambda'_0-2\rho(A')-2\rho(C')
\]
satisfy all conditions of the proposition, but $p,q$ are reduced to $p'=p-1$, $q'=p-2$. Moreover, $P_K(\Lambda)\neq 0$ is equivalent to $P_{K'}(\Lambda')\neq 0$. Therefore the inductive assumption implies that $A'=C'=\emptyset$, and that
$\Lambda'=\lambda_0'$. This implies the statement of the proposition for $A$, $C$ and $\Lambda$.
\epf

To compute the constant $c_1^{p,q}$, we have to compute $P_{L\cap K}(\lambda_0)$ where $\lambda_0$ is given by \eqref{la0 B1}, and $P_K(\Lambda)$ for $\Lambda$ determined in Proposition \ref{prop B1}. The main ingredients for this computation are given in the following lemma.
\begin{lemma}
\label{PK}
(i) Let $P^1_p$ be the Weyl dimension formula polynomial for $\frs\fro(2p)$, $p\geq 1$, and let $\lambda_p=(p-\half,p-\frac{3}{2},\dots,\half)$. Then 
\[
P^1_p(\lambda_p)=2^{p-1}.
\]
(ii) Let $P^2_q$ be the Weyl dimension formula polynomial for $\frs\fro(2q+1)$, $q\geq 1$, and let $P^2_0$ be the constant polynomial $1$. Furthermore, let 
$\mu_q=(q,q-1,\dots,1)$ if $q\geq 1$, and $\mu_0=0$. Then 
\[
P^2_q(\mu_q)=2^q.
\]
\end{lemma}
\pf
(i) Let $n_p^1$ be the numerator of $P^1_p$; the denominator is then 
\[
d^1_p=n^1_p(\rho_{\frs\fro(2p)})=n^1_p(p-1,p-2,\dots,1,0).
\]
The factors of $n_p^1(\lambda_p)$ that correspond to the roots $\eps_i-\eps_j$ clearly cancel with the corresponding factors of $d^1_p$. Denoting by $m_p^1(\lambda_p)$ respectively $e_p^1$ the product
of factors of $n^1_p(\lambda_p)$ respectively $d_p^1$ corresponding to the roots $\eps_i+\eps_j$, we have
\begin{align*}
&m_p^1(\lambda_p)=(2p-2)(2p-3)\dots(p+1)p\,m_{p-1}^1(\lambda_{p-1});\\ 
&e_p^1=(2p-3)(2p-4)\dots p(p-1)\, e_{p-1}^1.
\end{align*}
It follows that
\[
P^1_p(\lambda_p)=\frac{m^1_p(\lambda_p)}{e^1_p}=\frac{(2p-2)m_{p-1}^1(\lambda_{p-1})}{(p-1)e_{p-1}^1}=2P^1_{p-1}(\lambda_{p-1}).
\]
Since $P_1^1$ is the constant polynomial 1, this proves (i).

(ii) There is nothing to prove for $q=0$, and it is obvious that
\[
P_1^2(\mu_1)=\frac{1}{1/2}=2.
\]
For $q\geq 2$, let $n_q^2$ be the numerator of $P^2_q$; the denominator is then 
\[
d^2_q=n^2_q(\rho_{\frs\fro(2q+1)})=n^2_q(q-\half,q-\frac{3}{2},\dots,\frac{3}{2},\half).
\]
The factors of $n_q^2(\mu_q)$ that correspond to the roots $\eps_i-\eps_j$ clearly cancel with the corresponding factors of $d^2_q$. Denoting by $m_q^2(\mu_q)$ respectively $e_q^2$ the product
of factors of $n^2_q(\mu_q)$ respectively $d_q^2$ corresponding to the roots $\eps_i+\eps_j$ and $\eps_i$, we have
\begin{align*}
&m_q^2(\mu_q)=(2q-1)(2q-2)\dots(q+2)(q+1)q\,m_{q-1}^2(\mu_{q-1});\\ 
&e_q^2=(2q-2)(2q-3)\dots (q+1)q(q-\half)\, e_{q-1}^2.
\end{align*}
It follows that
\[
P^2_q(\mu_q)=\frac{m^2_q(\mu_q)}{e^2_q}=\frac{(2q-1)m_{q-1}^2(\mu_{q-1})}{(p-\half)e_{q-1}^2}=2P^2_{q-1}(\mu_{q-1}).
\]
The statement now follows by induction.
\epf

To compute $P_{L\cap K}(\lambda_0)$, we recall \eqref{roots lk B1}, which shows that $\frl\cap\frk$ is up to center typically a product of three factors: the $\fru(p-1)$ on coordinates $2,\dots,p$, the $\fru(p-1)$ on coordinates $p+1,\dots,2p-1$, and the $\frs\fro(2(q-p+1)+1)$ on coordinates $2p,\dots,p+q$. (If $p=1$, then the first two factors are missing, if $q=p-1$ then the third factor is missing, and if $p=1$ and $q=0$, then $\frl\cap\frk$ is one-dimensional. What we say below applies also to these cases with obvious modifications.)

It is clear from the definition \eqref{la0 B1} of $\lambda_0$ that for each of the first two factors, the corresponding coordinates of $\lambda_0$ differ from the $\rho$ of the factor by a weight orthogonal to the roots of the factor, so in the notation of Lemma \ref{PK},
\eq
\label{PLK B1}
P_{L\cap K}(\lambda_0)=P^2_{q-p+1}(\mu_{q-p+1})=2^{q-p+1}.
\eeq
To compute $P_K(\Lambda)$, we first write $\Lambda=(\Lambda_L\,|\,\Lambda_R)$ and note that 
\[
P_K(\Lambda)=P^1_p(\Lambda_L)P^2_q(\Lambda_R).
\] 
To use Lemma \ref{PK}, we have to rearrange coordinates of $\Lambda_L$ and $\Lambda_R$, and use the fact that $P^1_p$ is skew for the Weyl group of $\frs\fro(2p)$, while $P_q^2$ is skew for the Weyl group of $\frs\fro(2q+1)$. To rearrange $\Lambda_L$ to $\lambda_p$, we only need to bring the $\half$ from the first coordinate to the $p$-th coordinate, and hence
\eq
\label{PK1 B1}
P_p^1(\Lambda_L)=(-1)^{p-1}2^{p-1}.
\eeq
To bring $\Lambda_R$ to $\mu_q=(q,\dots,1)$, we need to change $p-1$ signs, and then bring coordinates $p-1,p-2,\dots,1$, in that order, all the way to the right. The sign produced in this way is
\[
(-1)^{(p-1)+(q-p+1)+(q-p+2)+\dots+(q-1)}=(-1)^{(p-1)(q-p+1)+\frac{(p-1)p}{2}}.
\]
Since 
\eq
\label{congruence}
\frac{(p-1)p}{2}\equiv[\frac{p}{2}]\quad\text{mod } 2,
\eeq 
Lemma \ref{PK} implies that 
\eq
\label{PK2 B1}
P_q^2(\Lambda_R)=(-1)^{(p-1)(q-p+1)+[\frac{p}{2}]}2^q.
\eeq
Now we substitute \eqref{PLK B1}, \eqref{PK1 B1} and \eqref{PK2 B1} into \eqref{def const}. Since 
\[
\# A+\# C=\#C=(p-1)(q-p+1),
\]
and since $N$ from \eqref{def N} is easily checked to satisfy
\eq
\label{N B1}
N\equiv p\quad\text{mod }2,
\eeq
we see that the total sign is $(-1)^{[\frac{p}{2}]+1}$,
and we conclude
\begin{theorem}
\label{thm B1}
Let $G_\bbR=SO_e(2p,2q+1)$, $q\geq p-1\geq 0$, and let $c_1^{p,q}$  be the constant corresponding to the first real form of $\caO^\bbC$. Then
\[
c_1^{p,q}=(-1)^{[\frac{p}{2}]+1}2^{2p-2}.
\]
\qed
\end{theorem}

\subsection{The second real form}
\label{B2} 
This real form exists for all $q\geq p-1\geq 0$. 
The corresponding $h$ is 
\[
h_2=(2,\unb{p-2}{1,\dots,1},-1\,|\,\unb{p-1}{1,\dots,1},\unb{q-p+1}{0,\dots,0}).
\]
This real form is conjugate to the first real form by the automorphism $\sigma=s_{\eps_p}$, the reflection with respect to the short noncompact root $\eps_p$. The automorphism $\sigma$ of $\frg$ clearly satisfies the conditions of Proposition \ref{auto}. Moreover, the number $n$ from Proposition \ref{auto} is $p-1$; the roots from $\Delta_n^+(\frl_1)$ that $\sigma$ sends to $-\Delta^+$ are 
\[
\eps_p-\eps_j,\qquad p+1\leq j\leq 2p-1.
\]
Moreover, by \eqref{N B1}, $N_1\equiv p\quad\text{mod }2$, and another short computation shows that $N_2$ is always even. The total sign in Proposition \ref{auto} is thus 
\[
(-1)^{n+N_1+N_2}=-1,
\]
so 
\eq
\label{const B2}
c_2^{p,q}=-c_1^{p,q}=(-1)^{[\frac{p}{2}]}2^{2p-2}.
\eeq

\subsection{The third real form}
\label{B3} 
This real form exists for $q\geq p\geq 1$, so we assume this condition in the following. 
The corresponding $h$ is 
\[
h_3=(\unb{p-1}{1,\dots,1},0\,|\,2,\unb{p-1}{1,\dots,1},\unb{q-p}{0,\dots,0}).
\]

Since $\frl=\frl_3$ is built from roots that vanish on $h_3$, we see that
\[
\Delta^+_n(\frl)=\{\eps_i-\eps_j\,\big|\, 1\leq i\leq p-1,\, p+2\leq j\leq 2p\}\cup\{\eps_p\pm\eps_j\,\big|\,2p+1\leq j\leq p+q\}\cup\{\eps_p\}.
\]
It follows that for any $A\subseteq\Delta^+_n(\frl)$,
\eq
\label{2rhoA B3}
2\rho(A)=(a_1,\dots,a_{p-1};x\,|\,0;-b_1,\dots,-b_{p-1};y_1,\dots,y_{q-p}),
\eeq
with 
\begin{align}
\label{ineq abxy B3}
&0\leq a_i,b_j\leq p-1;\qquad \textstyle{\sum_i a_i=\sum_j b_j}; \\
&\nonumber 0\leq x\leq 2(q-p)+1;\qquad -1\leq y_j\leq 1.
\end{align}
Furthermore, recall that $\Delta(\frp_1)$ consists of noncompact roots that are 1 on $h_3$. So 
\begin{multline*}
\Delta(\frp_1)= 
\{\eps_i\pm\eps_j\,\big|\, 1\leq i\leq p-1,\, 2p+1\leq j\leq p+q\}\cup \{\eps_i\,\big|\, 1\leq i\leq p-1\}\,\cup \\
\cup\, \{\eps_j\pm\eps_p\,\big|\, p+2\leq j\leq 2p\}\,\cup\,\{\eps_{p+1}-\eps_i\,\big|\,1\leq i\leq p-1\}.
\end{multline*}
It follows that for any $C\subseteq\Delta(\frp_1)$,
\eq
\label{2rhoC B3}
2\rho(C)=(c_1,\dots,c_{p-1};u\,|\,v;d_1,\dots,d_{p-1};e_1,\dots,e_{q-p}),
\eeq
with 
\begin{align}
\label{ineq cdeuv B3}
&-1\leq c_i\leq 2(q-p)+1;\qquad -(p-1)\leq u\leq p-1; \\ 
&0\leq v\leq p-1;\qquad 0\leq d_j\leq 2;\qquad -(p-1)\leq e_j\leq p-1. \nonumber
\end{align}
If we write \eqref{2rhoA B3} for $A=\Delta^+_n(\frl)$, we get
\[
\rho_n(\frl)=(p-1,\dots,p-1;2(q-p)+1\,|\, 0;-p+1,\dots,-p+1;0,\dots,0).
\]
This is clearly orthogonal to all roots of $\frl\cap\frk$, which are equal to
\begin{align}
\label{roots lk B3}
\Delta(\frl\cap\frk)=\{\eps_i-\eps_j\,\big|\,1\leq i,j\leq p-1\}\cup\{\eps_i-\eps_j\,\big|\,p+2\leq i,j\leq 2p\}\cup\\
\nonumber \cup\{\eps_i\pm\eps_j\,\big|\,2p+1\leq i,j\leq p+q\}\cup\{\eps_{2p+1},\dots,\eps_{p+q}\}.
\end{align}

By Proposition \ref{eq const}, this means that the constant $c=c_3^{p,q}$ satisfies \eqref{def const} for any $\lambda$, and we will compute $c_3^{p,q}$ by using this for $\lambda=\lambda_0$, where
\eq
\label{la0 B3}
\lambda_0=(q-\frac{3}{2},q-\frac{5}{2}\dots,q-p+\frac{1}{2};q-p+\half\,|\,p-1;0,-1,\dots,-(p-2);q-p,q-p-1,\dots,1).
\eeq
(Out of the 5 groups of coordinates separated by semicolons and the bar, the first and the fourth group are missing if  $p=1$, and the fifth group is missing if $q=p$.)
 
\begin{proposition}
\label{prop B3} 
Let $\Lambda=\lambda_0-2\rho(A)-2\rho(C)$, with $\lambda_0$ given by \eqref{la0 B3}, and with $A\subseteq\Delta_n^+(\frl)$ and $C\subseteq\Delta(\frp_1)$. Then $\Lambda_{p+1}=0$. In particular, $P_K(\Lambda)=0$.
\end{proposition}
\pf
By \eqref{la0 B3}, \eqref{2rhoA B3} and \eqref{2rhoC B3}, we have
\begin{align*}
\textstyle{\Lambda=(q-\frac{3}{2}-a_1-c_1,\dots,q-p+\half-a_{p-1}-c_{p-1};q-p+\half-x-u\,|}\\ 
|\,p-1-v;b_1-d_1,-1+b_2-d_2,\dots,-(p-2)+b_{p-1}-d_{p-1}; \\ q-p-y_1-e_1,q-p-1-y_2-e_2,\dots,1-y_{q-p}-e_{q-p}).  
\end{align*}
If $p=1$, then the situation is much simpler; in particular, since the coordinates of $h_3$ are 0 or 2, $\Delta(\frp_1)$ is empty, so $C=\emptyset$ and $2\rho(C)=0$. It follows that
$\Lambda_{p+1}=\Lambda_2=p-1=0$, and so $P_K(\Lambda)=0$ as claimed. 

So the proposition is true for $p=1$. We continue by induction on $p$. Let us assume that $p\geq 2$, that $q\geq p$ is arbitrary, and that the statement of the proposition is true for $G_\bbR=SO_e(2p-2,2p-1)$, i.e., when $p,q$ are replaced by $p'=p-1$, $q'=p-1$.

By \eqref{la0 B3}, \eqref{ineq abxy B3} and \eqref{ineq cdeuv B3}, the coordinates $\Lambda_{p+1},\dots,\Lambda_{p+q}$ are in the following intervals:
\begin{align}
\label{segR B3}
\textstyle{\Lambda=(\dots|\,\unb{[0,p-1]}{p-1-v};\unb{[-2,p-1]}{b_1-d_1},\unb{[-3,p-2]}{-1+b_2-d_2},\dots,\unb{[-p,1]}{-(p-2)+b_{p-1}-d_{p-1}};}  \\ \unb{[q-2p,q]}{q-p-y_1-e_1},\unb{[q-2p-1,q-1]}{q-p-1-y_2-e_2},\dots,\unb{[-(p-1),p+1]}{1-y_{q-p}-e_{q-p}}).  \nonumber
\end{align}
So $\Lambda_{p+1},\dots,\Lambda_{p+q}$ are $q$ integers between $-p$ and $q$. Moreover, $P_K(\Lambda)\neq 0$ implies that these integers are nonzero, and no two of them are equal or opposite to each other. It follows that $q,q-1,\dots,p+1$ must each be equal to some $\Lambda_i$, and the only possibility for that is
\eqn
\Lambda_{2p+1}=q,\ \Lambda_{2p+2}=q-1,\ \dots,\ \Lambda_{p+q}=p+1.
\eeqn
It follows that $y_1,\dots,y_{q-p}$ are all equal to $-1$, and that $e_1,\dots,e_{q-p}$ are all equal to $-(p-1)$. 
So 
\begin{align*}
& \eps_p-\eps_j\in A,\ \eps_p+\eps_j\notin A,\qquad 2p+1\leq j\leq p+q;\\
& \eps_i-\eps_j\in C,\ \eps_i+\eps_j\notin C,\qquad 1\leq i\leq p-1,\ 2p+1\leq j\leq p+q.
\end{align*}
This implies
\begin{align}
\label{ineq xc B3}
& q-p\leq x\leq q-p+1; \\
&\nonumber q-p-1\leq c_i\leq q-p+1,\qquad 1\leq i\leq p-1.
\end{align}
Note that $x=q-p$ if $\eps_p\notin A$ and $x=q-p+1$ if $\eps_p\in A$. Similarly, $c_i=q-p-1$ if $\eps_{p+1}-\eps_i\in C$, $\eps_i\notin C$; $c_i=q-p$ if $\eps_{p+1}-\eps_i\in C$, $\eps_i\in C$ or $\eps_{p+1}-\eps_i\notin C$, $\eps_i\notin C$; and $c_i=q-p+1$ if $\eps_{p+1}-\eps_i\notin C$, $\eps_i\in C$.

Using the same arguments as above, we can also conclude from 
\eqref{segR B3} that
\eqn
\Lambda_{2p}=-p.
\eeqn
This implies that $b_{p-1}=0$ and $d_{p-1}=2$. It follows that
\begin{align*}
&\eps_i-\eps_{2p}\notin A,\quad 1\leq i\leq p-1;\\  
&\eps_{2p}\pm\eps_p\in C,
\end{align*}
so
\begin{align}
\label{ineq au B3}
&0\leq a_i\leq p-2,\quad 1\leq i\leq p-1;\\ \nonumber 
&-(p-2)\leq u\leq p-2.
\end{align}
Using the improved inequalities \eqref{ineq xc B3} and \eqref{ineq au B3}, we see that $\Lambda_1,\dots,\Lambda_p$ are in the following intervals:
\begin{align*}
\Lambda=(\unb{[-\half,p-\half]}{\textstyle{q-\frac{3}{2}-a_1-c_1}},\unb{[-\frac{3}{2},p-\frac{3}{2}]}{\textstyle{q-\frac{5}{2}-a_2-c_2}},\dots,\unb{[-(p-\frac{3}{2}),\frac{3}{2}]}{\textstyle{q-p+\half-a_{p-1}-c_{p-1}}}; \\
\unb{[-(p-\frac{3}{2}),p-\frac{3}{2}]}{\textstyle{q-p+\half-x-u}}\,|\, \dots \nonumber
\end{align*}
So $\Lambda_1,\dots,\Lambda_p$ are $p$ half-integers between $-(p-\frac{3}{2})$ and $p-\half$, such that no two of them are equal or opposite to each other. It follows that 
\eqn
\Lambda_1=p-\half,
\eeqn
and consequently $a_1=0$, $c_1=q-p+1$. Therefore,
\begin{align*}
&\eps_1-\eps_{j}\notin A,\quad p+2\leq j\leq 2p;\\  
&\eps_{p+1}-\eps_1\in C,\ \eps_1\notin C,
\end{align*}
and we conclude that
\begin{align*}
&0\leq b_j\leq p-2,\quad 1\leq j\leq p-1;\\  
&1\leq v\leq p-1.
\end{align*}
We see that
\begin{align*}
\textstyle{\Lambda=(p-\half,q-\frac{5}{2}-a_2-c_2,\dots,q-p+\frac{1}{2}-a_{p-1}-c_{p-1};q-p+\half-x-u\,|}\\ 
|\,p-1-v;b_1-d_1,\dots,-(p-3)+b_{p-2}-d_{p-2},-p;  \\ 
q,q-1,\dots,p+1). 
\end{align*}
(If $q=p$, the coordinates $q,\dots,p+1$ are not there; if $p=2$ there are no coordinates involving $a_i,c_i,b_i$ or $d_i$.)

We now consider the subalgebra $\frg'\cong\frs\fro(2p-2,2p-1)$ of $\frg$ built on coordinates 
\[
\eps_2,\eps_3,\dots,\eps_p;\eps_{p+1},\dots,\eps_{2p-1},
\]
so the coordinates 1 and $2p,2p+1,\dots,p+q$ are deleted. We also consider the real form of $\caO_{K'}$ given by 
\[
h'_3=(\unb{p-2}{1,\dots,1},0\,|\,2,\unb{p-2}{1,\dots,1}),
\]
with centralizer $\frl'=\frl\cap\frg'$. Then
\begin{align*}
&\Delta^+_n(\frl')=\{\eps_i-\eps_j\,\big|\,2\leq i\leq p-1,\,p+2\leq j\leq 2p-1\}\cup\{\eps_p\};\\
&\Delta(\frp'_1)=\{\eps_2,\dots,\eps_{p-1}\}\cup \{\eps_j\pm\eps_p\,\big|\,p+2\leq j\leq 2p-1\}\cup\\
&\qquad\qquad\qquad\qquad\qquad\qquad\qquad\qquad\cup\{\eps_{p+1}-\eps_i\,\big|\,2\leq i\leq p-1\}.
\end{align*}
We set
\begin{align*}
& A'=A\cap\Delta^+_n(\frl')=A\setminus\{\eps_p-\eps_j\,\big|\,2p+1\leq j\leq p+q\};\\
&C'=C\cap\Delta(\frp'_1)=\\
&\qquad\qquad C\setminus\{\eps_{2p}\pm\eps_p;\,\eps_{p+1}-\eps_1;\,\eps_i-\eps_j\,\big|\,1\leq i\leq p-1,\ 2p+1\leq j\leq p+q\}.
\end{align*}
Then 
\begin{align*}
& 2\rho(A')=(a_2,\dots,a_{p-1};x-(q-p)\,|\,0;-b_1,\dots,-b_{p-2})=  \\  
&\qquad\qquad\qquad\qquad(a_1',\dots,a_{p-2}';x'\,|\,0;-b_1',\dots,-b_{p-2}') ;\\
&2\rho(C')=(c_2-(q-p),\dots,c_{p-1}-(q-p);u\,|\,v-1;d_1,\dots,d_{p-2})=\\
&\qquad\qquad\qquad\qquad\qquad(c_1',\dots,c_{p-2}';u'\,|\,v';d_1',\dots,d_{p-2}'),
\end{align*}
where we define
\begin{align*}
&a_i'=a_{i+1};\quad x'=x-(q-p);\quad b_i'=b_i;\\ 
&c_i'=c_{i+1}-(q-p);\quad u'=u;\quad v'=v-1;\quad d_i'=d_{i}.
\end{align*}
The numbers $a_i',x',b_i',c_i',u',v',d_i'$ satisfy analogues of \eqref{ineq abxy B3} and \eqref{ineq cdeuv B3}.

We define $\lambda_0'$ by \eqref{la0 B3}, but for $G_\bbR=SO_e(2p-2,2p-1)$, i.e.,
\[
\lambda'_0=(p-\frac{5}{2},p-\frac{7}{2},\dots,\half;\half;\,|\,p-2;0,-1,\dots,-(p-3)).
\]
Then $A'$, $C'$ and 
\[
\Lambda'=\lambda'_0-2\rho(A')-2\rho(C')
\]
satisfy all conditions of the proposition, but $p,q$ are reduced to $p'=p-1$, $q'=p-1$. Therefore the inductive assumption implies that $\Lambda'_p=0$. So $v'=p-2$, and therefore $v=p-1$ and $\Lambda_{p+1}=0$. It follows that $P_K(\Lambda)=0$, since $\Lambda$ is orthogonal to the compact root $\eps_{p+1}$.
\epf

Proposition \ref{prop B3} implies that the left hand side of \eqref{def const} is 0 in this case. On the other hand, 
\[
P_{L\cap K}(\lambda_0)\neq 0
\]
by \eqref{la0 B3} and \eqref{roots lk B3}. We conclude

\begin{theorem}
\label{thm B3}
For $G_\bbR=SO_e(2p,2q+1)$, $q\geq p\geq 1$, the constant corresponding to the third real form is
\[
c_3^{p,q}=0.
\]
\end{theorem}

\section{The case $G_\bbR=\Sp(2n,\bbR)$, $n\leq 1$} 
\label{sec sp}

The real forms of $\caO^\bbC$ correspond to integers $p$ such that $0\leq p\leq n$. We denote $n-p$ by $q$. The $h$ corresponding to $p$ is 
\eqn
h_p=(\unb{p}{1,\dots,1},\,\unb{q}{-1,\dots,-1}),\qquad  p=0,1,\dots,n.
\eeqn

Since $\frl=\frl_p$ is built from roots that vanish on $h_p$, we see that
\[
\Delta^+_n(\frl)=\{\eps_i+\eps_{p+j}\,\big|\, 1\leq i\leq p,\, 1\leq j\leq q\}.
\]
It follows that for any $A\subseteq\Delta^+_n(\frl)$,
\eq
\label{2rhoA C}
2\rho(A)=(a_1,\dots,a_p\,|\,b_1,\dots,b_q),
\eeq
with 
\begin{align}
\label{ineq ab C}
0\leq a_i\leq q, \qquad 0\leq b_j\leq p,\qquad \textstyle{\sum_ia_i=\sum_jb_j}.
\end{align}
In particular,
\[
\rho_n(\frl)=(q,\dots,q\,|\,p,\dots,p),
\]
and this is clearly orthogonal to the roots of $\frl\cap\frk$, which are given by
\eqn
\Delta^+(\frl\cap\frk)=\{\eps_i-\eps_j\,\big|\,1\leq i<j\leq p\}\cup\{\eps_{p+i}-\eps_{p+j}\,\big|\,1\leq i<j\leq q\}.
\eeqn
So the constants $c=c^n_p$ can be calculated from \eqref{def const}. Since it is clear that in our present case 
\[
\Delta(\frp_1)=\emptyset,
\]
\eqref{def const} becomes
\eq
\label{def const C}
\sum_{A\subseteq\Delta_n^+(\frl)} (-1)^{\#A}P_K(\lambda-2\rho(A))=c P_{L\cap K}(\lambda).
\eeq
We take $\lambda=\lambda_0$, where
\eq
\label{la0 C}
\lambda_0=(n,n-1,\dots,q+1\,|\,n,n-1,\dots,p+1),
\eeq
(If $p$ is 0 or $n$, then there is only one group of coordinates in the above expression, and $\lambda_0=(n,n-1,\dots,1)$.)

Since $\lambda_0$ differs from $\rho_{\frl\cap\frk}$ by a weight orthogonal to all roots of $\frl\cap\frk$, 
\[
P_{L\cap K}(\lambda_0)=1.
\]
So to compute $c_p^n$ we have to compute the left side of \eqref{def const C}. The following proposition describes the relevant $A$ and the corresponding $\Lambda$.

\begin{proposition}
\label{prop C} 
Let $\Lambda=\lambda_0-2\rho(A)$, with $\lambda_0$ given by \eqref{la0 C}, and with $A\subseteq\Delta_n^+(\frl)$. 

(i) If $p$ and $q$ are both odd, then $P_K(\Lambda)=0$ for all $\Lambda$ as above.

(ii) Suppose that at least one of $p,q$ is even, and suppose that for some $A$ the corresponding $\Lambda$ satisfies 
$P_K(\Lambda)\neq 0$. Then:
\begin{enumerate}
\item If $p=0$ or $q=0$, then $A=\emptyset$ and $\Lambda=\lambda_0=(n,n-1,\dots,1)$.
\item If $0<p<n$, let $r=[\frac{p}{2}]$ and $s=[\frac{q}{2}]$. Then there is a shuffle
\[
1\leq i_1<\dots<i_r\leq r+s;\qquad 1\leq j_1<\dots<j_s\leq r+s
\]
of $1,2,\dots,r+s$ such that 
\[
A=\{\alpha_{u,v},\beta_{u,v}\,\big|\, 1\leq u\leq r,\,1\leq v\leq s\}\cup B,
\]
where 
\[
\alpha_{u,v}=\eps_{p+1-u}+\eps_{n+1-v};\qquad \beta_{u,v}=\left\{
\begin{matrix} 
\eps_{p+1-u}+\eps_{p+v},&\qquad i_u<j_v;\cr
\eps_u+\eps_{n+1-v},&\qquad i_u>j_v.
\end{matrix}
\right.
\]
and
\[
B=\left\{
\begin{matrix}
\emptyset, \qquad p,q \text{ even};\cr 
\{\eps_{r+1}+\eps_{p+j}\,|\, s+1\leq j\leq q\},\qquad p \text{ odd};\cr
\{\eps_i+\eps_{p+s+1}\,|\, r+1\leq i\leq p\},\qquad q \text{ odd}.
\end{matrix}
\right.
\]
The corresponding $\Lambda$ has coordinates   
\begin{align*}
&\Lambda_1=n+1-i_1,\dots,\Lambda_r=n+1-i_r;\qquad \Lambda_{p-r+1}=i_r,\dots,\Lambda_p=i_1;\\
&\Lambda_{p+1}=n+1-j_1,\dots,\Lambda_{p+s}=n+1-j_s;\qquad \Lambda_{n-s+1}=j_s,\dots,\Lambda_n=j_1,
\end{align*}
and possibly in addition
\begin{align*}
\Lambda_{p-r}=n-r-s, \quad \text{if } p\text{ is odd};\qquad
\Lambda_{n-s}=n-r-s, \quad \text{if } q\text{ is odd}.
\end{align*}
\end{enumerate}
\end{proposition}
\pf
The statement is obviously true for any $n$ if $p=0$ or $q=0$. Hence it is true for $n=1$. If $n=2$ and $p=q=1$, there are two cases: 
\[
A=\emptyset,\qquad \text{or}\qquad A=\{\eps_1+\eps_2\}.
\]
If $A=\emptyset$, then $\Lambda=\Lambda_0=(n\,|\,n)$, so
$P_K(\Lamdba)=0$. If $A=\{\eps_1+\eps_2\}$, then 
$\Lambda=(n-1\,|\,n-1)$, and again $P_K(\Lambda)=0$. So the proposition is true for $n=2$.

We proceed by induction on $n$. Assume that $n>2$ and $p,q\geq 1$, and assume that the proposition is true for $n-2$.

Using the definitions and the inequalities \eqref{ineq ab C},
we see that
\begin{align*}
\Lambda=(\unb{[p,n]}{n-a_1},\unb{[p-1,n-1]}{n-1-a_2},\dots,\unb{[1,q+1]}{q+1-a_p}\,|\,\unb{[q,n]}{n-b_1},\unb{[q-1,n-1]}{n-1-b_2},\dots,\\
\unb{[1,p+1]}{p+1-b_q}).
\end{align*}
So the coordinates of $\Lambda$ are $n$ integers between $1$ and $n$, and assuming that $P_K(\Lamdba)\neq 0$, they have to be different from each other, i.e., $\Lambda$ has to be a permutation of $(n,\dots,1)$. In particular, some $\Lambda_i$ must be equal to $n$ and there are two possibilities:
\eq
\label{choice 1 C}
\Lambda_1=n\qquad\text{or}\qquad\Lambda_{p+1}=n.
\eeq
Assume first that $\Lambda_1=n$. Then
\eqn
a_1=0,
\eeqn
and it follows that
\eqn
\eps_1+\eps_{p+j}\notin A,\qquad 1\leq j\leq q.
\eeqn
This implies that
\eqn
0\leq b_j\leq p-1,\qquad 1\leq j\leq q,
\eeqn
and so
\begin{align*}
\Lambda=(n,\unb{[p-1,n-1]}{n-1-a_2},\dots,\unb{[1,q+1]}{q+1-a_p}\,|\,\unb{[q+1,n]}{n-b_1},\unb{[q,n-1]}{n-1-b_2},\dots,\\
\unb{[2,p+1]}{p+1-b_q}).
\end{align*}
If $p=1$, then there is only $\Lambda_1=n$ in the left group of coordinates, and we see there is no place to put the coordinate 1. Therefore, if $p=1$ then $\Lambda_1$ can not be $n$, hence $\Lamdba_{p+1}=n$, so we are in the second case which we treat below. If $p>1$, then there is exactly one place where 1 can be, i.e., 
\[
\Lambda_p=1.
\]
This implies 
\eqn
a_p=q,
\eeqn
and therefore
\eqn
\eps_p+\eps_{p+j}\in A,\qquad 1\leq j\leq q.
\eeqn
It follows that
\eqn
1\leq b_j\leq p-1,\qquad 1\leq j\leq q,
\eeqn
and so
\begin{align*}
\Lambda=(n,\unb{[p-1,n-1]}{n-1-a_2},\dots,\unb{[2,q+2]}{q+2-a_{p-1}},1\,|\,\unb{[q+1,n-1]}{n-b_1},\unb{[q,n-2]}{n-1-b_2},\dots,\\
\unb{[2,p]}{p+1-b_q}).
\end{align*}
Let now $\frg'\cong\frs\frp(2(n-2),\bbR)$ be the subalgebra of $\frg$ built on coordinates $2,\dots,p-1,p+1,\dots,n$, and let $\frl'=\frl\cap\frg'$. Then
\[
\Delta_n^+(\frl')=\Delta_n^+(\frl)\setminus\{\eps_1+\eps_{p+j},\eps_p+\eps_{p+j}\,\big|\,1\leq j\leq q\},
\]
and we set
\[
A'=A\setminus\{\eps_p+\eps_{p+j}\,\big|\,1\leq j\leq q\}.
\]
We define $\lambda_0$ as in \eqref{la0 C}, but with $n$ replaced by $n-2$ and $p$ replaced by $p-2$. Then $\Lambda'$ corresponding to $A'$ can be obtained from $\Lamdba$ by deleting coordinates $\Lambda_1$ and $\Lambda_p$, and decreasing all the other coordinates by $1$.
We now see that $\Lambda$ is a permutation of $(n,\dots,1)$ if and only if $\Lamdba'$ is a permutation of $(n-2,\dots,1)$. By inductive assumption, this is equivalent to $A'$ and $\Lambda'$ being defined by a shuffle as in the statement of the proposition, and this clearly implies the same statement for $A$ and $\Lamdba$.

The other possibility in \eqref{choice 1 C} is handled analogously: $\Lambda_{p+1}=n$ implies
\eqn
b_1=0,
\eeqn
and it follows that
\eqn
\eps_i+\eps_{p+1}\notin A,\qquad 1\leq i\leq p.
\eeqn
This implies that
\eqn
0\leq a_i\leq q-1,\qquad 1\leq i\leq p,
\eeqn
and so
\begin{align*}
\Lambda=(\unb{[p+1,n]}{n-a_1},\unb{[p,n-2]}{n-1-a_2},\dots,\unb{[2,q]}{q+1-a_p}\,|\,n,\unb{[q-1,n-1]}{n-1-b_2},\dots,\\
\unb{[1,p+1]}{p+1-b_q}).
\end{align*}
If $q=1$, then there is only $\Lambda_{p+1}=n$ in the right group of coordinates, and we see there is no place to put the coordinate 1. Therefore, if $q=1$, $\Lambda_{p+1}$ can not be $n$ and we are back to the first case that we already handled. If $q>1$, then there is exactly one place where 1 can be, i.e., 
\[
\Lambda_n=1.
\]
This implies 
\eqn
b_q=p,
\eeqn
and therefore
\eqn
\eps_i+\eps_n\in A,\qquad 1\leq i\leq p.
\eeqn
It follows that
\eqn
1\leq a_i\leq q-1,\qquad 1\leq i\leq p,
\eeqn
and so
\begin{align*}
\Lambda=(\unb{[p+1,n-1]}{n-a_1}\unb{[p,n-2]}{n-1-a_2},\dots,\unb{[2,q]}{q+1-a_p}\,|\,n,\unb{[q-1,n-1]}{n-1-b_2},\unb{[q-2,n-2]}{n-2-b_3},\dots,\\
\unb{[2,p+1]}{p+2-b_{q-1}},1).
\end{align*}
We now reason in the same way as in the first case, and conclude that the proposition follows from the inductive assumption for $n-2$ with $p$ staying the same and $q$ being replaced by $q-2$.
\epf

To finish the computation of the constant $c_p^n$, we first note that for every $A$ described in Proposition \ref{prop C}(ii)
\eq
\label{cardA C}
\# A=\frac{pq}{2}.
\eeq
Namely, the $\alpha_{u,v}$ and $\beta_{u,v}$ make for $2rs$ elements of $A$. In addition, the set $B$ has $0$ elements if $p$ and $q$ are even, $s$ elements if $p$ is odd, and $r$ elements if $q$ is odd. So the total number of elements is
\begin{align*}
2rs=\frac{pq}{2},\qquad p,q\text{ even};\\
2rs+s=(2r+1)s=p\frac{q}{2},\qquad p\text{ odd};\\
2rs+r=r(2s+1)=\frac{p}{2}q,\qquad q\text{ odd}.
\end{align*}
On the other hand, since $\Lambda$ is a permutation of $(n,\dots,1)$, $P_K(\Lambda)$ is equal to $\pm 1$. To compute the sign, we need to find the parity of the permutation bringing $\Lamdba$ to $(n,\dots,1)$. This parity can be found by counting the number of inversions in $\Lambda$ when compared with $(n,\dots,1)$, i.e., counting the number of pairs $(i,j)$, $1\leq i<j\leq n$, such that $\Lambda_i<\Lambda_j$.
We know from Proposition \ref{prop C} that 
\eq
\label{final La C}
\Lamdba=(n+1-i_1,\dots,n+1-i_r,i_r,\dots,i_1\,|\,n+1-j_1,\dots,n+1-j_s,j_s,\dots,j_1)
\eeq
if $p$ and $q$ are both even. It is clear that $i_r,\dots,i_1$ are in inversion with $n+1-j_1,\dots,n+1-j_s$; that is $rs$ inversions. The further inversions are possible only
between groups 
\eq
\label{inversions1 C}
n+1-i_1,\dots,n+1-i_r \qquad\text{and}\qquad n+1-j_1,\dots,n+1-j_s,
\eeq 
and
\eq
\label{inversions2 C}
i_r,\dots,i_1 \qquad\text{and}\qquad   j_s,\dots,j_1.
\eeq
If $i_u$ is in inversion with $j_v$, i.e., $i_u<j_v$, then
$n+1-i_u>n+1-j_v$, i.e., $n+1-i_u$ is not in inversion with $n+1-j_v$. The converse also holds, and it follows that the total number of inversions in groups \eqref{inversions1 C} and \eqref{inversions2 C} is again $rs$.
So the total number of inversions in case $p$ and $q$ are even is
\[
2rs=\frac{pq}{2}.
\]
If $p$ is odd, then $\Lamdba$ is again given by \eqref{final La C}, except that there is in addition $r+s+1$ between $n-i_r$ and $i_r$. This coordinate is in inversion with the coordinates $n+1-j_1,\dots,n+1-j_s$, and with no others, so the total number of inversions in this case is
\[
2rs+s=\frac{pq}{2}.
\]
Similarly, if $q$ is odd, then $\Lamdba$ is given by \eqref{final La C}, with the addition of $r+s+1$ between $n+1-j_s$ and $j_s$. This coordinate is in inversion with the coordinates $i_r,\dots,i_1$, and with no others, so the total number of inversions in this case is
\[
2rs+r=\frac{pq}{2}.
\]
So we have proved that for each $\Lambda$ from Proposition \ref{prop C},
\[
P_K(\Lambda)=(-1)^{\frac{pq}{2}}.
\]
Combined with \eqref{cardA C}, and with the fact that $N$ from \eqref{def N} is in this case
\[
N=\binom{p}{2}+pq+p\equiv [\frac{p+1}{2}]\quad\text{mod }2,
\]
this tells us that the nonzero contributions to the sum in \eqref{def const C} are all equal to $(-1)^{[\frac{p+1}{2}]}$. Since the number of nonzero summands is by Proposition \ref{prop C} equal to the number of $(r,s)$-shuffles of $r+s$, i.e., to $\binom{r+s}{r}$, we have proved:

\begin{theorem}
\label{thm C}
Let $G_\bbR=\Sp(2n,\bbR)$, $n\geq 1$, and let $p$, $0\leq p\leq n$, be an integer. Let $r=[\frac{p}{2}]$ and let $s=[\frac{n-p}{2}]$. Then the constant $c_p^n$ for the real form of $\caO^\bbC$ corresponding to $p$ is
\[
c_p^n=\left\{\begin{matrix} 
&0,&\qquad \text{if } n \text{ is even and } p \text{ is odd};\cr
&(-1)^{[\frac{p+1}{2}]}\binom{r+s}{r},&\qquad \text{if } n \text{ is odd, or if } n\text{ is even and } p \text{ is even}.
\end{matrix}\right.
\]
\end{theorem}

\bigskip

\section{The case $G_\bbR=SO_e(2p,2q)$, $q\geq p\geq 1$} 
\label{sec:D}

\bigskip

There are three real forms of $\caO^\bbC$ if $q>p>1$, four if $q=p>1$, and two if p=1. 

\subsection{The first real form}
\label{D1}
This real form is defined in all cases; it corresponds to
\eqn
h_1=(2,\unb{p-1}{1,\dots,1}\,|\,\unb{p-1}{1,\dots,1},\,\unb{q-p+1}{0,\dots,0}).
\eeqn

Since $\frl=\frl_1$ is built from roots that vanish on $h_1$, we see that
\[
\Delta^+_n(\frl)=\{\eps_i-\eps_j\,\big|\, 2\leq i\leq p,\, p+1\leq j\leq 2p-1\}.
\]
It follows that for any $A\subseteq\Delta^+_n(\frl)$,
\eq
\label{2rhoA D1}
2\rho(A)=(0;a_1,\dots,a_{p-1}\,|\,-b_1,\dots,-b_{p-1};0,\dots,0),
\eeq
with 
\eq
\label{ineq ab D1}
0\leq a_i,b_j\leq p-1;\qquad \textstyle{\sum_i a_i=\sum_j b_j}.
\eeq
Furthermore, recall that $\Delta(\frp_1)$ consists of noncompact roots that are 1 on $h_1$. So 
\[
\Delta(\frp_1)=\{\eps_1-\eps_j\,\big|\, p+1\leq j\leq 2p-1\}\,\cup \,\{\eps_i\pm\eps_j\,\big|\, 2\leq i\leq p,\, 2p\leq j\leq p+q\}.
\]
It follows that for any $C\subseteq\Delta(\frp_1)$,
\eq
\label{2rhoC D1}
2\rho(C)=(c;d_1,\dots,d_{p-1}\,|\,-c_1,\dots,-c_{p-1};e_1,\dots,e_{q-p+1}),
\eeq
with 
\begin{align}
\label{ineq cde D1}
&0\leq c_j\leq 1;\qquad 0\leq c\leq p-1;\qquad \textstyle{c=\sum_j c_j;} \\ 
&0\leq d_i\leq 2(q-p+1);\qquad -(p-1)\leq e_j\leq p-1. \nonumber
\end{align}
(If $p=1$, then
\[
h_1=(2\,|\,0,\dots,0),
\]
so $\Delta_n^+(\frl)=\Delta(\frp_1)=\emptyset$.)

By \eqref{2rhoA D1},
\[
\rho_n(\frl)=(0,p-1,\dots,p-1\,|\,-p+1,\dots,-p+1,0,\dots,0).
\]
This is clearly orthogonal to all roots of $\frl\cap\frk$, which are equal to
\begin{align}
\label{roots lk D1}
\Delta(\frl\cap\frk)=\{\eps_i-\eps_j\,\big|\,2\leq i,j\leq p\}\cup\{\eps_i-\eps_j\,\big|\,p+1\leq i,j\leq 2p-1\}\cup\\
\nonumber \cup\{\eps_i\pm\eps_j\,\big|\,2p\leq i,j\leq p+q\}.
\end{align}
By Proposition \ref{eq const}, this means that we can determine the constant $c=c_1^{p,q}$ from the equation \eqref{def const}.
We apply \eqref{def const} for $\lambda=\lambda_0$, where
\begin{align}
\label{la0 D1}
&\lambda_0=(\half;q-\half,q-\frac{3}{2}\dots,q-p+\frac{3}{2}\,|\, \\ 
& \qquad\qquad\qquad |\,-\frac{3}{2},-\frac{5}{2},\dots,-(p-\half);q-p+\half,\dots,\frac{3}{2},\half) 
\qquad  \text{if }\ p\geq 2;\nonumber \\
&\lambda_0=(\half\,|\, q-\half,q-\frac{3}{2},\dots,\half) \nonumber
\qquad \text{if }\ p=1.
\end{align}

\begin{proposition}
\label{prop D1} 
Let $\Lambda=\lambda_0-2\rho(A)-2\rho(C)$, with $\lambda_0$ given by \eqref{la0 D1}, and with $A\subseteq\Delta_n^+(\frl)$ and $C\subseteq\Delta(\frp_1)$. If $P_K(\Lambda)\neq 0$, then:
\begin{enumerate}
\item If $p=1$, then $A=C=\emptyset$ and $\Lambda=\lambda_0$.
\item If $p\geq 2$, then
\begin{align*}
&A=\emptyset; \\
&C=\{\eps_i-\eps_j\,\big|\, 2\leq i\leq p,\, 2p\leq j\leq p+q-1\}; \\
&\Lambda=(\half;p-\half,p-\frac{3}{2},\dots,\frac{3}{2}\,|\, -\frac{3}{2},-\frac{5}{2},\dots,-(p-\half);q-\half,q-\frac{3}{2},\dots,p+\half,\half).
\end{align*}
\end{enumerate}
\end{proposition}
\pf
The proposition is clear if $p=1$, since in that case 
$\Delta^+_n(\frl)=\Delta(\frp_1)=\emptyset$. We continue by induction on $p$. Let $p\geq 2$ and let $q\geq p$ be arbitrary. We assume that the proposition is true for $p'=p-1$ and $q'=p-1$, and we show it is then also true for $p$ and $q$.

By \eqref{la0 D1}, \eqref{2rhoA D1} and \eqref{2rhoC D1}, we have
\begin{align*}
%\label{La D2}
\textstyle{\Lambda=(\half-c;q-\half-a_1-d_1,q-\frac{3}{2}-a_2-d_2,\dots,q-p+\frac{3}{2}-a_{p-1}-d_{p-1}\,|}\\ 
\textstyle{|\,-\frac{3}{2}+b_1+c_1,-\frac{5}{2}+b_2+c_2,\dots,-(p-\half)+b_{p-1}+c_{p-1};} \nonumber  \\ 
\textstyle{q-p+\half-e_1,q-p-\half-e_2,\dots,\frac{3}{2}-e_{q-p},\half-e_{q-p+1}).} \nonumber 
\end{align*}

Using \eqref{ineq ab D1} and \eqref{ineq cde D1}, we see that the coordinates $\Lambda_{p+1},\dots,\Lambda_{p+q}$ are in the following intervals:
\begin{align*}
\Lambda=(\dots|\,\unb{[-\frac{3}{2},p-\frac{3}{2}]}{\textstyle{-\frac{3}{2}+b_1+c_1}},\unb{[-\frac{5}{2},p-\frac{5}{2}]}{\textstyle{-\frac{5}{2}+b_2+c_2}},\dots,\unb{[-(p-\half),\half]}{\textstyle{-(p-\half)+b_{p-1}+c_{p-1}}};  \\ \unb{[q-2p+\frac{3}{2},q-\half]}{\textstyle{q-p+\half-e_1}},\unb{[q-2p+\half,q-\frac{3}{2}]}{\textstyle{q-p-\half-e_2}},\dots,\unb{[-(p-\frac{5}{2}),p+\half]}{\textstyle{\frac{3}{2}-e_{q-p}}}\unb{[-(p-\frac{3}{2}),p-\half]}{\textstyle{\half-e_{q-p+1}}}).  
\end{align*}
So $\Lambda_{p+1},\dots,\Lambda_{p+q}$ are $q$ half-integers between $-(p-\half)$ and $q-\half$. Moreover, $P_K(\Lambda)\neq 0$ implies that no two of these half-integers are equal or opposite to each other. If $q>p$, it follows that $q-\half,q-\frac{3}{2},\dots,p+\half$ must each be equal to some $\Lambda_i$, and the only possibility for that is
\eqn
\Lambda_{2p}=q-\half,\ \Lambda_{2p+1}=q-\frac{3}{2},\ \dots,\ \Lambda_{p+q-1}=p+\half.
\eeqn
So $e_1,\dots,e_{q-p}$ are all equal to $-(p-1)$, and hence
\eqn
\eps_i-\eps_j\in C,\ \eps_i+\eps_j\notin C,\qquad 2\leq i\leq p,\ 2p\leq j\leq p+q-1.
\eeqn
(If $q=p$, the above says nothing and should be skipped.)

This implies
\eq
\label{ineq d D1}
q-p\leq d_i\leq q-p+2,\qquad 1\leq i\leq p-1,
\eeq
with $d_i$ being $q-p$ if $\eps_{i+1}\pm\eps_{p+q}\notin C$,
$d_i=q-p+2$ if $\eps_{i+1}\pm\eps_{p+q}\in C$, and $d_i=q-p+1$ if one of the roots $\eps_{i+1}\pm\eps_{p+q}$ is in $C$ while the other is not in $C$. 
(If $q=p$, this gives no new information about the $d_i$.
The following arguments all work also in case $q=p$ if we delete the group of coordinates from place $2p$ to place $p+q-1$.) 

Looking at the bounds for coordinates $\Lambda_{p+1},\dots,\Lambda_{2p-1}$ and $\Lambda_{p+q}$, we see that they are $p$ half-integers between $-(p-\half)$ and $p-\half$, such that no two of them are equal or opposite to each other. It follows that some of these $\Lambda_j$ must be equal to
$\pm(p-\half)$. There are two possibilities:
\[
\Lambda_{2p-1}=-(p-\half)\qquad\text{or}\qquad\Lambda_{p+q}=p-\half.
\]
Let us first examine the possibility that $\Lambda_{p+q}=p-\half$. If this is true, then $e_{p+q}=-(p-1)$, so
\[
\eps_{i+1}-\eps_{p+q}\in C,\quad\eps_{i+1}+\eps_{p+q}\notin C,\qquad 1\leq i\leq p-1,
\]
and it follows that
\[
d_i=q-p+1,\qquad 1\leq i\leq p-1.
\]
Using this together with the inequalities \eqref{ineq ab D1}, \eqref{ineq cde D1} for $a_i$ and $c$, we see
\[
\Lambda=(\unb{[-(p-\frac{3}{2}),\half]}{\textstyle{\half-c}};
\unb{[-\half,p-\frac{3}{2}]}{\textstyle{p-\frac{3}{2}-a_1}},
\unb{[-\frac{3}{2},p-\frac{5}{2}]}{\textstyle{p-\frac{5}{2}-a_2}},
\dots,
\unb{[-(p-\frac{3}{2}),\half]}{\textstyle{\half-a_{p-1}}}\,|\dots).
\]
So $\Lambda_1,\dots,\Lambda_p$ are $p$ half-integers between $-(p-\frac{3}{2})$ and $p-\frac{3}{2}$. Moreover, $P_K(\Lambda)\neq 0$ implies that no two of these half-integers are equal or opposite to each other. This is impossible, and so $\Lambda_{p+q}$ can not be $p-\half$.

It follows that 
\[
\Lambda_{2p-1}=-(p-\half),
\]
and hence
\eqn
b_{p-1}=0;\qquad c_{p-1}=0.
\eeqn
This implies that
\begin{align*}
&\eps_i-\eps_{2p-1}\notin A,\qquad 2\leq i\leq p; \\
& \eps_1-\eps_{2p-1}\notin C,
\end{align*}
and therefore 
\begin{align}
\label{ineq ac1 D1}
&0\leq a_i\leq p-2,\quad 1\leq i\leq p-1;\\
\nonumber & 0\leq c\leq p-2.
\end{align}
Using \eqref{ineq ac1 D1} and \eqref{ineq d D1}, we see that 
$\Lambda_1,\dots,\Lambda_p$ are in the following intervals:
\begin{align*}
\Lambda=(\unb{[-(p-\frac{5}{2}),\half]}{\textstyle{\half-c}};\unb{[-\half,p-\half]}{\textstyle{q-\half-a_1-d_1}},\unb{[-\frac{3}{2},p-\frac{3}{2}]}{\textstyle{q-\frac{3}{2}-a_2-d_2}},\dots,\unb{[-(p-\frac{3}{2}),\frac{3}{2}]}{\textstyle{q-p+\frac{3}{2}-a_{p-1}-d_{p-1}}}\,|\dots).
\end{align*}
So $\Lambda_1,\dots,\Lambda_p$ are $p$ half-integers between $-(p-\frac{3}{2})$ and $p-\half$. As before, these half-integers must be different from each other, and no two of them are opposite, so one of them must be equal to $p-\half$, and the only possibility is
\eqn
\Lambda_2=p-\half.
\eeqn
So 
\eqn
a_1=0;\qquad d_1=q-p.
\eeqn 
It follows that 
\begin{align*}
&\eps_2-\eps_j\notin A,\qquad p+1\leq j\leq 2p-1; \\ 
& \eps_2\pm\eps_{p+q}\notin C,
\end{align*}
and hence 
\begin{align*}
& 0\leq b_j\leq p-2,\qquad 1\leq j\leq p-1; \\
& -(p-2)\leq e_{q-p+1}\leq p-2.
\end{align*}

So we see that
\begin{align*}
\textstyle{\Lambda=(\half-c;p-\half,q-\frac{3}{2}-a_2-d_2,\dots,q-p+\frac{3}{2}-a_{p-1}-d_{p-1}\,|}\\ 
|\,\textstyle{-\frac{3}{2}+b_1+c_1,\dots,-(p-\frac{3}{2})+b_{p-2}+c_{p-2},-(p-\half);}  \\ 
\textstyle{q-\half,q-\frac{3}{2},\dots,p+\half,\half-e_{q-p+1}).} 
\end{align*}
(The coordinates $q-\half,q-\frac{3}{2},\dots,p+\half$ are not there if $q=p$.)

We now consider the subalgebra $\frg'\cong\frs\fro(2p-2,2p-3)$ of $\frg$ built on coordinates 
\[
\eps_1;\eps_3,\dots,\eps_p;\eps_{p+1},\dots,\eps_{2p-2};\eps_{p+q},
\]
so the coordinates 2 and $2p-1,2p,\dots,p+q-1$ are deleted. We also consider the real form of $\caO_{K'}$ given by 
\[
h'_1=(2,\unb{p-2}{1,\dots,1}\,|\,\unb{p-2}{1,\dots,1},0),
\]
with centralizer $\frl'=\frl\cap\frg'$.
Then
\begin{align*}
&\Delta^+_n(\frl')=\{\eps_i-\eps_j\,\big|\, 3\leq i\leq p,\,p+1\leq j\leq 2p-2\};\\
&\Delta(\frp'_1)=\{\eps_1-\eps_j\,\big|\,p+1\leq j\leq 2p-2\}\cup\{\eps_i\pm\eps_{p+q}\,\big|\,3\leq i\leq p\}.
\end{align*}
We set
\begin{align*}
& A'=A\cap\Delta^+_n(\frl')=A;\\
&C'=C\cap\Delta(\frp'_1)=C\setminus\{\eps_i-\eps_j\,\big|\,2\leq i\leq p,\ 2p\leq j\leq p+q-1\}.
\end{align*}
Then 
\begin{align*}
& 2\rho(A')=(0;a_2,\dots,a_{p-1}\,|\,-b_1,\dots,-b_{p-2};0)=   \\ 
& \qquad\qquad\qquad\qquad\qquad\qquad(0;a_1',\dots,a_{p-2}'\,|\,-b_1',\dots,-b_{p-2}';0) ;\\
&2\rho(C')=(c;d_2-(q-p),\dots,d_{p-1}-(q-p)\,|\,-c_1,\dots,-c_{p-2};e_{q-p+1})=\\
&\qquad\qquad\qquad\qquad\qquad\qquad(c';d_1',\dots,d_{p-1}'\,|\,-c_1',\dots,-c_{p-2}';e_{q-p+1}'),
\end{align*}
where we define
\begin{align*}
&a_i'=a_{i+1};\quad b_i'=b_i;\quad c_i'=c_i;\quad c'=c;\\ &d_i'=d_{i+1}-(q-p);\quad e_{q-p+1}'=e_{q-p+1}.
\end{align*}
The numbers $a_i',b_i',c_i',c',d_i'$ satisfy analogues of \eqref{ineq ab D1} and \eqref{ineq cde D1}.
We define $\lambda_0'$ by \eqref{la0 D1}, but for $G_\bbR=SO_e(2p-2,2p-1)$, i.e.,
\[
\lambda'_0=(\half;p-\frac{3}{2},\dots,\frac{3}{2}\,|\, -\frac{3}{2},-\frac{5}{2},\dots,-(p-\frac{3}{2});\half).
\]
Then $A'$, $C'$, and 
\[
\Lambda'=\lambda'_0-2\rho(A')-2\rho(C')
\]
satisfy all conditions of the proposition, but $p,q$ are reduced to $p'=p-1$, $q'=p-1$. Moreover, $P_K(\Lambda)\neq 0$ is equivalent to $P_{K'}(\Lambda')\neq 0$. Therefore the inductive assumption implies that $A'=C'=\emptyset$, and that
$\Lambda'=\lambda_0'$. This implies the statement of the proposition for $A$, $C$ and $\Lambda$.
\epf

In view of \eqref{def const}, to compute the constant $c=c_1^{p,q}$ we need to compute
$P_{L\cap K}(\lambda_0)$ and $P_K(\Lambda)$, where 
$\lambda_0$ is given by \eqref{la0 D1}, and $\Lambda$ is given by Proposition \ref{prop D1}.

To compute $P_{L\cap K}(\lambda_0)$, we note that we described $\frl\cap\frk$ in \eqref{roots lk D1}; it has up to three factors, two of which are $\fru(p-1)$, and the third is $\frs\fro(2(q-p+1))$. From the shape of $\lambda_0$ it now follows that, in the notation of Lemma \ref{PK},
\[
P_{L\cap K}(\lambda_0)=P^1_{q-p+1}(\lambda_{q-p+1}),
\]
and we see that Lemma \ref{PK}(i) implies that
\[
P_{L\cap K}(\lambda_0)=2^{q-p}.
\]
To compute $P_K(\Lambda)$, with $\Lambda$ as in Proposition \ref{prop D1}, we first write $\Lambda=(\Lambda_L\,|\,\Lambda_R)$ and note that 
\[
P_K(\Lambda)=P^1_p(\Lambda_L)P^1_q(\Lambda_R).
\] 
To use Lemma \ref{PK}, we have to rearrange coordinates of $\Lambda_L$ and $\Lambda_R$, using the fact that $P^1_p$ is skew for the Weyl group of $\frs\fro(2p)$, while $P_q^1$ is skew for the Weyl group of $\frs\fro(2q)$. Moreover, both polynomials are invariant under sign changes of the variables; this follows since the sign change of the $j$-th coordinate switches roots $\eps_i-\eps_j$ and $\eps_i+\eps_j$.

To rearrange $\Lambda_L$ to $\lambda_p$, we only need to bring the $\half$ from the first coordinate to the $p$-th coordinate, and hence
\[
P_p^1(\Lambda_L)=(-1)^{p-1}2^{p-1}.
\]
To bring $\Lambda_R$ to $\mu_q=(q,\dots,1)$, after removing the signs which does not change the expression, we need to bring coordinates $p-\half,p-\frac{3}{2},\dots,\frac{3}{2}$, in that order, to the right of $p+\half$, leaving $\half$ at the end. The sign produced in this way is
\[
(-1)^{(q-p)+(q-p+1)+\dots+(q-2)}=(-1)^{(p-1)(q-p-1)+\frac{(p-1)p}{2}}=(-1)^{(p-1)(q-p-1)+[\frac{p}{2}]},
\]
and it follows from Lemma \ref{PK} that 
\[
P_q^1(\Lambda_R)=(-1)^{(p-1)(q-p-1)+[\frac{p}{2}]}2^{q-1}.
\]
Putting this together with the fact that
\[
\# A+\# C=\#C=(p-1)(q-p),
\]
that $N$ of \eqref{def N} satisfies
\eq
\label{N D1}
N\equiv p-1\quad\text{mod }2,
\eeq
and that 
\[
[\frac{p}{2}]+p-1\equiv [\frac{p-1}{2}]\quad\text{mod }2,
\]
we see that \eqref{def const} implies
\begin{theorem}
\label{thm D1}
For $G_\bbR=SO_e(2p,2q)$, $q\geq p\geq 1$, the constant $c_1^{p,q}$ is
\[
c_1^{p,q}=(-1)^{[\frac{p-1}{2}]}2^{2p-2}.
\]
\qed
\end{theorem}

\subsection{The second real form}
\label{D2} 
This real form exists for $q\geq p\geq 2$. It corresponds to
\[
h_2=(2,\unb{p-2}{1,\dots,1},-1\,|\,\unb{p-1}{1,\dots,1},\unb{q-p+1}{0,\dots,0}).
\]
This real form is conjugate to the first real form by the automorphism $\sigma$ which acts on $\frh$ by changing the sign of the $p$-th coordinate, and leaving all other coordinates the same. On the level of $\frg$, this is an outer automorphism, which become the standard one if we compose it with the isomorphism $\frs\fro(2p,2q)\cong\frs\fro(2q,2p)$. The automorphism $\sigma$ satisfies the conditions of
Lemma \ref{auto}, and we just have to compute the sign. 
The number $n$ from Proposition \ref{auto} is, as in Subsection \ref{B1}, equal to $p-1$. The number $N_1$ is by \eqref{N D1} congruent to $p-1$ modulo 2. Finally, 
$N_2$ is easily seen to be always even. The conclusion is that there is no sign in Proposition \ref{auto}, so
\eq
\label{const D2}
c_2^{p,q}=c_1^{p,q}=(-1)^{[\frac{p-1}{2}]}2^{2p-2}.
\eeq

\subsection{The third real form}
\label{D3} 
This real form exists for all $q\geq p\geq 1$. It corresponds to
\[
h_3=(\unb{p-1}{1,\dots,1},0\,|\,2,\unb{p-1}{1,\dots,1},\unb{q-p}{0,\dots,0}).
\]

Since $\frl=\frl_3$ is built from roots that vanish on $h_3$, we see that
\[
\Delta^+_n(\frl)=\{\eps_i-\eps_j\,\big|\, 1\leq i\leq p-1,\, p+2\leq j\leq 2p\}\cup\{\eps_p\pm\eps_j\,\big|\,2p+1\leq j\leq p+q\}.
\]
It follows that for any $A\subseteq\Delta^+_n(\frl)$,
\eq
\label{2rhoA D3}
2\rho(A)=(a_1,\dots,a_{p-1};x\,|\,0;-b_1,\dots,-b_{p-1};y_1,\dots,y_{q-p}),
\eeq
with 
\begin{align}
\label{ineq abxy D3}
&0\leq a_i,b_j\leq p-1;\qquad \textstyle{\sum_i a_i=\sum_j b_j}; \\
&\nonumber 0\leq x\leq 2(q-p);\qquad -1\leq y_j\leq 1.
\end{align}
Furthermore, recall that $\Delta(\frp_1)$ consists of noncompact roots that are 1 on $h_3$. So 
\begin{multline*}
\Delta(\frp_1)= 
\{\eps_i\pm\eps_j\,\big|\, 1\leq i\leq p-1,\, 2p+1\leq j\leq p+q\}\cup  \\
\cup\, \{\eps_j\pm\eps_p\,\big|\, p+2\leq j\leq 2p\}\,\cup\,\{\eps_{p+1}-\eps_i\,\big|\,1\leq i\leq p-1\}.
\end{multline*}
It follows that for any $C\subseteq\Delta(\frp_1)$,
\eq
\label{2rhoC D3}
2\rho(C)=(c_1,\dots,c_{p-1};u\,|\,v;d_1,\dots,d_{p-1};e_1,\dots,e_{q-p}),
\eeq
with 
\begin{align}
\label{ineq cdeuv D3}
&-1\leq c_i\leq 2(q-p);\qquad -(p-1)\leq u\leq p-1; \\ 
&0\leq v\leq p-1;\qquad 0\leq d_j\leq 2;\qquad -(p-1)\leq e_j\leq p-1. \nonumber
\end{align}
If we write \eqref{2rhoA D3} for $A=\Delta^+_n(\frl)$, we get
\[
\rho_n(\frl)=(p-1,\dots,p-1;2(q-p)\,|\, 0;-p+1,\dots,-p+1;0,\dots,0).
\]
This is clearly orthogonal to all roots of $\frl\cap\frk$, which are equal to
\begin{align}
\label{roots lk D3}
\Delta(\frl\cap\frk)=\{\eps_i-\eps_j\,\big|\,1\leq i,j\leq p-1\}\cup\{\eps_i-\eps_j\,\big|\,p+2\leq i,j\leq 2p\}\cup\\
\nonumber \cup\{\eps_i\pm\eps_j\,\big|\,2p+1\leq i,j\leq p+q\}.
\end{align}
By Proposition \ref{eq const}, this means that we can determine the constant $c=c_3^{p,q}$ from the equation \eqref{def const}.
We will use this for $\lambda=\lambda_0$, where
\begin{align}
\label{la0 D3}
\lambda_0=(q-\frac{3}{2},q-\frac{5}{2}\dots,q-p+\frac{1}{2};q-p+\half\,|\,p-\frac{3}{2};\half,-\half,\dots,-p+\frac{5}{2};\\
\nonumber q-p-\half,q-p-\frac{3}{2},\dots,\half).
\end{align}
(Out of the 5 groups of coordinates separated by semicolons and the bar, the first and the fourth group are missing if  $p=1$, and the fifth group is missing if $q=p$.) 

\begin{proposition}
\label{prop D3} 
Let $\Lambda=\lambda_0-2\rho(A)-2\rho(C)$, with $\lambda_0$ given by \eqref{la0 D3}, and with $A\subseteq\Delta_n^+(\frl)$ and $C\subseteq\Delta(\frp_1)$. If $P_K(\Lambda)\neq 0$, then
\begin{enumerate}
\item If $q>p\geq 2$, then
\begin{align*}
&A=\{\eps_p-\eps_j\,\big|\,2p+1\leq j\leq p+q\}; \\
&C=\{\eps_i-\eps_j\,\big|\, 1\leq i\leq p-1,\, 2p+1\leq j\leq p+q\}\cup  \\
&\qquad\qquad\cup\, \{\eps_j\pm\eps_p\,\big|\, p+2\leq j\leq 2p\}\,\cup\,\{\eps_{p+1}-\eps_i\,\big|\,1\leq i\leq p-1\};\\
&\Lambda=(p-\half,p-\frac{3}{2},\dots,\frac{3}{2};\half\,|\, -\half;-\frac{3}{2},-\frac{5}{2},\dots,-(p-\half);\\
&\qquad\qquad\qquad\qquad\qquad\qquad\qquad\qquad\qquad q-\half,q-\frac{3}{2},\dots,p+\half).
\end{align*}
\item If $q>p=1$, then $A$ is as in (1), $C=\emptyset$, and 
\[
\Lambda=(\half\,|\,-\half;q-\half,q-\frac{3}{2},\dots,\frac{3}{2}).
\]
\item If $q=p\geq 2$, then 
\begin{align*}
&A=\emptyset; \\
&C=\{\eps_j\pm\eps_p\,\big|\, p+2\leq j\leq 2p\}\,\cup\,\{\eps_{p+1}-\eps_i\,\big|\, 1\leq i\leq p-1\}; \\
&\Lambda=(p-\half,p-\frac{3}{2},\dots,\frac{3}{2};\half\,|\, -\half;-\frac{3}{2},-\frac{5}{2},\dots,-(p-\half)).
\end{align*}
\item If $q=p=1$, then $A=C=\emptyset$ and $\Lambda=\lambda_0=(\half\,|\,-\half)$.
\end{enumerate}
\end{proposition}
\pf
By \eqref{la0 D3}, \eqref{2rhoA D3} and \eqref{2rhoC D3}, we have
\begin{align*}
\textstyle{\Lambda=(q-\frac{3}{2}-a_1-c_1,\dots,q-p+\half-a_{p-1}-c_{p-1};q-p+\half-x-u\,|}\\ 
|\,\textstyle{p-\frac{3}{2}-v;\half+b_1-d_1,-\half+b_2-d_2,\dots,-p+\frac{5}{2}+b_{p-1}-d_{p-1};}   \\ 
\textstyle{q-p-\half-y_1-e_1,q-p-\frac{3}{2}-y_2-e_2,\dots,\half-y_{q-p}-e_{q-p}).} 
\end{align*}
There are five groups of coordinates separated by semicolons and the bar. If $p=1$, then the first and the fourth group of coordinates are missing, and if $q=p$, then the fifth group of coordinates is missing.

Using \eqref{ineq abxy D3} and \eqref{ineq cdeuv D3}, we see 
\begin{align}
\label{segR D3}
\Lambda=(\dots|\,\unb{[-\half,p-\frac{3}{2}]}{\textstyle{p-\frac{3}{2}-v}};\unb{[-\frac{3}{2},p-\half]}{\textstyle{\half+b_1-d_1}},\unb{[-\frac{5}{2},p-\frac{3}{2}]}{\textstyle{-\half+b_2-d_2}},\dots,\unb{[-(p-\half),\frac{3}{2}]}{\textstyle{-p+\frac{5}{2}+b_{p-1}-d_{p-1}}};  \\ \unb{[q-2p-\half,q-\half]}{\textstyle{q-p-\half-y_1-e_1}},\unb{[q-2p-\frac{3}{2},q-\frac{3}{2}]}{\textstyle{q-p-\frac{3}{2}-y_2-e_2}},\dots,\unb{[-(p-\half),p+\half]}{\textstyle{\half-y_{q-p}-e_{q-p}}}).  \nonumber
\end{align}
So $\Lambda_{p+1},\dots,\Lambda_{p+q}$ are $q$ half-integers between $-(p-\half)$ and $q-\half$. Since $P_K(\Lambda)\neq 0$, no two of them are equal or opposite to each other. If $q>p$, it follows that $q-\half,q-\frac{3}{2},\dots,p+\half$ must each be equal to some $\Lambda_i$, and the only possibility for that is
\eq
\label{Lalast D3}
\Lambda_{2p+1}=q-\half,\ \Lambda_{2p+2}=q-\frac{3}{2},\ \dots,\ \Lambda_{p+q}=p+\half.
\eeq
It follows that $y_1,\dots,y_{q-p}$ are all equal to $-1$, and that $e_1,\dots,e_{q-p}$ are all equal to $-(p-1)$. 

In case $q>p=1$, this implies that $A$ is as stated in the proposition, and it is also clear that $C=\emptyset$. Moreover, it follows that $x=q-p$, and so $\Lambda$ is as stated in the proposition. Since the case $q=p=1$ is obvious, this proves the proposition for $p=1$ and any $q\geq p$.

We proceed by induction on $p$. Let $p\geq 2$ and let $q\geq p$ be arbitrary. Assuming that the proposition is true for $p'=q'=p-1$, we will prove it for $p$ and $q$. 

If $q>p$, we get back to \eqref{Lalast D3} and see that it implies
\begin{align*}
& \eps_p-\eps_j\in A,\ \eps_p+\eps_j\notin A,\qquad 2p+1\leq j\leq p+q;\\
& \eps_i-\eps_j\in C,\ \eps_i+\eps_j\notin C,\qquad 1\leq i\leq p-1,\ 2p+1\leq j\leq p+q. 
\end{align*} 
This implies
\begin{align}
\label{ineq xc D3}
& x= q-p; \\
&\nonumber q-p-1\leq c_i\leq q-p,\qquad 1\leq i\leq p-1.
\end{align}
Note that $c_i=q-p-1$ if $\eps_{p+1}-\eps_i\in C$, and $c_i=q-p$ if $\eps_{p+1}-\eps_i\notin C$.

If $q=p$, then the above discussion does not apply; in this case \eqref{ineq xc D3} is true, but this information is already contained in \eqref{ineq abxy D3} and \eqref{ineq cdeuv D3}. The following discussion applies equally to $q>p$ and $q=p$, but in the latter case the last group of coordinates should be deleted.

Using \eqref{ineq xc D3} together with the inequalities \eqref{ineq abxy D3} for $a_i$ and \eqref{ineq cdeuv D3} for $u$, we see 
\begin{align*}
\Lambda=(\unb{[-\half,p-\half]}{\textstyle{q-\frac{3}{2}-a_1-c_1}},\unb{[-\frac{3}{2},p-\frac{3}{2}]}{\textstyle{q-\frac{5}{2}-a_2-c_2}},\dots,\unb{[-(p-\frac{3}{2}),\frac{3}{2}]}{\textstyle{q-p+\half-a_{p-1}-c_{p-1}}}; \\
\unb{[-(p-\frac{3}{2}),p-\half]}{\textstyle{\half-u}}\,|\, \dots 
\end{align*}
So $\Lambda_1,\dots,\Lambda_p$ are $p$ half-integers between $-(p-\frac{3}{2})$ and $p-\half$, such that no two of them are equal or opposite to each other. There are two possibilities:
\eqn
\Lambda_1=p-\half,\qquad\text{or}\qquad \Lambda_p=p-\half
\eeqn
Let us first assume that $\Lambda_p=p-\half$. Then 
\[
u=-(p-1)
\] 
and it follows that 
\[
\eps_j-\eps_p\in C,\,\eps_j+\eps_p\notin C,\qquad p+2\leq j\leq 2p.
\]
This implies
\[
d_j=1,\qquad p+2\leq j\leq 2p,
\]
and we see that \eqref{segR D3} becomes
\begin{align*}
%\label{segR false D3}
\Lambda=(\dots|\,
\unb{[-\half,p-\frac{3}{2}]}{\textstyle{p-\frac{3}{2}-v}};\unb{[-\half,p-\frac{3}{2}]}{\textstyle{-\half+b_1}},
\unb{[-\frac{3}{2},p-\frac{5}{2}]}{\textstyle{-\frac{3}{2}+b_2}},\dots,
\unb{[-(p-\frac{3}{2}),\half]}{\textstyle{-p+\frac{3}{2}+b_{p-1}}}; \textstyle{q-\half,q-\frac{3}{2},\dots,
p+\half).} 
\end{align*}
So $\Lambda_{p+1},\dots,\Lambda_{2p}$ are $p$ half-integers between $-(p-\frac{3}{2})$ and $p-\frac{3}{2}$, such that no two of them are equal or opposite to each other, and this is impossible.

We conclude that
\[
\Lambda_1=p-\half,
\]
and consequently 
\eqn
a_1=0;\qquad c_1=q-p-1.
\eeqn
It follows that
\begin{align*}
&\eps_1-\eps_j\notin A,\quad p+2\leq j\leq 2p;\\ 
&\eps_{p+1}-\eps_1\in C, \nonumber
\end{align*}
and therefore
\begin{align*}
&0\leq b_j\leq p-2,\quad 1\leq j\leq p-1; \\
&1\leq v\leq p-1. \nonumber
\end{align*}
This implies
\begin{align*}
\Lambda=(\dots|\,\unb{[-\half,p-\frac{5}{2}]}{\textstyle{p-\frac{3}{2}-v}};\unb{[-\frac{3}{2},p-\frac{3}{2}]}{\textstyle{\half+b_1-d_1}},\unb{[-\frac{5}{2},p-\frac{5}{2}]}{\textstyle{-\half+b_2-d_2}},\dots,\unb{[-(p-\half),\frac{1}{2}]}{\textstyle{-p+\frac{5}{2}+b_{p-1}-d_{p-1}}}; \\ 
\textstyle{q-\half,q-\frac{3}{2},\dots,p+\half).}  
\end{align*}
So $\Lambda_{p+1},\dots,\Lambda_{2p}$ are $p$ half-integers between $-(p-\half)$ and $p-\frac{3}{2}$, such that no two of them are equal or opposite to each other. It follows that 
\[
\Lambda_{2p}=-(p-\half),
\]
and consequently
\eqn
b_{p-1}=0;\qquad d_{p-1}=2.
\eeqn
It follows that
\begin{align*}
&\eps_i-\eps_{2p}\notin A,\quad 1\leq i\leq p-1;\\ 
&\eps_{2p}\pm\eps_p\in C, 
\end{align*}
and therefore
\begin{align*}
&0\leq a_i\leq p-2,\quad 1\leq i\leq p-1; \\
&-(p-2)\leq u\leq p-2. 
\end{align*}
We see that
\begin{align*}
\textstyle{\Lambda=(p-\half,q-\frac{5}{2}-a_2-c_2,\dots,q-p+\frac{1}{2}-a_{p-1}-c_{p-1};q-p+\half-x-u\,|}\\ 
|\,\textstyle{p-\frac{3}{2}-v;\half+b_1-d_1,-\half+b_2-d_2,\dots,-(p-\frac{7}{2})+b_{p-2}-d_{p-2},-(p-\half);}   \\
\textstyle{q-\half,q-\frac{3}{2},\dots,p+\half). }
\end{align*}
(If $q=p$, the coordinates $q-\half,\dots,p+\half$ are not there; if $p=2$ there are no coordinates involving $a_i,c_i,b_i$ or $d_i$.)

We now consider the subalgebra $\frg'\cong\frs\fro(2p-2,2p-1)$ of $\frg$ built on coordinates 
\[
\eps_2,\eps_3,\dots,\eps_p;\eps_{p+1},\dots,\eps_{2p-1},
\]
so the coordinates 1 and $2p,2p+1,\dots,p+q$ are deleted. We also consider the real form of $\caO_{K'}$ given by 
\[
h'_3=(\unb{p-2}{1,\dots,1},0\,|\,2,\unb{p-2}{1,\dots,1}),
\]
with centralizer $\frl'=\frl\cap\frg'$. Then
\begin{align*}
&\Delta^+_n(\frl')=\{\eps_i-\eps_j\,\big|\,2\leq i\leq p-1,\,p+2\leq j\leq 2p-1\};\\
&\Delta(\frp'_1)=\{\eps_j\pm\eps_p\,\big|\,p+2\leq j\leq 2p-1\}\cup\{\eps_{p+1}-\eps_i\,\big|\,2\leq i\leq p-1\}.
\end{align*}
We set
\begin{align*}
& A'=A\cap\Delta^+_n(\frl')=A\setminus\{\eps_p-\eps_j\,\big|\,2p+1\leq j\leq p+q\};\\
&C'=C\cap\Delta(\frp'_1)=\\
&\qquad\qquad C\setminus\{\eps_{2p}\pm\eps_p;\,\eps_{p+1}-\eps_1;\,\eps_i-\eps_j\,\big|\,1\leq i\leq p-1,\ 2p+1\leq j\leq p+q\}.
\end{align*}
(If $q=p$, then $A'=A$ and $C'=C\setminus\{\eps_{2p}\pm\eps_p;\,\eps_{p+1}-\eps_1\}$.)
Then 
\begin{align*}
& 2\rho(A')=(a_2,\dots,a_{p-1};0\,|\,0;-b_1,\dots,-b_{p-2})=  \\  
&\qquad\qquad\qquad\qquad(a_1',\dots,a_{p-2}';0\,|\,0;-b_1',\dots,-b_{p-2}') ;\\
&2\rho(C')=(c_2-(q-p),\dots,c_{p-1}-(q-p);u\,|\,v-1;d_1,\dots,d_{p-2})=\\
&\qquad\qquad\qquad\qquad\qquad(c_1',\dots,c_{p-2}';u'\,|\,v';d_1',\dots,d_{p-2}'),
\end{align*}

where we define
\begin{align*}
&a_i'=a_{i+1};\quad b_i'=b_i;\quad c_i'=c_{i+1}-(q-p); \\ 
&u'=u;\quad v'=v-1;\quad d_i'=d_{i}.
\end{align*}
The numbers $a_i',b_i',c_i',u',v',d_i'$ satisfy analogues of \eqref{ineq abxy D3} and \eqref{ineq cdeuv D3}.

We define $\lambda_0'$ by \eqref{la0 D3}, but for $G_\bbR=SO_e(2p-2,2p-1)$, i.e.,
\[
\lambda'_0=(p-\frac{5}{2},p-\frac{7}{2},\dots,\half;\half;\,|\,p-\frac{5}{2};\half,-\half,\dots,-(p-\frac{7}{2})).
\]
Then $A'$, $C'$ and 
\[
\Lambda'=\lambda'_0-2\rho(A')-2\rho(C')
\]
satisfy all conditions of the proposition, but $p,q$ are reduced to $p'=p-1$, $q'=p-1$. Therefore the inductive assumption implies that $A'=\emptyset$, that
\[
C'=\{\eps_j\pm\eps_p\,\big|\, p+2\leq j\leq 2p-1\}\,\cup\,\{\eps_{p+1}-\eps_i\,\big|\,2\leq i\leq p-1\},
\]
and that
\[
\Lambda'=(p-\frac{3}{2},\dots,\frac{3}{2};\half\,|\,-\half;-\frac{3}{2},-\frac{5}{2},\dots,-(p-\frac{3}{2})).    
\]
This implies the statement of the proposition for $A$, $C$ and $\Lambda$.
\epf

In view of \eqref{def const}, to compute the constant $c=c_3^{p,q}$ we need to compute
$P_{L\cap K}(\lambda_0)$ and $P_K(\Lambda)$, where 
$\lambda_0$ is given by \eqref{la0 D3}, and $\Lambda$ is given by Proposition \ref{prop D3}.

To compute $P_{L\cap K}(\lambda_0)$, we note that we described $\frl\cap\frk$ in \eqref{roots lk D3}; it has up to three factors, two of which are $\fru(p-1)$, and the third is $\frs\fro(2(q-p))$. From the shape of $\lambda_0$ it now follows that, in the notation of Lemma \ref{PK},
\[
P_{L\cap K}(\lambda_0)=P^1_{q-p}(\lambda_{q-p}),
\]
and in case $q>p$, we see that Lemma \ref{PK}(i) implies that
\[
P_{L\cap K}(\lambda_0)=2^{q-p-1}.
\]
If $q=p$, then $P_{L\cap K}(\lambda_0)=1$, which is not covered by the above formula. (In Lemma \ref{PK}, we could have defined $P^1_0=1$ and $\lambda_0=0$, but the formula in Lemma \ref{PK}(i) would not work for $p=0$.)

To compute $P_K(\Lambda)$ for $\Lambda$ as in Proposition \ref{prop D3}, we first write $\Lambda=(\Lambda_L\,|\,\Lambda_R)$ and note that 
\eq
\label{D3 PK factor}
P_K(\Lambda)=P^1_p(\Lambda_L)P^1_q(\Lambda_R).
\eeq 
By Lemma \ref{PK}(i), 
\eq
\label{D3 PK1}
P^1_p(\Lambda_L)=P^1_p(\lambda_p)=2^{p-1}.
\eeq
To apply Lemma \ref{PK} also for $\Lambda_R$, we must first rearrange coordinates of $\Lambda_R$, using the fact that $P_q^1$ is skew for the Weyl group of $\frs\fro(2q)$, and invariant under sign changes of the variables.

To bring $\Lambda_R$ to $\mu_q=(q,\dots,1)$, after removing the signs which does not change the expression, we need to bring coordinates 
\[
p-\half,p-\frac{3}{2},\dots,\frac{3}{2},\half,
\] 
in that order, to the right of $p+\half$. The sign produced in this way is
\[
(-1)^{(q-p)+(q-p+1)+\dots+(q-1)}=(-1)^{p(q-p)+[\frac{p}{2}]},
\]
and it follows from Lemma \ref{PK} that 
\[
P_q^1(\Lambda_R)=(-1)^{p(q-p)+[\frac{p}{2}]}2^{q-1}.
\]
Putting this together with \eqref{D3 PK factor}, \eqref{D3 PK1}, the fact that
\[
\# A=q-p;\qquad \# C=3(p-1)+(p-1)(q-p),
\]
and the fact that $N$ of \eqref{def N} satisfies
\eq
\label{N D3}
N \equiv p\quad\text{mod }2
\eeq
we get
\begin{theorem}
\label{thm D3}
For $G_\bbR=SO_e(2p,2q)$, the constant $c_3^{p,q}$ corresponding to the third real form of $\caO^\bbC$ is
\[
c_3^{p,q}=\left\{
\begin{matrix}
(-1)^{[\frac{p}{2}]+1}2^{2p-1}, & \quad q>p;\cr
(-1)^{[\frac{p}{2}]+1}2^{2p-2}, & \quad q=p.
\end{matrix}\right.
\]
\qed
\end{theorem}

\subsection{The fourth real form}
\label{D4} 
This real form exists if $q= p\geq 2$. It corresponds to
\[
h_4=(\unb{p-1}{1,\dots,1},0\,|\,2,\unb{p-2}{1,\dots,1},-1).
\]
This real form is conjugate to the third real form by the automorphism $\sigma$ which acts on $\frh$ by changing the sign of the last coordinate, and leaving all other coordinates the same. On the level of $\frg$, this is the standard outer automorphism. It satisfies the conditions of
Proposition \ref{auto}, and we just have to compute the sign. 
The number $n$ from Proposition \ref{auto} is equal to 0, since $\sigma$ preserves $\Delta^+$. The number $N_1$ is by \eqref{N D3} congruent to $p$ modulo 2. Finally, 
$N_2$ is the same as $N_1$ because $\sigma$ preserves $\Delta^+$ (or one can do a computation). 
So there is no sign in Proposition \ref{auto}, and
\eq
\label{const D4}
c_4^{p,q}=(-1)^{[\frac{p}{2}]+1}2^{2p-2}
\eeq

\bigskip

\section{The case $G_\bbR=SO^*(2n)$, $n\geq 1$} 
\label{sec so*}

\bigskip

\subsection{The case of even $n$}
\label{D5}
If $n$ is even, the real forms of $\caO^\bbC$ correspond to even integers $p$ such that $0\leq p\leq n$. We denote $n-p$ by $q$. The $h$ corresponding to $p$ is 
\eqn
h_p=(\unb{p}{1,\dots,1}\,|\,\unb{q}{-1,\dots,-1}).
\eeqn
Since $\frl=\frl_p$ is built from roots that vanish on $h_p$, we see that
\[
\Delta^+_n(\frl)=\{\eps_i+\eps_{p+j}\,\big|\, 1\leq i\leq p,\, 1\leq j\leq q\}.
\]
It follows that for any $A\subseteq\Delta^+_n(\frl)$,
\eq
\label{2rhoA D5}
2\rho(A)=(a_1,\dots,a_p\,|\,b_1,\dots,b_q),
\eeq
with 
\begin{align}
\label{ineq ab D5}
0\leq a_i\leq q, \qquad 0\leq b_j\leq p,\qquad \textstyle{\sum_ia_i=\sum_jb_j}.
\end{align}
In particular,
\[
\rho_n(\frl)=(q,\dots,q\,|\,p,\dots,p),
\]
and this is clearly orthogonal to the roots of $\frl\cap\frk$, which are given by
\eq
\label{roots lk D5}
\Delta^+(\frl\cap\frk)=\{\eps_i-\eps_j\,\big|\,1\leq i<j\leq p\}\cup\{\eps_{p+i}-\eps_{p+j}\,\big|\,1\leq i<j\leq q\}.
\eeq
So the constants $c=c^n_p$ can be calculated from \eqref{def const}. Since it is clear that in the present case 
\[
\Delta(\frp_1)=\emptyset,
\]
\eqref{def const} becomes
\eq
\label{def const D5}
\sum_{A\subseteq\Delta_n^+(\frl)} (-1)^{\#A}P_K(\lambda-2\rho(A))=c P_{L\cap K}(\lambda).
\eeq
We take $\lambda=\lambda_0$, where
\eq
\label{la0 D5}
\lambda_0=(n,n-1,\dots,q+1\,|\,n,n-1,\dots,p+1),
\eeq
(If $p$ is 0 or $n$, then there is only one group of coordinates in the above expression, and $\lambda_0=(n,n-1,\dots,1)$.)

Since $\lambda_0$ differs from $\rho_{\frl\cap\frk}$ by a weight orthogonal to all roots of $\frl\cap\frk$, 
\[
P_{L\cap K}(\lambda_0)=1.
\]
So to compute $c_p^n$ we have to compute the left side of \eqref{def const D5}. The following proposition describes the relevant $A$ and the corresponding $\Lambda$.

\begin{proposition}
\label{prop D5} 
Let $\Lambda=\lambda_0-2\rho(A)$, with $\lambda_0$ given by \eqref{la0 D5}, and with $A\subseteq\Delta_n^+(\frl)$. 

Suppose that for some $A$ the corresponding $\Lambda$ satisfies $P_K(\Lambda)\neq 0$. Then:
\begin{enumerate}
\item If $p=0$ or $q=0$, then $A=\emptyset$ and $\Lambda=\lambda_0=(n,n-1,\dots,1)$.
\item If $p,q>0$, let $r=\frac{p}{2}$ and $s=\frac{q}{2}$. Then there is a shuffle
\[
1\leq i_1<\dots<i_r\leq r+s;\qquad 1\leq j_1<\dots<j_s\leq r+s
\]
of $1,2,\dots,r+s$ such that 
\[
A=\{\alpha_{u,v},\beta_{u,v}\,\big|\, 1\leq u\leq r,\,1\leq v\leq s\},
\]
where 
\[
\alpha_{u,v}=\eps_{p+1-u}+\eps_{n+1-v};\qquad \beta_{u,v}=\left\{
\begin{matrix} 
\eps_{p+1-u}+\eps_{p+v},&\qquad i_u<j_v;\cr
\eps_u+\eps_{n+1-v},&\qquad i_u>j_v.
\end{matrix}
\right.
\]
The corresponding $\Lambda$ is
\begin{align*}
\Lambda=(n+1-i_1,\dots,n+1-i_r,i_r,\dots,i_1\,|\,
n+1-j_1,\dots,n+1-j_s,j_s,\dots,j_1).
\end{align*}
\end{enumerate}
\end{proposition}
\pf
The situation is combinatorially exactly the same as for $G_\bbR=\Sp(2n,\bbR)$, with $n$, $p$ and $q$ even. Therefore the proof of Proposition \ref{prop C} applies verbatim; the only difference is that the present proof is simpler because $p$ and $q$ are even.
\epf

The complete parallel with the case of $G_\bbR=\Sp(2n,\bbR)$, with $n$, $p$ and $q$ even extends also to the computation of $P_K(\Lamdba)$ for any $\Lambda$ from Proposition \ref{prop D5}, and the constant $c_p^n$. The only difference is that in the present case, $N$ of \eqref{def N} is 
\[
N=\binom{p}{2}+pq\equiv \frac{p}{2}\quad\text{mod }2,
\]
so the sign is now $(-1)^{\frac{p}{2}}$. We conclude

\begin{theorem}
\label{thm D5}
Let $G_\bbR=SO^*(2n)$, with $n\geq 2$ even. Let $p$, $0\leq p\leq n$, be an even integer. Let $r=\frac{p}{2}$ and let $s=\frac{n-p}{2}$. Then the constant $c_p^n$ for the real form of $\caO^\bbC$ corresponding to $p$ is
\[
c_p^n=(-1)^{\frac{p}{2}}\binom{r+s}{r}.
\]
\end{theorem}

\subsection{The case of odd $n$}
\label{D6}
For odd $n$ the real forms of $\caO^\bbC$ correspond to even integers $p$ such that $0\leq p\leq n-1$. We denote $n-1-p$ by $q$, so $q$ is another even integer. The $h$ corresponding to $p$ is 
\eqn
h_p=(\unb{p}{1,\dots,1}\,|\,0;\unb{q}{-1,\dots,-1}).
\eeqn
Since $\frl=\frl_p$ is built from roots that vanish on $h_p$, we see that
\[
\Delta^+_n(\frl)=\{\eps_i+\eps_{p+1+j}\,\big|\, 1\leq i\leq p,\, 1\leq j\leq q\}.
\]
It follows that for any $A\subseteq\Delta^+_n(\frl)$,
\eq
\label{2rhoA D6}
2\rho(A)=(a_1,\dots,a_p\,|\,0;b_1,\dots,b_q),
\eeq
with 
\begin{align}
\label{ineq ab D6}
0\leq a_i\leq q, \qquad 0\leq b_j\leq p,\qquad \textstyle{\sum_ia_i=\sum_jb_j}.
\end{align}
In particular,
\[
\rho_n(\frl)=(q,\dots,q\,|\,0;p,\dots,p),
\]
and this is clearly orthogonal to the roots of $\frl\cap\frk$, which are given by
\eq
\label{roots lk D6}
\Delta^+(\frl\cap\frk)=\{\eps_i-\eps_j\,\big|\,1\leq i<j\leq p\}\cup\{\eps_{p+1+i}-\eps_{p+1+j}\,\big|\,1\leq i<j\leq q\}.
\eeq
So the constants $c=c^n_p$ can be calculated from \eqref{def const}. 

The set $\Delta(\frp_1)$ consisting of noncompact roots that are 1 on $h_p$ is 
\[
\Delta(\frp_1)= 
\{\eps_i+\eps_{p+1}\,\big|\, 1\leq i\leq p\}\cup\{-\eps_{p+1}-\eps_{p+1+j}\,\big|\, 1\leq j\leq q\}.
\]
So for any $C\subseteq\Delta(\frp_1)$,
\eq
\label{2rhoC D6}
2\rho(C)=(c_1,\dots,c_p\,|\,d;-e_1,\dots,-e_q),
\eeq
with 
\begin{align}
\label{ineq cde D6}
&0\leq c_i\leq 1;\qquad -q\leq d\leq p;\qquad 0\leq e_j\leq  1;\qquad \textstyle{d=\sum_i c_i-\sum_je_j.}
\end{align}
To compute the constant $c_p^n$ using \eqref{def const}, we take $\lambda=\lambda_0$, where
\begin{align}
\label{la0 D6}
&\lambda_0=(n,n-1,\dots,q+2\,|\,p+1;n-1,n-2,\dots,p+1),\qquad p,q>0;\\ \nonumber
&\lambda_0=(\,|\,1;n-1,n-2,\dots,1),\qquad p=0,q>0;\\ \nonumber
&\lambda_0=(n,n-1,\dots,2\,|\,p+1),\qquad p>0,q=0;\\ \nonumber
&\lambda_0=(\,|\,1),\qquad p=q=0.
\end{align}
Using \eqref{roots lk D6}, we see that $\lambda_0$ differs from $\rho_{\frl\cap\frk}$ by a weight orthogonal to all roots of $\frl\cap\frk$, and hence 
\eq
\label{Plk D6}
P_{L\cap K}(\lambda_0)=1.
\eeq
So to compute $c_p^n$ we have to compute the left side of \eqref{def const}. The following proposition describes the relevant $A$ and $C$, and the corresponding $\Lambda$.

\begin{proposition}
\label{prop D6} 
Let $\Lambda=\lambda_0-2\rho(A)-2\rho(C)$, with $\lambda_0$ given by \eqref{la0 D6}, and with $A\subseteq\Delta_n^+(\frl)$, $C\subseteq\Delta(\frp_1)$ as above. 

Suppose that for some $A$ and $C$ the corresponding $\Lambda$ satisfies $P_K(\Lambda)\neq 0$. Let $r=\frac{p}{2}$ and $s=\frac{q}{2}$. Then there is an $(r,s)$ shuffle
\[
1\leq i_1<\dots<i_r\leq r+s;\qquad 1\leq j_1<\dots<j_s\leq r+s
\]
of $1,2,\dots,r+s$ such that 
\[
A=\{\alpha_{u,v},\beta_{u,v}\,\big|\, 1\leq u\leq r,\,1\leq v\leq s\},
\]
where 
\[
\alpha_{u,v}=\eps_{p+1-u}+\eps_{n+1-v};\qquad \beta_{u,v}=\left\{
\begin{matrix} 
\eps_{p+1-u}+\eps_{p+1+v},&\qquad i_u<j_v;\cr
\eps_u+\eps_{n+1-v},&\qquad i_u>j_v,
\end{matrix}
\right.
\]
and
\[
C=\{\eps_i+\eps_{p+1}\,\big|\, r+1\leq i\leq p\}\cup\{-\eps_{p+1}-\eps_{p+1+j}\,\big|\, 1\leq j\leq s\}.
\]
The corresponding $\Lambda$ is
\[
\Lambda=(n+1-i_1,\dots,n+1-i_r,i_r,\dots,i_1\,|\, r+s+1;
n+1-j_1,\dots,n+1-j_s,j_s,\dots,j_1).
\]
If $p=0$, then the shuffle is necessarily trivial, i.e., there are no $i_u$ and $(j_1,\dots,j_s)=(1,\dots,s)$. This means that 
\begin{align*}
&A=\emptyset; \\
&C=\{-\eps_{p+1}-\eps_{p+1+j}\,\big|\, 1\leq j\leq s\};\\
&\Lambda=(\,|\, s+1; n,\dots,n+1-s,s,\dots,1)=(\,|\, s+1; 2s+1,\dots,s+2,s,\dots,1).
\end{align*}
Similarly, if $q=0$ than $(i_1,\dots,i_r)=(1,\dots,r)$, there are no $j_v$, and
\begin{align*}
&A=\emptyset; \\
&C=\{\eps_i+\eps_{p+1}\,\big|\, r+1\leq i\leq p\};\\
&\Lambda=(n,\dots,n+1-r,r,\dots,1\,|\,r+1)=(2r+1,\dots,r+2,r,\dots,1\,|\, r+1).
\end{align*}
Finally, if $p=q=0$, i.e., $n=1$, then the shuffle contains no $i_u$ or $j_v$, $A=C=\emptyset$, and $\Lamdba=\lambda_0=(\,|\,1)$.
\end{proposition}
\pf
The statement is obviously true if $n=1$, i.e., if $p=q=0$.
We proceed by induction on $n$. 

So let us assume that $n\geq 3$ is odd, and let $0\leq p\leq n-1$ be an even integer. We assume that the statement is true for $n-2$ and for any even integer $p'$ between $0$ and $n-3$.

Using the definitions and the inequalities 
\eqref{ineq ab D6}, we see that
\begin{align*}
\Lambda=(\unb{[p,n]}{n-a_1-c_1},\unb{[p-1,n-1]}{n-1-a_2-c_2},\dots,\unb{[1,q+2]}{q+2-a_p-c_p}\,|\,\unb{[1,n]}{p+1-d};\\\unb{[q,n]}{n-1-b_1+e_1},\unb{[q-1,n-1]}{n-2-b_2+e_2},\dots,
\unb{[1,p+2]}{p+1-b_q+e_q}).
\end{align*}
So the coordinates of $\Lambda$ are $n$ integers between $1$ and $n$, and assuming that $P_K(\Lamdba)\neq 0$, they have to be different from each other, i.e., $\Lambda$ has to be a permutation of $(n,\dots,1)$. In particular, some $\Lambda_i$ must be equal to $n$ and there are three possibilities:
\eq
\label{choice 1 D6}
\Lambda_1=n\qquad\text{or}\qquad\Lambda_{p+1}=n\qquad\text{or}\qquad\Lambda_{p+2}=n.
\eeq
Assume first that $\Lambda_1=n$; this is only possible if $p>0$, i.e., $p\geq 2$. Then
\eqn
a_1=0,\qquad c_1=0,
\eeqn
and it follows that
\begin{align*}
&\eps_1+\eps_{p+1+j}\notin A,\qquad 1\leq j\leq q;\\ 
&\eps_1+\eps_{p+1}\notin C.
\end{align*}
This implies that
\begin{align*}
&0\leq b_j\leq p-1,\qquad 1\leq j\leq q; \\ 
&-q\leq d\leq p-1,
\end{align*}
and so
\begin{align*}
\Lambda=(n,\unb{[p-1,n-1]}{n-1-a_2-c_2},\dots,\unb{[1,q+2]}{q+2-a_p-c_p}\,|\,\unb{[2,n]}{p+1-d};\unb{[q+1,n]}{n-1-b_1+e_1},\\
\unb{[q,n-1]}{n-2-b_2+e_2},\dots,\unb{[2,p+2]}{p+1-b_q+e_q}).
\end{align*}
We see that there is exactly one place where 1 can be, i.e., 
\[
\Lambda_p=1.
\]
This implies 
\eqn
a_p=q,\qquad c_p=1,
\eeqn
and therefore
\begin{align*}
&\eps_p+\eps_{p+1+j}\in A,\qquad 1\leq j\leq q;\\
&\eps_p+\eps_{p+1}\in C.
\end{align*}
It follows that
\begin{align*}
&1\leq b_j\leq p-1,\qquad 1\leq j\leq q,\\ \nonumber
&-q+1\leq d\leq p-1
\end{align*}
and so
\begin{align*}
\Lambda=(n,\unb{[p-1,n-1]}{n-1-a_2-c_2},\dots,\unb{[2,q+3]}{q+3-a_{p-1}-c_{p-1}},1\,|\,\unb{[2,n-1]}{p+1-d};\\
\unb{[q+1,n-1]}{n-1-b_1+e_1},\unb{[q,n-2]}{n-2-b_2+e_2},\dots,
\unb{[2,p+1]}{p+1-b_q+e_q}).
\end{align*}
Let now $\frg'\cong\frs\fro^*(2(n-2))$ be the subalgebra of $\frg$ built on coordinates $2,\dots,p-1,p+1,\dots,n$, and let $\frl'=\frl\cap\frg'$. 
We consider the real form of the corresponding $\caO^\bbC$ given by $h=h_{p-2}$.

Then
\begin{align*}
&\Delta_n^+(\frl')=\Delta_n^+(\frl)\setminus\{\eps_1+\eps_{p+1+j},\eps_p+\eps_{p+1+j}\,\big|\,1\leq j\leq q\};\\
&\Delta(\frp'_1)=\Delta(\frp_1)\setminus\{\eps_1+\eps_{p+1},\eps_p+\eps_{p+1}\},
\end{align*}
and we set
\begin{align*}
&A'=A\setminus\{\eps_p+\eps_{p+1+j}\,\big|\,1\leq j\leq q\};\\
&C'=C\setminus\{\eps_p+\eps_{p+1}\}.
\end{align*}
We define $\lambda_0$ as in \eqref{la0 D6}, but with $n$ replaced by $n-2$ and $p$ replaced by $p-2$. Then $\Lambda'$ corresponding to $A'$ and $C'$ can be obtained from $\Lamdba$ by deleting coordinates $\Lambda_1$ and $\Lambda_p$, and decreasing all the other coordinates by $1$. More precisely,
deleting the first and the $p$-th coordinate, we have
\begin{align*}
&2\rho(A')=(a_2,\dots,a_{p-1}\,|\,0;b_1-1,\dots,b_q-1)=(a_1',\dots,a'_{p-2}\,|\,0;b'_1,\dots,b'_q);\\
&2\rho(C')=(c_2,\dots,c_{p-1}\,|\,d-1;-e_1\dots,-e_q)=(c_1',\dots,c_{p-2}'\,|\,d';-e_1'\dots,-e_q');\\
&\Lambda'=(n-2-a_1'-c_1',\dots,q+2-a'_{p-2}-c'_{p-2}\,|\,p-1-d';\\
&\qquad\qquad\qquad\qquad\qquad\qquad\qquad\qquad n-3-b'_1+e'_1,\dots,p-1-b'_q+e'_q)=\\
&\qquad\qquad\qquad\qquad\qquad\qquad\qquad(\Lambda_2-1,\dots,\Lambda_{p-1}-1\,|\,\Lambda_{p+1}-1,\dots,\Lambda_n-1).
\end{align*}
We now see that $\Lambda$ is a permutation of $(n,\dots,1)$ if and only if $\Lamdba'$ is a permutation of $(n-2,\dots,1)$. By inductive assumption, this is equivalent to $A'$ and $\Lambda'$ being defined by a shuffle as in the statement of the proposition, and this clearly implies the same statement for $A$ and $\Lamdba$.

The second possibility in \eqref{choice 1 D6} is 
$\Lambda_{p+1}=n$, which implies $d=-q$ and hence
\begin{align*}
&-\eps_{p+1}-\eps_{p+1+j}\in C,\qquad 1\leq j\leq q;\\
&\quad \ \eps_i+\eps_{p+1}\notin C,\qquad 1\leq i\leq p.
\end{align*}
This implies
\begin{align*}
&c_i=0,\qquad 1\leq i\leq p;\\
&e_j=1,\qquad 1\leq j\leq q,
\end{align*}
and so
\begin{align*}
\Lambda=(\unb{[p+1,n]}{n-a_1},\unb{[p,n-1]}{n-1-a_2},\dots,\unb{[2,q+2]}{q+2-a_p}\,|\,n;\unb{[q+1,n]}{n-b_1},\unb{[q,n-1]}{n-1-b_2},\dots,
\unb{[2,p+2]}{p+2-b_q}).
\end{align*}
We see that there is no place where 1 could be, so this case is impossible if $P_K(\Lamdba)\neq 0$.

The third possibility in \eqref{choice 1 D6} is 
$\Lambda_{p+2}=n$; this is possible only if $q>0$, i.e., $q\geq 2$. It follows that
\eqn
b_1=0,\qquad e_1=1,
\eeqn
and so
\begin{align*}
&\eps_i+\eps_{p+2}\notin A,\qquad 1\leq i\leq p;\\ 
&-\eps_{p+1}-\eps_{p+2}\in C.
\end{align*}
This implies
\begin{align*}
0\leq a_i\leq q-1,\qquad 1\leq i\leq p;\\ 
-q\leq d\leq p-1,
\end{align*}
and hence
\begin{align*}
\Lambda=(\unb{[p+1,n]}{n-a_1-c_1},\unb{[p,n-1]}{n-1-a_2-c_2},\dots,\unb{[2,q+2]}{q+2-a_p-c_p}\,|\,\unb{[2,n]}{p+1-d};\\
n,\unb{[q-1,n-1]}{n-2-b_2+e_2},\dots,\unb{[1,p+2]}{p+1-b_q+e_q}).
\end{align*}
We see that there is exactly one place where 1 can be, i.e., 
\[
\Lambda_n=1.
\]
This implies 
\eqn
b_q=p,\qquad e_q=0
\eeqn
and therefore
\begin{align*}
&\eps_i+\eps_n\in A,\qquad 1\leq i\leq p; \\ 
&-\eps_{p+1}-\eps_n\notin C.
\end{align*}
It follows that
\begin{align*}
&1\leq a_i\leq q-1,\qquad 1\leq i\leq p;\\ 
&-q+1\leq d\leq p-1,
\end{align*}
and so
\begin{align*}
\Lambda=(\unb{[p+1,n-1]}{n-a_1-c_1}\unb{[p,n-2]}{n-1-a_2-c_2},\dots,\unb{[2,q+1]}{q+2-a_p-c_p}\,|\,\unb{[2,n-1]}{p+1-d};
n,\unb{[q-1,n-1]}{n-2-b_2+e_2},\\
\unb{[q-2,n-2]}{n-3-b_3+e_3},\dots,
\unb{[2,p+3]}{p+2-b_{q-1}+e_{q-1}},1).
\end{align*}
We now reason in the same way as in the first case, and conclude that the proposition follows from the inductive assumption for $n-2$ with $p$ staying the same and $q$ being replaced by $q-2$.
\epf

To finish the computation of the constant $c_p^n$, we first note that for every $A$ and $C$ described in Proposition \ref{prop D6}
\eq
\label{cardAC D6}
\# A=2rs;\qquad \# C=r+s.
\eeq
On the other hand, since $\Lambda$ is a permutation of $(n,\dots,1)$, $P_K(\Lambda)$ is equal to $\pm 1$. To compute the sign, we need to find the parity of the permutation bringing $\Lamdba$ to $(n,\dots,1)$. As in type C, we find this parity by counting the number of inversions in $\Lambda$ when compared with $(n,\dots,1)$.
We know from Proposition \ref{prop D6} that 
\eqn
\Lamdba=(n+1-i_1,\dots,n+1-i_r,i_r,\dots,i_1\,|\,r+s+1;n+1-j_1,\dots,n+1-j_s,j_s,\dots,j_1).
\eeqn
Clearly $i_r,\dots,i_1$ are in inversion with $n+1-j_1,\dots,n+1-j_s$; that is $rs$ inversions. Arguing as in type C, we get further $rs$ inversions from the groups 
\eqn
n+1-i_1,\dots,n+1-i_r \qquad\text{and}\qquad n+1-j_1,\dots,n+1-j_s,
\eeqn 
and
\eqn
i_r,\dots,i_1 \qquad\text{and}\qquad   j_s,\dots,j_1.
\eeqn
Finally, the coordinate $\Lambda_{p+1}=r+s+1$ is in inversion
with 
\[
i_r,\dots,i_1\qquad\text{ and }\qquad n+1-j_1,\dots,n+1-j_s.
\]
So the total number of inversions is $2rs+r+s$, and combined
with \eqref{cardAC D6} this implies that 
the nonzero contributions to the sum in \eqref{def const C}, which we know come from $A$, $C$ and $\Lambda$ as in Proposition \ref{prop D6}, are all equal to
\[
(-1)^{\#A+\#C}P_K(\Lambda)=1.
\]
Furthermore, the number $N$ of \eqref{def N} satisfies
\[
N \equiv \frac{p}{2}\quad\text{mod }2.
\]
Since the number of nonzero summands is by Proposition \ref{prop D6} equal to the number of $(r,s)$-shuffles of $r+s$, i.e., to $\binom{r+s}{r}$, we have proved:

\begin{theorem}
\label{thm D6}
Let $G_\bbR=SO^*(2n)$, for an odd $n\geq 1$, and let $p$, $0\leq p\leq n-1$, be an even integer. Let $r=\frac{p}{2}$ and let $s=\frac{n-1-p}{2}$. Then the constant $c_p^n$ for the real form of $\caO^\bbC$ corresponding to $p$ is
\[
c_p^n=(-1)^{\frac{p}{2}}\binom{r+s}{r}.
\]
\end{theorem}
\bigskip
For the convenience % DV ience
of the reader, below is a table that gives the value of the constant for each real form in every case.
\newpage

 \small{\begin{table}[ht]
 \addtolength{\tabcolsep}{-6pt}
%\caption{On Springer assumption} 
\centering 
\scalebox{0.81}{
\begin{tabular}{c c c c c} 
\hline\hline 
% & & & &\\
$\begin{matrix}{\rule{0pt}{2.6ex}G_{\bbR}}\\{K_{\bbR}}\end{matrix}$ && ${{\mathcal O}^{\mathbb C}}$ & {Real forms}& { Constants} \\ % [0.5ex]
% & & & &\\
% inserts table heading 
\hline\hline 
% & & & &\\
% DV $SU(p,q)$& $S(U(p)\times U(q))$ &%\tiny
% {$\lbrack 2^p,1^{q-p}\rbrack$} &%\tiny
% {$(\underbrace{1,1,\cdots,1}_{k},\underbrace{-1,-1,\cdots,-1}_{p-k},\underbrace{1,1,\cdots,1}_{p-k},\underbrace{0,0,\cdots,0}_{q-p},\underbrace{-1,-1,\cdots,-1}_{k})$}
$\begin{matrix}\rule{0pt}{2.6ex}SU(p,q)\\S(U(p)\times U(q))\end{matrix}$&{$q\geq p\geq 1$} &%\tiny 
{$\lbrack 2^p,1^{q-p}\rbrack$} &%\tiny
{$\begin{matrix}\rule{0pt}{2.6ex}(\underbrace{1,\cdots,1}_{k},\underbrace{-1,\cdots,-1}_{p-k},\\ \underbrace{1,\cdots,1}_{p-k},\underbrace{0,\cdots,0}_{q-p},\underbrace{-1,\cdots,-1}_{k})\end{matrix}$}
&%\tiny
{$(-1)^{k(p+q-k)}{{p}\choose{k}}$}\\[2ex] 
% &  & & & \\[2ex] 
%\tiny
%{$q\geq p\geq 1$}
&  & &  $k=0,1,\cdots,p$ &  \\% [1ex] 
% &  &  & & \\ 
\hline\hline
% & & & & \\
$\begin{matrix}SO_{e}(2p,2q+1)\\ %\;\;\;\; % DV S(O(2p)\times O(2q+1))$ &
SO(2p)\times SO(2q+1) \end{matrix}$ &\quad {$q\geq p\geq 1$}&
%\tiny
{\;\;\;\;\;$\lbrack 3, 2^{2p-2},1^{2(q-p+1)}\rbrack$} &%\tiny
{$\begin{matrix}\rule{0pt}{2.6ex}(2,\underbrace{1,\cdots,1}_{p-1},\\ \underbrace{1,\cdots,1}_{p-1},\underbrace{0,\cdots,0}_{q-p+1})\end{matrix}$}
&%\tiny
{$(-1)^{\lbrack(p/2)\rbrack+1}2^{2p-2}$}\\ \cline{4-5} %[12ex] 
%\tiny
% {$q\geq p\geq 1$}&  & & & \\ % [-4ex]
& &  &%\tiny % DV
{$\begin{matrix}\rule{0pt}{2.6ex}(2,\underbrace{1,\cdots,-1}_{p-1},\\ \underbrace{1,\cdots,1}_{p-1},\underbrace{0,\cdots,0}_{q-p+1})\end{matrix}$}& %\tiny
{$(-1)^{\lbrack(p/2)\rbrack}2^{2p-2}$}  \\ \cline{4-5} %[9ex] 
& &\tiny $\begin{matrix}\text{3rd real form}\\ \text{only if }q>p-1\end{matrix}$ &  %\tiny
{$\begin{matrix}\rule{0pt}{2.6ex}(\underbrace{1,\cdots,1}_{p-1},0,\\ 2,\underbrace{1,\cdots,1}_{p-1},\underbrace{0,\cdots,0}_{q-p})\end{matrix}$}& %\tiny
{$0$}   \\%[5ex] 
% & & &%\tiny
% {$\text{ third real form only if }q>p-1$}&  \\[2ex] 
% & & & &\\
\hline\hline
% & & & &\\
$\begin{matrix}\rule{0pt}{2.6ex}\Sp(2n,{\mathbb R})\\ U(n)\end{matrix}$ & $n\geq 1$ &%\tiny
{$\lbrack 2^n\rbrack$} &%\tiny
{$(\underbrace{1,\cdots,1}_{k},\underbrace{-1,\cdots,-1}_{n-k})$}
&%\tiny
{$(-1)^{\lbrack(k+1)/2\rbrack}{{r+s}\choose{r}}$}\\[-3ex] 
& & & & %\tiny
{($n$ odd) or ($n$ and $k$ even)}\\[2ex] 
 %\tiny
% {$n\geq 1$}
& & &   $k=0,1,\cdots,n$ & %\tiny
{$0$}\\[-1ex] 
& & & & %\tiny
{$n$ even and $k$ odd}\\[5ex] 
& & & & \\[-4ex] 
&  & & & %\tiny
{$r=\lbrack \frac{k}{2}\rbrack$ and $s=\lbrack \frac{n-k}{2}\rbrack$}
\\[1ex] 
% & & & &\\
\hline\hline
% & & & &\\
$\begin{matrix}SO_{e}(2p,2q)\\ % & % $S(O(2p)\times O(2q))$& %\tiny
SO(2p)\times SO(2q)\end{matrix}$& $q\geq p\geq 1$ &
{$\lbrack 3, 2^{2p-2},1^{2(q-p)+1}\rbrack$} &%\tiny
{$\begin{matrix}\rule{0pt}{2.6ex}(2,\underbrace{1,\cdots,1}_{p-1},\\ \underbrace{1,\cdots,1}_{p-1},\underbrace{0,\cdots,0}_{q-p+1})\end{matrix}$}
&%\tiny
{$(-1)^{\lbrack(p-1)/2\rbrack}2^{2p-2}$}\\ \cline{4-5}% [8ex] 
%\tiny
% {$q\geq p\geq 1$}
&  & &  %\tiny
{$\begin{matrix}\rule{0pt}{2.6ex}(2,\underbrace{1,\cdots,-1}_{p-1},\\ \underbrace{1,\cdots,1}_{p-1},\underbrace{0,\cdots,0}_{q-p+1})\end{matrix}$}& %\tiny
{$(-1)^{\lbrack(p-1)/2\rbrack}2^{2p-2}$} \\ \cline{4-5}% [8ex] 
&  & & %\tiny
{$\begin{matrix}\rule{0pt}{2.6ex}(\underbrace{1,\cdots,1}_{p-1},0,\\2,\underbrace{1,\cdots,1}_{p-1},\underbrace{0,\cdots,0}_{q-p})\end{matrix}$}& %\tiny
{$(-1)^{\lbrack p/2\rbrack+1}2^{2p-1}$} \\[-4ex] 
& & & & %\tiny
{if $q>p$}\\[3ex] 
&  &  & & %\tiny
{$(-1)^{\lbrack p/2\rbrack+1}2^{2p-2}$}  \\[0ex]
& & & & %\tiny
{if $q=p$}\\ \cline{4-5}%[8ex] 
&  & \tiny $\begin{matrix}\text{ fourth real form}\\ \text{only if }q=p\end{matrix}$ & %\tiny
{$\begin{matrix}\rule{0pt}{2.6ex}(\underbrace{1,\cdots,1}_{p-1},0,\\ 2,\underbrace{1,\cdots,1}_{p-1},\underbrace{0,\cdots,0}_{q-p})\end{matrix}$}& %\tiny
{$(-1)^{\lbrack p/2\rbrack+1}2^{2p-2}$}  \\[1.5ex] 
% &  &  &%\tiny
% {$\text{ fourth real form only if }q=p$}&   \\[2ex] 
% & & & & \\
\hline\hline
% & & & & \\
$\begin{matrix}\rule{0pt}{2.6ex}SO^*(2n)\\U(n)\end{matrix}$ &$n\geq 1$& %\tiny
{$\lbrack 2^n\rbrack$} &%\tiny
{$(\underbrace{1,\cdots,1}_{p},\underbrace{-1,\cdots,-1}_{n-p})$}
&%\tiny
{$(-1)^{p/2}{{n/2}\choose{p/2}}$}\\% [7ex] 
%\tiny
% {$n\geq 1$}
&  & %\tiny
{$n$ even}&  $0\leq p\leq n$ with $p$ even&   \\ \cline{3-5}% [8ex] 
& &%\tiny
{$\lbrack 2^{n-1},1^2\rbrack$} &%\tiny
{\rule{0pt}{2.6ex}$(\underbrace{1,\cdots,1}_{p},0,\underbrace{-1,\cdots,-1}_{n-1-p})$}
&%\tiny
{$(-1)^{p/2}{{(n-1)/2}\choose{p/2}}$}\\% [7ex] 
&  &%\tiny
{$n$ odd}&  $0\leq p\leq n-1$ with $p$ even &  \\% [1ex] 
\hline\hline
\end{tabular} }\label{tableSpringer}
\end{table}}

\end{document}